%% file: mhdSISC.tex
\documentclass[final]{siamart1116}

\pdfoutput=1

\usepackage{color}

\usepackage{soul}

\usepackage{amssymb,amsmath}

\usepackage{calc}
\usepackage[margin=1.in]{geometry}

\input texSISC/defs

\usepackage{tikz}
\input texSISC/trimFig.tex

\usepackage{algorithm} 
\usepackage{algpseudocode} 

\numberwithin{theorem}{section}

\newcommand{\TheTitle}{A high-order finite difference WENO scheme for ideal magnetohydrodynamics on curvilinear meshes} 

\title{{\TheTitle}}

\author{
  Andrew J.~Christlieb\thanks{Department of Computational Mathematics, Science and Engineering, 
      Department of Mathematics and Department of Electrical and Computer Engineering, 
      Michigan State University, East Lansing, MI (\email{christli@msu.edu}).
Research is supported in part by AFOSR grants FA9550-12-1-0343, FA9550-12-1-0455, FA9550-15-1-0282 and  NSF grant DMS-1418804.}
  \and
  Xiao Feng\thanks{The MathWorks, Inc., Natick, MA (\email{xiao.feng@mathworks.com}).}
  \and
  Yan Jiang\thanks{Department of Mathematics, Michigan State University, East Lansing, MI (\email{jiangyan@math.msu.edu}).}
  \and
  Qi Tang\thanks{Department of Mathematical Sciences, Rensselaer Polytechnic Institute, Troy, NY (\email{tangq3@rpi.edu}).
Research is supported by the Eliza Ricketts Postdoctoral Fellowship.
  }
}

\begin{document}

\maketitle

\begin{abstract}
A high-order finite difference numerical scheme is developed for the ideal magnetohydrodynamic equations based on  
{an alternative} flux formulation of the weighted essentially non-oscillatory (WENO) scheme. 
It computes a {high-order} numerical flux by a Taylor expansion in space, 
with the lowest-order term solved from a Riemann solver 
and the higher-order terms constructed from physical fluxes by limited central differences.
The scheme coupled with several Riemann solvers, including a Lax-Friedrichs solver and HLL-type solvers,
is developed on general curvilinear meshes in two dimensions and verified on a number of benchmark problems.
In particular, a HLLD solver on Cartesian meshes is extended to curvilinear meshes with proper modifications.
A numerical boundary condition for the perfect electrical conductor (PEC) boundary 
is derived for general geometry and verified through a bow shock flow.
Numerical results also confirm the advantages of using low dissipative Riemann solvers in the current framework.

\end{abstract}


\input texSISC/intro

\input texSISC/governing

\input texSISC/alternative

\input texSISC/ct

\input texSISC/pplimiter

\input texSISC/boundaryConditions

\input texSISC/numerical

\input texSISC/conclusions

\vskip .5\baselineskip
\noindent{\bf Acknowledgements.} We would like to thank Professor J.W.~Banks for valuable discussions and comments.

\appendix
\input texSISC/WENO

\input texSISC/riemannSolver

\bibliographystyle{siamplain}
\bibliography{mhd}

\end{document}

%% file: texSISC/defs.tex
\renewcommand{\url}[1]{}
\newcommand{\citeCount}[1]{}
\def\hf{\frac{1}{2}}

\newcommand{\norm}[1]{\left\lVert{#1}\right\rVert}
\newcommand{\Alfven}{{Alfv\'{e}n}}
\newcommand{\divg}{\nabla \cdot}
\DeclareMathOperator{\sign}{sign}

\newcommand{\abs}[1]{\left|{#1}\right|}
\newcommand{\LL}{\text{L}}
\newcommand{\RR}{\text{R}}
\newcommand{\MM}{\text{M}}  

\newcommand{\hll}{\text{HLL}}
\newcommand{\hL}[1]{{#1}_{\LL}}
\newcommand{\hR}[1]{{#1}_{\RR}}
\newcommand{\hA}[1]{{#1}_{\alpha}}
\newcommand{\hLstar}[1]{{#1}_{\LL}^{*}}
\newcommand{\hRstar}[1]{{#1}_{\RR}^{*}}
\newcommand{\hAstar}[1]{{#1}_{\alpha}^{*}}
\newcommand{\hLstartwo}[1]{{#1}_{\LL}^{**}}
\newcommand{\hRstartwo}[1]{{#1}_{\RR}^{**}}
\newcommand{\hAstartwo}[1]{{#1}_{\alpha}^{**}}
\newcommand{\ptot}{p_{\text{tot}}}

\newcommand{\dx}{\Delta x}
\newcommand{\dy}{\Delta y}

\newcommand{\bogus}[1]{{}}

\newenvironment{myIndent}%
 {\list{}{\leftmargin=0.1in\rightmargin=0.1in}\item[]}%
  {\endlist}

{\noindent\textbf{Procedure~}{#1}\begin{myIndent}\em}
{\end{myIndent}}

\newcommand{\av}{\mathbf{ a}}

\newcommand{\fv}{\mathbf{ f}}
\newcommand{\gv}{\mathbf{ g}}

\newcommand{\iv}{\mathbf{ i}}

\newcommand{\nv}{\mathbf{ n}}

\newcommand{\qv}{\mathbf{ q}}
\newcommand{\rv}{\mathbf{ r}}

\newcommand{\tv}{\mathbf{ t}}
\newcommand{\uv}{\mathbf{ u}}
\newcommand{\vv}{\mathbf{ v}}

\newcommand{\xv}{\mathbf{ x}}

\newcommand{\Av}{\mathbf{ A}}
\newcommand{\Bv}{\mathbf{ B}}

\newcommand{\Ev}{\mathbf{ E}}
\newcommand{\Fv}{\mathbf{ F}}
\newcommand{\Gv}{\mathbf{ G}}

\newcommand{\curl}{\nabla \times}
\newcommand{\calo}{{\cal O}}

\newcommand{\Ec}{{\mathcal E}}

\newcommand{\tableFont}{\footnotesize}

\newcommand{\num}[2]{#1e#2} 

\clearpage

\newcommand{\omegaBar}{\bar{\omega}}
\newcommand{\alphaBar}{\bar{\alpha}}
\newcommand{\betaBar}{\bar{\beta}}
\newcommand{\dBar}{\bar{d}}

\newcommand\tnv{{\mathbf{t}}}

\renewcommand\tv{\tnv}

%% file: texSISC/trimFig.tex
\newlength{\tfwidth}
\newlength{\tfheight}
\newlength{\tfxa}
\newlength{\tfxb}
\newlength{\tfya}
\newlength{\tfyb}
\newcommand{\trimFigWithBox}[6]{%
\setlength\fboxsep{0pt}%
\setlength\fboxrule{1.0pt}
\fbox{\includegraphics[width=#2, clip, trim=#3 #4 #5 #6]{#1}}%
}
\newcommand{\trimFigNoBox}[6]{%
\setlength\fboxsep{1pt}
\setlength\fboxrule{0.0pt}
\fbox{\includegraphics[width=#2, clip, trim=#3 #4 #5 #6]{#1}}%
}
\newcommand{\trimFigHeightWithBox}[6]{%
\setlength\fboxsep{0pt}%
\setlength\fboxrule{1.0pt}
\fbox{\includegraphics[height=#2, clip, trim=#3 #4 #5 #6]{#1}}%
}
\newcommand{\trimFigHeightNoBox}[6]{%
\setlength\fboxsep{1pt}
\setlength\fboxrule{0.0pt}
\fbox{\includegraphics[height=#2, clip, trim=#3 #4 #5 #6]{#1}}%
}

\newcommand{\trimFig}[6]{%
\setlength{\tfwidth}{(#2+#2*\real{#3})+#2*\real{#4}}
\setlength{\tfheight}{(#2+#2*\real{#5})+#2*\real{#6}}%
\setlength{\tfxa}{\tfwidth*\real{#3}}%
\setlength{\tfxb}{\tfwidth*\real{#4}}%
\setlength{\tfya}{\tfheight*\real{#5}}%
\setlength{\tfyb}{\tfheight*\real{#6}}%
\trimFigNoBox{#1}{#2}{\tfxa}{\tfya}{\tfxb}{\tfyb}%
}

\newsavebox\figBox

\newcommand{\trimw}[6]{%
\sbox\figBox{\includegraphics{#1}}
\setlength{\tfwidth}{\the\wd\figBox}
\setlength{\tfheight}{\the\ht\figBox}
\setlength{\tfxa}{\tfwidth*\real{#3}}%
\setlength{\tfxb}{\tfwidth*\real{#4}}%
\setlength{\tfya}{\tfheight*\real{#5}}%
\setlength{\tfyb}{\tfheight*\real{#6}}%
\trimFigNoBox{#1}{#2}{\tfxa}{\tfya}{\tfxb}{\tfyb}%
}

\newcommand{\trimwb}[6]{%

\sbox\figBox{\includegraphics{#1}}
\setlength{\tfwidth}{\the\wd\figBox}
\setlength{\tfheight}{\the\ht\figBox}
\setlength{\tfxa}{\tfwidth*\real{#3}}%
\setlength{\tfxb}{\tfwidth*\real{#4}}%
\setlength{\tfya}{\tfheight*\real{#5}}%
\setlength{\tfyb}{\tfheight*\real{#6}}%
\trimFigWithBox{#1}{#2}{\tfxa}{\tfya}{\tfxb}{\tfyb}%
}

\newcommand{\trimh}[6]{%
\sbox\figBox{\includegraphics{#1}}
\setlength{\tfwidth}{\the\wd\figBox}
\setlength{\tfheight}{\the\ht\figBox}
\setlength{\tfxa}{\tfwidth*\real{#3}}%
\setlength{\tfxb}{\tfwidth*\real{#4}}%
\setlength{\tfya}{\tfheight*\real{#5}}%
\setlength{\tfyb}{\tfheight*\real{#6}}%
\trimFigHeightNoBox{#1}{#2}{\tfxa}{\tfya}{\tfxb}{\tfyb}%
}

\newcommand{\trimhb}[6]{%

\sbox\figBox{\includegraphics{#1}}
\setlength{\tfwidth}{\the\wd\figBox}
\setlength{\tfheight}{\the\ht\figBox}
\setlength{\tfxa}{\tfwidth*\real{#3}}%
\setlength{\tfxb}{\tfwidth*\real{#4}}%
\setlength{\tfya}{\tfheight*\real{#5}}%
\setlength{\tfyb}{\tfheight*\real{#6}}%
\trimFigHeightWithBox{#1}{#2}{\tfxa}{\tfya}{\tfxb}{\tfyb}%
}

%% file: texSISC/intro.tex
\section{Introduction}

The ideal magnetohydrodynamic (MHD) equations are a fluid model to describe the dynamics of a perfectly conducting quasi-neutral plasma.
The equations are a system of nonlinear hyperbolic conservation laws with the constraint that its magnetic field is divergence free. 
In this work, we describe a high-order finite difference schemes for ideal MHD
based on an alternative flux formulation of the weighted essentially
non-oscillatory (WENO) scheme~\cite{Jiang2013, Jiang2014,Shu1988}.
The resultant scheme is applicable to general curvilinear meshes, {which can be obtained by a smooth or non-smooth mapping,} and compatible with many approximate Riemann solvers.

In recent years, high-order numerical schemes using essentially
non-oscillatory (ENO) and WENO approaches have been extended to the ideal MHD equations
in~\cite{balsara2009divergence,Balsara2013,Christlieb2016,Christlieb2014,IVAN2015157,Jiang1999, LondrilloDelZanna00,shen2012cusp,SUSANTO2013141} for example.
Many of those ENO/WENO approaches use the idea of \emph{reconstruction}, 
in which a numerical flux is typically reconstructed from the physical flux. 
We refer the reader to the review paper~\cite{Shu2009} for details. 
However, the results in~\cite{nonomura2010freestream,visbal2002use} show that when a 
standard finite difference WENO scheme is applied to curvilinear meshes, 
the free-stream preservation condition is not satisfied, 
which will cause large errors and even lead to numerical instabilities for high-order schemes. 
This issue can be resolved by an alternative 
flux formulation for the conservative finite difference WENO scheme in~\cite{Shu1988}. 
In this formulation, a WENO interpolation procedure is applied to the solution rather than to the flux functions. 
In~\cite{Jiang2014}, it has been theoretically proved and numerically demonstrated that 
this scheme can preserve free-stream solutions on both stationary and dynamically generalized coordinate systems, 
hence giving much better performance than the standard finite difference WENO schemes on curvilinear meshes. 
In addition, the alternative flux formulation takes advantage of monotone fluxes for the scalar case 
and  approximate Riemann solvers for the system case, 
while the standard finite difference WENO schemes 
can only use certain fluxes since its nonlinear stability relies on a smooth flux splitting.
Note that the most commonly used flux splitting in finite difference schemes is a Lax-Friedrichs flux splitting, 
which is one of the most diffusive Riemann solvers. 
Therefore, in this work we rely on the alternative flux formulation of the WENO schemes 
to solve the MHD equations on curvilinear meshes.

The alternative flux formulation requires an approximated Riemann solver in the low-order terms.
In this work, a Lax-Friedrichs Riemann solver and HLL-type Riemann solvers are used. 
The HLL Riemann solver, first proposed in~\cite{Harten1983}, 
solves a Riemann problem by an approximate solution  consisting of one intermediate state that is connected to the left and right states by discontinuities.
This intermediate state is obtained by exploiting the conservation of the equations, commonly referred to as the \emph{consistency condition}.  
When applied to the Euler equations in hydrodynamics, the HLL solver exhibits excessive dissipations in the presence of contact discontinuities.  
To remedy this, a HLLC (C stands for contact) solver is proposed  for the Euler equations in~\cite{Toro1994}.  
It assumes two intermediate states in the approximate solution, which are connected to each other by a contact discontinuity and connected to the left and right states by shocks.  
The \emph{Rankine-Hugoniot condition}, in addition to the consistency condition, is used to determine the intermediate states.
Similar ideas were later used in designing Riemann solvers for the ideal MHD equations in~\cite{Gurski2004,Li2005,Miyoshi2005}. 
The solvers in~\cite{Gurski2004, Li2005} were both named HLLC solvers for ideal MHD, 
because two intermediate states are assumed to be connected to each other by a contact discontinuity.  
The solver in~\cite{Miyoshi2005} was named HLLD (D stands for discontinuities) solver for ideal MHD, 
since the solver involving four intermediate states can exactly resolve most types of discontinuities, the only exception being the slow shocks.  
Note that the HLL and HLLC solvers can be directly applied to general curvilinear meshes while the HLLD solver is previously designed for Cartesian meshes.
{Besides the aforementioned one-dimensional Riemann solvers, we note that there are recent
developments of multidimensional Riemann solvers for ideal MHD,
see~\cite{balsara2014, balsara2015, balsara2017} for instance.}
In the current work, we experiment the HLL-type solvers when exploring the effects of 
the choice of Riemann solver in the current framework. In particular, the HLLD solver is extended to curvilinear meshes.

One of the main numerical difficulties for simulating the ideal MHD equations is to control divergence errors in the magnetic field.
Failure to control the divergence error creates an unphysical force parallel
to the magnetic field (see~\cite{Brackbill1980} for instance), 
which may eventually result into numerical instabilities as its effects accumulate.  
There are mainly four types of numerical approaches to address this issue, including (1) the non-conservative eight-wave
method~\cite{Gombosi1994}, (2) the projection method~\cite{Brackbill1980},
(3) the hyperbolic divergence cleaning method~\cite{Dedner2002}, and (4) the
various {constrained} transport methods~\cite{Balsara2004,Balsara1999a,Christlieb2014,Dai1998,Evans1988,Fey2003,Helzel2011,Helzel2013,Rossmanith2006}.
See the review paper~\cite{Toth2000} for more discussions on the advantages and disadvantages of those approaches.
In this work, we use a finite difference constrained transport method proposed in~\cite{Christlieb2014} to address this issue.


The main motivation for using curvilinear meshes in the current work is
that in certain MHD applications, the complex geometry is easier to describe using boundary-fitted grids in a curvilinear coordinate system, 
see~\cite{Chacon2004,Dhaeseleer1991, Golovin2011, IVAN2015157,Rosenbluth1979} for instance.
In a curvilinear coordinate, it is also relatively easier for finite difference methods to impose boundary conditions in the presence of curved surfaces,
compared to some other approaches such as cut-cell methods.
In this work, we impose several boundary conditions in the numerical tests for the ideal MHD equations, 
including the inflow and outflow boundary conditions, and a perfect electrical conductor (PEC) condition.
In particular, we derive a numerical compatibility boundary condition for the PEC boundary for the both conserved quantities and magnetic potential,
based on the previous work in the Euler equations~\cite{Henshaw2006moving} and Maxwell's equations~\cite{max2006b}.
Some discussions of the numerical boundary conditions for the ideal MHD equations can be found 
in~\cite{ivan2011thesis,SUSANTO2013141},
in which the PEC boundary is implemented differently through a least-squares reconstruction to set the normal components of the velocity
and magnetic field to zero.
Instead, we rely on the governing equations and a local characteristic analysis to derive the numerical boundary condition.

The remaining sections of the paper are organized as follows.
The governing equations are reviewed in Section~\ref{sec:IdealMHD}.
The details of the alternative flux formulation are presented in Section~\ref{sec:alternative}, 
including an outline of the base scheme, an additional limiter applied to higher order terms and its extensions to ideal MHD and curvilinear meshes.
Section~\ref{sec:numericalApproach} presents some numerical approaches applied to the ideal MHD simulations,
which include a brief outline of the constrained transport method and a positivity-preserving limiter, 
and a derivation of numerical boundary conditions for the PEC boundary.
Numerical results are presented in Section~\ref{sec:numericalResults}.
Conclusions and future directions are given in Section~\ref{sec:conclusions}. 
A WENO interpolation is described in Appendix~\ref{sec:WENO}.
Several HLL-type Riemann solvers are detailed in Appendix~\ref{sec:riemannSolver}, including our version of the HLLD solver.

%% file: texSISC/governing.tex
\section{Governing equations}
  \label{sec:IdealMHD}
  In this section we briefly review the ideal MHD equations, with an emphasis on the hyperbolicity and discontinuities of the system.
  In a conservation form, the ideal MHD equations are
  \begin{gather}
    \label{eq:MHD}
    {\partial_t}
    \begin{bmatrix}
      \rho \\
      \rho \uv \\
      \Ec \\
      \Bv
    \end{bmatrix}
    + \divg
    \begin{bmatrix}
      \rho \uv \\
      \rho \uv \otimes \uv + \ptot\mathbb{I} - \Bv \otimes \Bv \\
      \uv(\Ec + \ptot) - \Bv(\uv\cdot\Bv) \\
      \uv \otimes \Bv - \Bv \otimes \uv
    \end{bmatrix}
    = 0, \\
    \label{eq:DivergenceFree}
    \divg \Bv = 0,
  \end{gather}
  where $\rho$ is the mass density, $\rho \uv = (\rho u, \rho v, \rho w)^T$ is the momentum
  density, $\Ec$ is the total energy density, $\Bv  = (B_1, B_2, B_3)^T$ is the magnetic
  field, $p$ is the thermal pressure, $\left\|\cdot\right\|$ is the
  Euclidean vector norm, and $\ptot = p + \frac{1}{2}\lVert\Bv\rVert^2$
  is the total pressure.
  Let $\gamma = 5/3$ be the ideal gas constant, and
  the pressure satisfies the equation of state
  \begin{equation*}
    \Ec = \frac{p}{\gamma - 1} + \frac{\rho\left\| \uv \right\|^2}{2} + \frac{\left\|\Bv\right\|^2}{2}.
  \end{equation*}


  \subsection{Waves in the ideal MHD equations}

  The wave speeds of the ideal MHD system~\eqref{eq:MHD} {in some arbitrary direction $\nv$ ($\|\nv\|=1$)} are
  \begin{subequations}
  \begin{alignat}{3}
    \label{eq:WaveSpeedFast}
    \lambda_{1,8}&=\uv\cdot\nv \mp c_f  && \qquad \text{(fast magnetosonic waves),}\\
    \label{eq:WaveSpeedAlfven}            
    \lambda_{2,7}&=\uv\cdot\nv \mp c_a  && \qquad \text{(\Alfven~waves),}\\
    \label{eq:WaveSpeedSlow}              
    \lambda_{3,6}&=\uv\cdot\nv \mp c_s  && \qquad \text{(slow magnetosonic waves),}\\
    \label{eq:WaveSpeedEctropy}           
    \lambda_{4}=\lambda_{5}&=\uv\cdot\nv&& \qquad \text{(entropy and divergence waves),}
  \end{alignat}
  \end{subequations}
  where
  \begin{alignat*}{3}
    a &= \sqrt{\frac{\gamma p}{\rho}}  &&\qquad \text{(sound speed)}, \\
    c_a &= \sqrt{\frac{{(\Bv\cdot\nv)}^2}{\rho}} && \qquad \text{(\Alfven~speed)}, \\
    c_f &={\left[\frac{1}{2}{\left(a^2 + \frac{{\norm{\Bv}}^2}{\rho}  +  \sqrt{{\left(a^2+\frac{{\norm{\Bv}}^2}{\rho}\right)}^2-4a^2\frac{{\left(\Bv\cdot\nv\right)}^2}{\rho}}\right)}\right]}^{\frac{1}{2}}  && \qquad \text{(fast magnetosonic speed),} \\
    c_s &={\left[\frac{1}{2}{\left(a^2 + \frac{{\norm{\Bv}}^2}{\rho}  -  \sqrt{{\left(a^2+\frac{{\norm{\Bv}}^2}{\rho}\right)}^2-4a^2\frac{{\left(\Bv\cdot\nv\right)}^2}{\rho}}\right)}\right]}^{\frac{1}{2}}  && \qquad \text{(slow magnetosonic speed).}
  \end{alignat*}
  The eigen-decomposition of the Jacobian matrix for ideal MHD equations
  is complicated and has its own subtleties, see~\cite{Barth1999} for instance.
  More details of MHD waves can be found in many MHD literature, see~\cite{Jeffrey1964,Barth1999} for example.

  \subsection{Discontinuities in the ideal MHD equations}
  \label{sec:MHDDiscontinuities}

  The different types of discontinuities in the ideal MHD equations~(\ref{eq:MHD}) are reviewed in this section.
  Those discontinuities are used in the later discussion of the Riemann solvers in Appendix~\ref{sec:riemannSolver}.
  Let $\qv(t,\xv)$ denote the conserved quantities of the system,
  and assume the Riemann problem has an initial condition given by
  \begin{equation}
    \label{eq:RiemannProblemInitialCondition}
    \qv(0, \xv) =
    \begin{cases}
      \qv_{\LL}, &\text{if } \nv\cdot \xv  < 0, \\
      \qv_{\RR}, &\text{if } \nv\cdot \xv \geq 0,
    \end{cases}
  \end{equation}
  where $\qv_{\LL}$ and $\qv_{\RR}$ are constant vectors, and
  $\nv$ is an arbitrary direction.  
  The solution to such a problem is a
  function $\qv$ that depends only on $t$ and $\nv\cdot\xv$.
  We are interested in the case when the solution consists of a single
  moving discontinuity given by
  \begin{equation}
    \qv(t, \xv) =
    \begin{cases}
      \qv_{\LL}, &\text{if } (\nv\cdot\xv)/t < S, \\
      \qv_{\RR}, &\text{if } (\nv\cdot\xv)/t \geq S,
    \end{cases}
  \end{equation}
  where $S$ is the speed at which the discontinuity moves.
Let $\Fv$ denote the flux in the direction $\nv$.
  The \emph{Rankine-Hugoniot (RH) condition} of the hyperbolic conservation law is 
  \begin{equation}
    \label{eq:RHCondition}
    S(\qv_{\RR} - \qv_{\LL}) = \Fv(\qv_{\RR}) - \Fv(\qv_{\LL}).
  \end{equation}
Note that  the divergence condition~(\ref{eq:DivergenceFree}) in the ideal MHD equations implies that the magnetic field
  in the initial conditions~(\ref{eq:RiemannProblemInitialCondition}) satisfy
  \begin{equation}
    \label{eq:NormalMagneticFieldsEqual}
    \nv \cdot \Bv_{\LL} = \nv \cdot \Bv_{\RR}.
  \end{equation}
The RH condition~(\ref{eq:RHCondition}) and the constraint~(\ref{eq:NormalMagneticFieldsEqual}) imply that a
single moving discontinuity in ideal MHD must be one from the
following list:
\begin{enumerate}
\item A (fast or slow) \emph{shock}.  In this case, the solutions satisfy
    \begin{alignat*}{2}
    \nv\cdot\uv_{\alpha} &\neq S, && \qquad \alpha = \LL, \RR, \\
    \nv\cdot\uv_{\LL} &\neq \nv\cdot\uv_{\RR}, && \\
    \rho_{\LL} &\neq \rho_{\RR}. &&
  \end{alignat*}
\item A \emph{rotational discontinuity}.  In this case, the solutions satisfy
    \begin{alignat*}{2}
    \nv\cdot\uv_{\alpha} &\neq S, &&\qquad \alpha = \LL, \RR, \\
    \nv\cdot\uv_{\LL} &= \nv\cdot\uv_{\RR}, && \\
    \rho_{\LL} &= \rho_{\RR}, && \\
   \tv_i \cdot (\uv_{\RR} - \uv_{\LL}) &= \frac{1}{\sqrt{\rho}} \, \tv_i \cdot  (\Bv_{\RR} - \Bv_{\LL}), && \qquad i=1,2,
  \end{alignat*}
  where $\tv_1$ and $\tv_2$ are the tangential vectors with respect to 
  the discontinuity interface.
\item A \emph{contact discontinuity}.  In this case, the solutions satisfy
    \begin{alignat*}{2}
    \nv \cdot \Bv_\alpha & \neq 0,  && \qquad \alpha = \LL, \RR,\\
    \Bv_{\LL} & = \Bv_{\RR}, && \\
    \uv_{\LL} & = \uv_{\RR}, && \\
    p_{\LL} & = p_{\RR}. &&
  \end{alignat*}
\item A \emph{tangential discontinuity}.  In this case, the solutions satisfy
    \begin{alignat*}{2}
    \nv \cdot \Bv_\alpha & = 0,  && \qquad \alpha = \LL, \RR,\\
    \hL{\ptot}  &= \hR{\ptot}. &&
  \end{alignat*}
  Here the jumps in the tangential velocities and tangential magnetic fields can be arbitrary.
\end{enumerate}

    Note that  shocks are the only types of
    discontinuities that can possibly admit jumps in the normal
    velocities or the total pressures.
    Rotational, contact, and tangential discontinuities are
    \emph{linearly degenerate}.  Rotational discontinuities correspond
    to~\Alfven~waves, while contact and tangential discontinuities
    correspond to entropy and divergence waves.  
    None of the discontinuities in this list is genuinely nonlinear.
    The discussions on discontinuities in the ideal MHD equations can be also found in~\cite{Takahashi2013}.


The divergence-free condition suggests $\nv \cdot \Bv$ must be identical on 
both sides of the discontinuity.  However, this relation does not hold in multiple dimensions in numerical simulations.
Therefore, Riemann solvers need special treatments for such cases.
The treatment in our version of the HLLD solver will be discussed in Appendix~\ref{sec:riemannSolver}.

%% file: texSISC/alternative.tex
  \section{An alternative flux formulation of the WENO scheme}
  \label{sec:alternative}

In this section, we describe a WENO scheme based on an alternative flux formulation from~\cite{Jiang2013, Jiang2014,Shu1988}. 
The basic scheme is first given for a system of conservation laws. 
Numerical experiments indicate that its direct extension to the MHD equations causes oscillations in some benchmark problems. 
Hence, a limiter is introduced in Section~\ref{sec:AdditionalLimiter} to control those oscillations. 
The extensions of the base scheme to curvilinear coordinates are given in Sections~\ref{sec:Curvilinear}.

  \subsection{Basic scheme}
  \label{sec:1DSystem}
   
  A one-dimensional system of conservation law takes the form
  \begin{equation}
  \label{eq:1DSystem}
  {\partial_t \qv} + {\partial_x \fv(\qv)} = 0,
  \end{equation}
  where the conserved variables
  $\qv = (q_1(t, x), \ldots, q_n(t, x))^T$
  is a vector function of $t$ and $x$ and
  $\fv(\qv) = (f_1(\qv), \ldots, f_n(\qv))^T$
  is a flux function.  
  A hyperbolic system~(\ref{eq:1DSystem}) indicates the Jacobian $\partial \fv / \partial \qv$ has $n$
  real eigenvalues satisfying $\lambda_1(\qv) \leq \cdots \leq \lambda_n(\qv)$
  and a set of $n$ independent (right) eigenvectors, $\{\rv_1(\qv), \ldots, \rv_n(\qv)\}$.
  Defining a matrix 
  \begin{equation*}
  R(\qv) = (\rv_1(\qv), \ldots, \rv_n(\qv)),
  \end{equation*}
  we note that the Jacobian matrix satisfies
  \begin{equation*}
  R^{-1}(\qv) \, \dfrac{\partial \fv}{\partial \qv} \, R(\qv) = {\rm diag}( \lambda_1(\qv), \ldots, \lambda_n(\qv)).
  \end{equation*}
  
  The system~\eqref{eq:1DSystem} is solved by a semi-discrete conservative finite difference scheme of the form
  \begin{equation}
  \label{eq:1DSystemSemidiscrete}
  {\partial_t \qv_i} + \dfrac{1}{\Delta x} (\hat{\fv}_{i+1/2} - \hat{\fv}_{i-1/2}) = 0,
  \end{equation}
 on a uniform mesh with $x_i = i \Delta x$.
 {Here $\qv_i(t)$ sits on $x_i$, which is is the numerical approximation to the point value $\qv(x_i,t)$,}
  and $\hat{\fv}$ is some numerical flux that sits on half grid points. Note that $\hat{\fv}$ is a vector function and 
  the $k$-th component of the numerical flux satisfies
  \begin{equation}
  \label{eq:1DSystemNumericalorder}
  \dfrac{1}{\Delta x} \left( \hat{f}_{k}|_{i+1/2} - \hat{f}_{k}|_{i-1/2}  \right) = \partial_x f_{k}(\qv(x))|_{x=x_i} + \calo({\Delta x}^{m}),
  \end{equation}
  where $m$ is the spatial order of accuracy of the scheme.
  The semi-discrete form~(\ref{eq:1DSystemSemidiscrete}) is then integrated in time
  using a time-stepping method, such as Runge-Kutta (RK) methods.  In the current
  work, a third-order TVD-RK method is used.

The alternative flux formulation of the WENO scheme  $\hat{\fv}$ at $x_{i+1/2}$, first proposed in~\cite{Shu1988},  
is given by
{ 
   \begin{equation}
  \label{eq:1DNumericalFluxTaylorExpansion}
  \hat{\fv}_{i+1/2} = \fv_{i+1/2} + \sum_{k=1}^{[(m-1)/2]} a_{2k}{\Delta x}^{2k}\, {\partial_x^{2k} \fv}|_{i+1/2}
  \end{equation}
to guarantees $m$-th order accuracy in~\eqref{eq:1DSystemNumericalorder},}
  where $a_{2k}$'s are some constants obtained by Taylor expansions and the accuracy constraint.  
  In the current work, a truncation at $m = 5$ is used and $\hat{\fv}$ is therefore approximated by
  \begin{equation}
  \label{eq:1DNumericalFluxTruncated}
  \hat{\fv}_{i+1/2} = \fv_{i+1/2} - \dfrac{1}{24} {\Delta x}^2 \, \partial_x^2 \fv |_{i+1/2} + \dfrac{7}{5760} {\Delta x}^4 \, \partial_x^4 \fv |_{i+1/2}.
   \end{equation}
   The first term in~\eqref{eq:1DNumericalFluxTruncated}
   is approximated by
   \begin{equation}
   \fv_{i+1/2} = \Fv (\qv_{i+1/2}^{-}, \qv_{i+1/2}^{+}),
   \end{equation}
   where $\Fv$ is a Riemann solver and  $\qv_{i+1/2}^{\pm}$ are sufficiently high-order one-sided approximations to $\qv$ at
   $x_{i+1/2}$.  The WENO interpolation is used to obtain
   $\qv_{i+1/2}^{\pm}$, 
   and the formulation of a fifth-order WENO interpolation is given in Appendix~\ref{sec:WENO}. 
Approximate Riemann solvers are used in the current work, 
including a Lax-Friedrich solver and HLL-type Riemann solvers.
For instance, a Lax-Friedrichs solver gives
\begin{equation}
    \label{eq:LFFlux}
    \Fv(\qv^{-},\qv^{+})=\frac{1}{2}\left[\fv(\qv^{-}) + \fv(\qv^{+}) -\alpha \left( \qv^{+} - \qv^{-}\right) \right],
\end{equation} 
with $\alpha=\max_{1\leq k\leq n}|\lambda_{k}(\qv)|$ taken over the relevant range of $\qv$. 
Depending on the region where the maximum is taken, there are two variations of the solver, 
the global and local Lax-Friedrichs solvers, both of which are used in the numerical section.
Note that its local version is commonly refered to as the Rusanov flux.
The HLL-type solvers are more complicated and the details are described in Appendix~\ref{sec:riemannSolver}.
The HLL, HLLC and HLLD solvers are all experimented in this framework, 
and unsurprisingly the HLL and HLLC solvers are found to produce solutions better than the Lax-Friedrichs solver (the most dissipative one we test)
and worse than the HLLD solver (the least dissipative one we test) for most benchmark problems.
Therefore, to save the space, only the results of the Lax-Friedrichs and HLLD solvers are presented in the numerical section.

The remaining higher order terms in~\eqref{eq:1DNumericalFluxTruncated} are constructed using 
the physical flux $\fv_{i}$ at the grid points $x_i$. 
For instance, if some central differences are used, the approximations become
    \begin{subequations}
    \label{eq:1DNumericalFluxHigherOrderCentralDifferences}
\begin{alignat}{2}
     {\Delta x}^2 \partial_x^2 \fv |_{i+1/2} & \approx \dfrac{1}{48} \left( -5 \fv_{i-2} + 39 \fv_{i-1} - 34 \fv_{i} - 34 \fv_{i+1} + 39 \fv_{i+2} - 5 \fv_{i+3}  \right), \\
     {\Delta x}^4 \partial_x^4 \fv |_{i+1/2} & \approx \dfrac{1}{2}
     \left( \fv_{i-2} - 3 \fv_{i-1} + 2 \fv_{i} + 2 \fv_{i+1} - 3 \fv_{i+2} +
     \fv_{i+3} \right).
    \end{alignat}
    \end{subequations}
{ Both approximations in  \eqref{eq:1DNumericalFluxHigherOrderCentralDifferences} give a truncation error of $\calo(\Delta x^6)$, which guarantees a fifth order accuracy of the numerical flux \eqref{eq:1DNumericalFluxTruncated}.} 
More discussions on those terms are given in Section~\ref{sec:AdditionalLimiter}.
    
     The extension to multiple dimensions  can be treated in a dimension-by-dimension fashion.
     For example, a system of hyperbolic conservation law in two
     dimensions takes the form
     \begin{equation}
     \label{eq:2DScalarEquation}
     {\partial_t \qv} + {\partial_x \fv(\qv)} + {\partial_y \gv(\qv)}  = 0,
     \end{equation}
     where $\qv$ is a vector function of $t$, $x$ and $y$, and
     $\fv$ and $\gv$ are the fluxes in the $x$ and $y$ directions,
     respectively.  On a uniform mesh with $x_i = i\Delta x$ and
     $y_j = j\Delta y$, the system~(\ref{eq:2DScalarEquation}) can be
     solved by using a semi-discrete scheme,
     \begin{equation}
     \label{eq:2DScalarEquationSemidiscrete}
     {\partial_t \qv_\iv} + \dfrac{1}{\Delta x} (\hat{\fv}_{i+1/2, j} - \hat{\fv}_{i-1/2, j}) + \dfrac{1}{\Delta y} (\hat{\gv}_{i, j+1/2} - \hat{\gv}_{i, j-1/2}) = 0,
     \end{equation}
     where $\iv = (i,j)$ is a multi-index, and $\hat \fv$ and $\hat \gv$ are numerical fluxes approximated similarly as in 1D case,
     \begin{subequations}
         \begin{alignat}{2}
      & \hat{\fv}_{i+1/2,j} = \Fv(\qv^{-}_{i+1/2,j},\qv^{+}_{i+1/2,j}) - \dfrac{1}{24} {\Delta x}^2 \partial_x^2 \fv |_{i+1/2,j} + \dfrac{7}{5760} {\Delta x}^4 \partial_x^4 \fv |_{i+1/2,j},\\
      & \hat{\gv}_{i+1/2} = \Gv(\qv^{-}_{i,j+1/2},\qv^{+}_{i,j+1/2}) - \dfrac{1}{24} {\Delta y}^2 \partial_y^2 \gv |_{i,j+1/2} + \dfrac{7}{5760} {\Delta y}^4 \partial_y^4 \gv |_{i,j+1/2}.
      \end{alignat}
     \end{subequations} 
 Here $\qv^{\pm}_{i+1/2,j}$ and $\qv^{\pm}_{i,j+1/2}$ are obtained by one-dimensional WENO interpolations,
 and $\Fv$ and $\Gv$ are Riemann solvers corresponding to $\fv$ and $\gv$.

\subsection{Limiting the higher-order terms}
\label{sec:AdditionalLimiter}

In~\cite{Jiang2013} the higher-order derivatives
$\partial_x^2 \fv |_{i+1/2}$ and $\partial_x^4 \fv |_{i+1/2}$ are
first expanded in terms of the derivatives of $\fv$ with respect to
$\qv$ and spatial derivatives of $\qv$, and central differences are then used for approximations.
In~\cite{Jiang2014} the central differences given in~(\ref{eq:1DNumericalFluxHigherOrderCentralDifferences}) are used to approximate the higher-order terms directly.  
In~\cite{Jiang2013, Jiang2014}, good performance of the resulting schemes is demonstrated through benchmark problems of the compressible Euler equations in hydrodynamics.

  
  However, during the numerical experiment of the MHD equations, the linear approximations~\eqref{eq:1DNumericalFluxHigherOrderCentralDifferences} are found to  
  cause oscillations near a strong discontinuity, which could result into instabilities. 
  Hence, an extra limiting procedure is needed to switch the high-order numerical flux~\eqref{eq:1DNumericalFluxTaylorExpansion} to a first-order flux for such a case, 
  while the resulting scheme is required to retain high-order accuracy in smooth regions.
  In this work, the high order derivatives in~\eqref{eq:1DNumericalFluxTaylorExpansion} are multiplied by an additional limiter $\sigma$ and the numerical flux {(denoted by $\hat{\fv}^{\sigma}$)} becomes
{ 
 \begin{equation}
\label{eq:1DNumericalFluxTaylorExpansionModified}
\hat{\fv}^{\sigma}_{i+1/2} = \fv_{i+1/2} + \sigma_{i+1/2}\sum_{k=1}^{[(m-1)/2]} a_{2k}{\Delta x}^{2k}\, {\partial_x^{2k} \fv}|_{i+1/2}.
\end{equation}}
Note that for the case of $m = 5$ used in the current work, $\sigma$ needs to satisfy
\begin{align}
\label{eq:LimiterOrder}
\sigma_{i+1/2} = \left\{ \begin{array}{ll}
    1 + \calo({\Delta x}^3), & \text{when $\qv$ is smooth in the stencil $ S_{i+1/2} = \{ x_{i-2}, \ldots, x_{i+3} \}$}, \\
  \calo({\Delta x}^2), & \text{when $\qv$ contains a
    strong discontinuity in $S_{i+1/2}$.}\\
\end{array}
\right.
\end{align}
{
To see the effect of such a limiter, we note that for smooth problems,
$$\hat{\fv}_{i+1/2}-\hat{\fv}_{i+1/2}^{\sigma}=\left(1-\sigma_{i+1/2}\right)\sum_{k=1}^{[(m-1)/2]} a_{2k}{\Delta x}^{2k}\, {\partial_x^{2k} \fv}|_{i+1/2}=C_{i+1/2} \Delta x^5.$$
Moreover, the high order derivatives ${\partial_x^{2k} \fv}$ and the coefficient of the $\calo(\dx^3)$ term in $\sigma$ lead to a smooth coefficient of $C_{i+1/2}$, meaning $|C_{i+1/2}-C_{i-1/2}|=\calo(\Delta x)$. 
Hence, a difference of~\eqref{eq:1DNumericalFluxTaylorExpansionModified} gives
\begin{align*}
\frac{\hat{\fv}^{\sigma}_{i+1/2}-\hat{\fv}^{\sigma}_{i-1/2}}{\Delta x} =&
\frac{\hat{\fv}_{i+1/2}-\hat{\fv}_{i-1/2}}{\Delta x} - \frac{C_{i+1/2} \Delta x^5 - C_{i-1/2} \Delta x^5}{\Delta x} \\
=& \frac{\hat{\fv}_{i+1/2}-\hat{\fv}_{i-1/2}}{\Delta x} + \calo(\Delta x^5),
\end{align*}
and it appears to be enough to maintain the fifth-order accuracy in \eqref{eq:1DSystemNumericalorder}. In the following, we use $\hat{\fv}_{i+1/2}$ to represent $\hat{\fv}^{\sigma}_{i+1/2}$ without special declaration.}

Here we employ the parameter introduced in~\cite{christlieb2017kernel} to control the oscillations. 
The parameter, based on the idea of the WENO-Z scheme~\cite{castro2011high},
is constructed from the smoothness indicators $\beta_k$ in~\eqref{eq:beta} from the WENO interpolation.  
In the process to obtain $\qv^{-}_{i+1/2}$, we set
  \begin{equation*}
      \sigma_{\max}  = 1 + \dfrac{|\beta_0 - \beta_2|}{\epsilon + \min\{\beta_0, \beta_2\}}, \qquad
      \sigma_{\min}  = 1 + \dfrac{|\beta_0 - \beta_2|}{\epsilon +
        \max\{\beta_0,  \beta_2\}},
  \end{equation*}
where $\epsilon$ is a small positive number (taken to be $10^{-6}$
in all the numerical examples) to avoid division by zero.  We can thus obtain
a candidate for the coefficient $\sigma$ by
\begin{equation*}
  \sigma^{-} = \dfrac{\sigma_{\min}}{\sigma_{\max}}.
\end{equation*}
A similar formula for $\qv^{+}_{i+1/2}$ gives rise to another
candidate $\sigma^{+}$.  We finally set
\begin{equation}
\label{eq:Limiter}
  \sigma_{i+1/2} = \min\{ \sigma^{-}, \sigma^{+} \}.
\end{equation}
Using Taylor expansions, it is easy to verify the definition~\eqref{eq:Limiter} satisfies the constraint~\eqref{eq:LimiterOrder}. 
More details can be found in~\cite{castro2011high, christlieb2017kernel}.

\subsection{Curvilinear coordinates}
\label{sec:Curvilinear}

Here {we provide a brief discussion of using 
curvilinear coordinates as the computational domain to solve a general hyperbolic system.}
Assume the coordinates $\xv = (x, y)$ is related to the {curvilinear coordinates} $\rv = (\xi, \eta)$ via a continuous coordinate
  transformation $\xv  = \xv(\rv)$.  
As illustrated in Figure~\ref{fig:CurvilinearMeshDomain}, a uniform mesh in the computational domain is typically used in our implementation.

{
\newcommand{\figWidth}{10cm}
\newcommand{\trimfig}[2]{\trimFig{#1}{#2}{.0}{.0}{.0}{.05}}
\begin{figure}[htb]
\begin{center}
\begin{tikzpicture}[scale=1]
  \useasboundingbox (0.0,0) rectangle (12,3.9);  
  \draw(1,-.5) node[anchor=south west,xshift=0pt,yshift=+0pt] {\trimfig{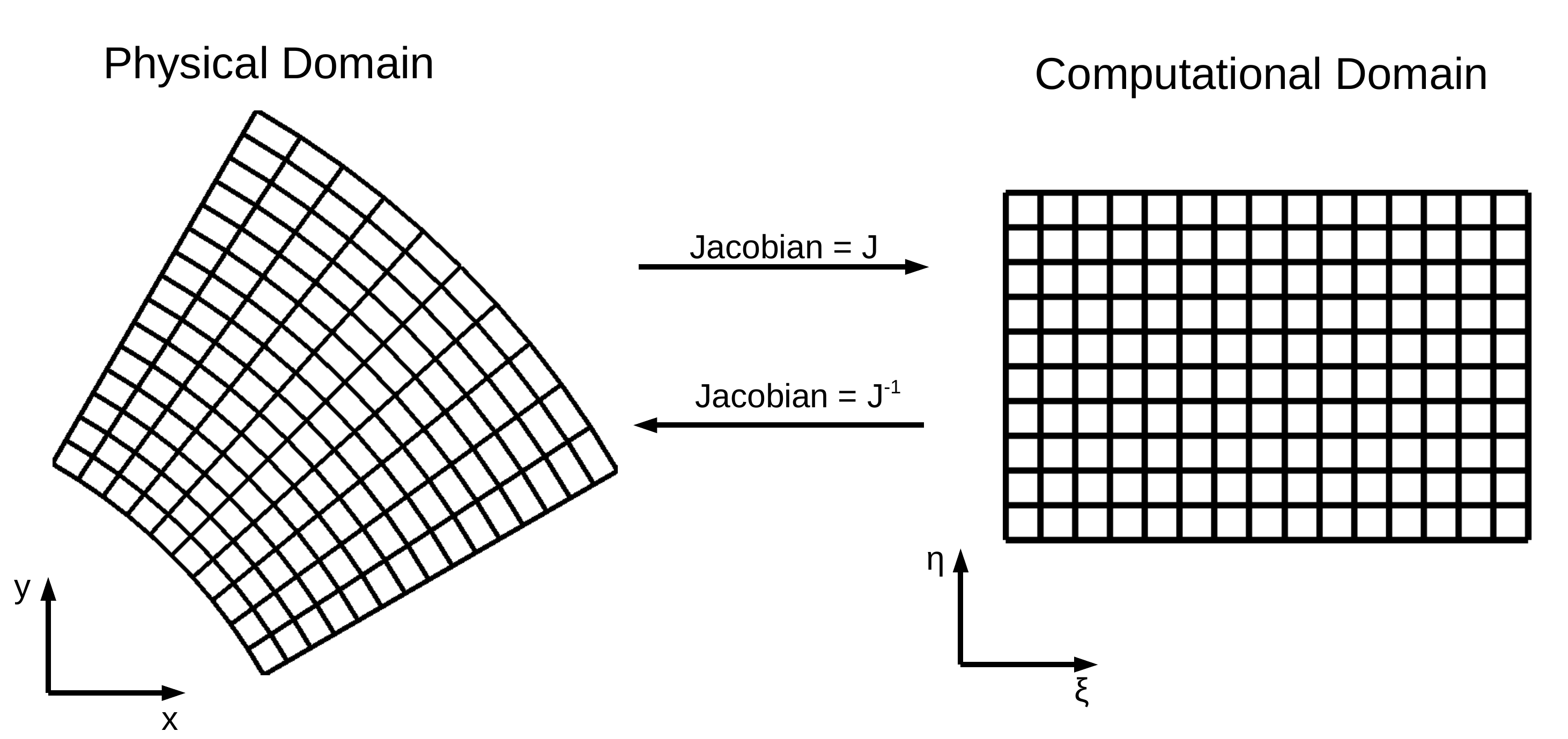}{\figWidth}};
\end{tikzpicture}
\end{center}
    \caption{
A schematic diagram of the transformations between the physical and computational domain.
}
\label{fig:CurvilinearMeshDomain}
\end{figure}
}

The two-dimensional system~\eqref{eq:2DScalarEquation} in the curvilinear coordinates has a conservative form given by 
  \begin{equation}
    \label{eq:2DScalarEquationXiEta}
    {\partial_t \widetilde{\qv}} + {\partial_\xi \widetilde{\fv}} + {\partial_\eta\widetilde{\gv}}  = 0,
  \end{equation}
  where
  \begin{equation}
      \label{eq:2DclFlux}
    \widetilde{\qv}  = \frac{\qv}{J}, \qquad
    \widetilde{\fv}  =\frac 1 J( \partial_x\xi \,\fv  + \partial_y \xi \, \gv), \qquad 
    \widetilde{\gv}  =\frac 1 J (\partial_x {\eta} \, \fv + \partial_y {\eta}\, \gv),
  \end{equation}
  and $J$ is the determinant of the Jacobian matrix defined by
  \begin{equation*}
      J := \text{det}\left[\frac{\partial \rv}{\partial \xv}\right],
  \end{equation*}
  which indicates $J^{-1}=\partial_\xi  x \,  \partial_\eta y -  \partial_\eta x \, \partial_\xi y$.
  The standard metrics satisfy
  \begin{alignat*}{3}
        & \frac{\partial_x \xi}  J  = \partial_\eta y, && \qquad 
         \frac{\partial_y \xi} J   = -\partial_\eta x, \\
        & \frac{\partial_x \eta} J  = -\partial_\xi y, && \qquad 
         \frac{\partial_y \eta} J  = \partial_\xi x.
  \end{alignat*}
  Note that the equation~(\ref{eq:2DScalarEquationXiEta}) in curvilinear coordinates is still hyperbolic.
  Therefore, the base scheme and other numerical treatments discussed previously can be applied straightforwardly on the uniform divided computational domain $(\xi_{i},\eta_{j})$
  after the numerical fluxes, $\tilde \fv$ and $\tilde \gv$, are defined properly through~\eqref{eq:2DclFlux}. 
  It has been suggested in~\cite{Jiang2014} that to preserve the freestream condition, 
   WENO interpolations should be applied to $\qv(\rv)$ instead of $\tilde{\qv}$.
   In the current work, we adopt this approach in the ideal MHD.
  In addition,  the metrics $\partial_\xi x$, $\partial_\xi y$, $\partial_\eta x$ and $\partial_\eta y$ at half point are approximated by central differences given by
  \begin{align*}
  w |_{i+1/2,j} = \frac{1}{256}  \left( 3w_{i-2,j} -25 w_{i-1,j} +150 w_{i,j} +150 w_{i+1,j} -25 w_{i+2,j} +3 w_{i+3,j} \right), \\
  w |_{i,j+1/2} = \frac{1}{256}  \left( 3w_{i,j-2} -25 w_{i,j-1} +150 w_{i,j} +150 w_{i,j+1} -25 w_{i,j+2} +3 w_{i,j+3} \right),
  \end{align*} 
  with $w$ stands the metrics, 
   {and have truncation errors $\calo(\Delta \xi^6)$ and $\calo(\Delta \eta^6)$, respectively.}
  We note that the HLLD flux requires some non-trivial extensions on curvilinear meshes 
  and the details are given in Appendix~\ref{sec:riemannSolver}.

%% file: texSISC/ct.tex
\section{Numerical approach in the ideal MHD} 
\label{sec:numericalApproach}
\subsection{Constrained transport}
  \label{sec:MHDConstrainedTransport}
A constrained transport framework is used to control the divergence error of the magnetic field.  
In this framework, alongside evolving the conserved quantities of the ideal MHD equations, 
a magnetic vector potential~$\Av$, satisfying $\Bv = \curl \Av$,
is evolved by
\begin{equation}
    \label{eq:AEquation3D}
  \partial_t \Av + ( \curl \Av  )\times \uv = 0.
\end{equation}
This evolution equation is derived from the magnetic induction equation, see~\cite{Rossmanith2006} for details.
In the case of two dimensions considered in this work,  the divergence-free condition~\eqref{eq:DivergenceFree} becomes
\begin{equation*}
  \divg \Bv = \partial_x B_1 + \partial_y B_2 = 0.
\end{equation*}
It therefore suffices to only account for $B_1$ and $B_2$ in terms of controlling divergence errors in two dimensions.
This leads to a nice property that only the third component of $\Av$ needs to be evolved.
For ease of presentation, a scalar quantity, $A$, is used to denote the third component of $\Av$.
In 2D, the equation~\eqref{eq:AEquation3D} leads to an evolution equation for $A$ given by
\begin{equation}
    \label{eq:AEquation2D}
  \partial_t A + u\, \partial_x A + v\, \partial_y A = 0
\end{equation}
and $\Bv = \curl \Av$ relates $B_1$ and $B_2$ with $A$ by
\begin{equation*}
  B_1 = \partial_y A, \qquad B_2 = - \partial_x A.
\end{equation*}
Many previous works show such a procedure can control the divergence error in $\Bv$ and improve numerical stabilities of base schemes.

\def\PhiWENO{\phi_\text{WENO}}
Same as our previous work~\cite{Christlieb2014}, a WENO method designed for Hamilton-Jacobi equations~\cite{jiang2000} is modified to solve the potential equation~(\ref{eq:AEquation2D}).
The approximation form is given by
    \begin{equation}
        \label{eq:discretizedAEquation}
        \partial_t  A_\iv = -u_\iv\left(\frac{\partial_x A^{-}_\iv + \partial_x A^{+}_\iv}{2}\right)
                             -v_{\iv}\left(\frac{\partial_y A^{-}_\iv + \partial_y A^{+}_\iv}{2}\right)
                         +\alpha_1 \left(\frac{\partial_x A^{+}_{\iv} - \partial_x A^{-}_{\iv}}{2}\right)
                             +\alpha_2 \left(\frac{\partial_y A^{+}_{\iv} - \partial_y A^{-}_{\iv}}{2}\right),
    \end{equation}
where $\alpha_1 = \max_{\iv} \lvert u_{\iv} \rvert$ and $\alpha_2 = \max_{\iv} \lvert v_{\iv} \rvert$.
Here $\partial_{x_m}  A^{\pm}_{ \iv}$ are defined by WENO reconstructions through
    \begin{equation*}
        \begin{aligned}
            \partial_{x_m} A^{-}_{\iv}
        &:=
        \PhiWENO\left(
            {D_{+x_m}A_{i-3, j}},
            {D_{+x_m}A_{i-2, j}},
            {D_{+x_m}A_{i-1, j}},
            {D_{+x_m}A_{i  , j}},
            {D_{+x_m}A_{i+1, j}}
            \right), \\
            \partial_{x_m}  A^{+}_{\iv}
        &:=
        \PhiWENO\left(
            {D_{+x_m}A_{i+2, j}},
            {D_{+x_m}A_{i+1, j}},
            {D_{+x_m}A_{i,   j}},
            {D_{+x_m}A_{i-1, j}},
            {D_{+x_m}A_{i-2, j}}
            \right), 
        \end{aligned}
    \end{equation*}
    and $D_{+x_m}$ is the standard forward difference defined by $D_{+x}A_{\iv} := (A_{i+1, j} - A_{i,j})/\Delta x$ and $D_{+y}A_{\iv} := (A_{i, j+1} - A_{i,j})/\Delta y$.
The function $\PhiWENO$ is the classical fifth-order WENO reconstruction whose coefficients can be found in many previous works such as~\cite{Christlieb2014,jiang2000,Jiang1996}.
In the current multistage setting, the constrained transport is implemented through a predictor-corrector strategy, i.e.,
after the $k$-stage of the time integrator at time $t^{n+1}$, the magnetic field is corrected by
\begin{align*}
    B^{n+1,(k)}_{1,\iv} & = \frac{A^{n+1,(k)}_{i,j-2} - 8 A^{n+1,(k)}_{i,j-1} + 8 A^{n+1,(k)}_{i,j+1} - A^{n+1,(k)}_{i,j+2}}{12\dy}, \\
    B^{n+1,(k)}_{2,\iv} & =-\frac{A^{n+1,(k)}_{i-2,j} - 8 A^{n+1,(k)}_{i-1,j} + 8 A^{n+1,(k)}_{i+1,j} - A^{n+1,(k)}_{i+2,j}}{12\dx}.
\end{align*}
Since the fourth-order central differences are used in the constrained transport step, the resulting scheme in this work
is fourth-order.
More details on the constrained transport using this approach can be found in the previous work~\cite{Christlieb2016, Christlieb2014}.


The potential $A$ on  curvilinear meshes are solved similarly.
For instance, the evolution equation~(\ref{eq:AEquation2D}) in the curvilinear coordinates becomes
\begin{equation*}
  \partial_t A + (u\, \partial_x {\xi} + v\, \partial_y{\xi} ) \partial_{\xi} A + (u\, \partial_x{\eta} + v\, \partial_y{\eta} ) \partial_{\eta} A = 0,
\end{equation*}
which is solved by a discretization similar to~\eqref{eq:discretizedAEquation} but approximated in the curvilinear coordinates.
We note that while such discretization only guarantees the divergence-free condition of
magnetic field to truncation errors on general curvilinear meshes, in practice we find this is
sufficient to suppress the unphysical oscillations associated with the divergence error of $\Bv$.


%% file: texSISC/pplimiter.tex
\subsection{A positivity-preserving limiter}
\label{sec:MHDPositivityPreserving}

During numerical experiments, we find that some HLL-type fluxes may cause
the density or pressure becoming negative in some hard problems (such as the cloud-shock interaction in Section~\ref{sec:cloudShock}) even with the constrained transport step turned on.
It appears to be related to the enhanced resolutions provided by the HLL-type solvers,
since numerical solutions remain positive if the Lax-Friedrichs flux is used in the scheme.
For such a case, a positivity-preserving limiter for ideal MHD equation in~\cite{Christlieb2015} is applied.

The limiter replaces the high-order flux $\hat \fv_{i+\hf}$ constructed in~\eqref{eq:1DNumericalFluxTaylorExpansionModified} 
with
a corrected flux $\hat{\fv}^\text{new}_{i+\hf}$ given by
\begin{align*}
    \hat{\fv}^\text{new}_{i+\hf} =  \theta_{i+\hf} \hat{\fv}_{i+\hf} +(1- \theta_{i+\hf})\hat{\fv}^\text{low}_{i+\hf},
\end{align*}
where $\hat \fv^\text{low}_{i+\hf}$ is a low-order flux {which could preserve positive density and pressure,} and the parameter $\theta_{i+\hf}\in[0,1]$.
The parameter $\theta_{i+\hf}$ is determined through solving a single optimization problem, which is derived from
guaranteeing positivity of both density and pressure in the whole domain.
It is very efficient to solve this optimization problem compared to the base scheme due to the construction of the parameter.
The Lax-Friedrichs flux $\fv_{i+\hf}$ given in~\eqref{eq:LFFlux} 
is typically used as the low-order flux $\hat{\fv}^\text{low}_{i+\hf}$, with $\qv^{-}=\qv_{i}$ and $\qv^{+}=\qv_{i+1}$.
It has been proved that the limiter will guarantee positive numerical solutions
if the low-order flux is positivity-preserving.
For some simple problems such as scalar cases, it has been proved that a limiter using similar ideas achieves the designed accuracy of based schemes for smooth problems, see~\cite{Christlieb2015a} for instance.
For the ideal MHD equations, it has been numerically demonstrated in~\cite{Christlieb2016,Christlieb2015} 
that the corrected flux maintains the designed order of accuracy of the high-order flux.
For more details on the positivity-preserving limiter, see~\cite{Christlieb2016,Christlieb2015,Seal2014c}.


The implementation of the limiter on curvilinear meshes requires minor modifications.
For instance, curvilinear meshes require solving the ideal MHD equation in the form of~\eqref{eq:2DScalarEquationXiEta}.
Therefore, some steps of the limiter need to be modified accordingly. 
In the current multistage scheme, the limiter is applied at each stage. 
While in~\cite{Christlieb2015} it was sufficient to apply the limiter only at the final stage of each
time step, we find it necessary to apply it at each stage in the current scheme.  
This, again, is possibly due to the enhanced resolution provided by the HLL-type fluxes.

%% file: texSISC/boundaryConditions.tex
\subsection{Numerical boundary condition for PEC}
\label{sec:boundaryConditions}

In this section we derive a numerical boundary condition on a PEC boundary for the conserved quantities and the magnetic potential.
A slip-wall numerical boundary condition for the Euler equations  was given in~\cite{Henshaw2006moving}
and a PEC numerical boundary condition for Maxwell's equations was given in~\cite{max2006b}.
Following the ideas in~\cite{max2006b, Henshaw2006moving}, we derive a numerical PEC boundary condition for the ideal MHD equations.
Note that the idea presented here is not limited to the PEC boundary condition for ideal MHD and it is extendible to some non-ideal MHD cases and other physical boundary conditions.

Without loss of generality, assume the PEC boundary is located at $\xi = \xi_0$. Its normal direction at the boundary is $\nv = -\nabla_{\xv} \xi/\|\nabla_{\xv} \xi\|$. 
We further define the following quantities
\begin{alignat*}{3}
    \bar u_1 & := \nabla_{\xv}\xi \cdot (u,v)^T,           &  \bar u_2 & := \nabla_{\xv}\eta \cdot (u,v)^T, \\
    \bar B_1 & := \nabla_{\xv}\xi \cdot (B_1,B_2)^T, \quad &  \bar B_2 & := \nabla_{\xv}\eta \cdot (B_1,B_2)^T,
\end{alignat*}
which are proportional to the components normal to curves $\xi$ or $\eta$ being constant, respectively.
On a PEC boundary two normal components satisfy
\begin{align}
    \bar u_1 (\xi_0, \eta, t) = 0, \label{eqn:analyticalPECV} \\ 
    \bar B_1 (\xi_0, \eta, t) = 0. \label{eqn:analyticalPECB}
\end{align}
Those are the only analytical boundary conditions needed on a PEC boundary.
A local characteristic analysis reveals that there are one characteristic (forward fast magnetosonic waves) moving into the domain and one characteristic (backward fast magnetosonic waves)
moving out, with the rest  propagating along the boundary.
Note that the magnetic boundary condition~\eqref{eqn:analyticalPECB} is involved to guarantee four characteristics (Alfv\'en and slow magnetosonic waves) moving along the boundary.
Therefore, these boundary conditions are still consistent with the local characteristic analysis,
although there are two boundary conditions applied but only one characteristic moving in.
Due to many characteristic moving parallel to the boundary, some careful treatments are needed to implement the conditions numerically.

Following the ideas in~\cite{Henshaw2006moving}, we rely on the extrapolations and compatibility conditions to derive the numerical boundary condition with the help of ghost points.
Based on the boundary condition~\eqref{eqn:analyticalPECV} and \eqref{eqn:analyticalPECB} and the ideal MHD equations~\eqref{eq:MHD},
a compatibility condition for the total pressure $p_\text{tot}$ can be derived as
\begin{align}
    \partial_n p_\text{tot} = -\rho \bar u_2 \nv \cdot \partial_\eta \uv + \bar B_2 \nv \cdot \partial_\eta \Bv.
    \label{eqn:compatibility}
\end{align}
It is a direct extension of the well-known compatibility condition for the Euler equations involving the normal derivative of the pressure along a curved slip wall.
The details of the derivations are therefore omitted here.

Before discussing the numerical boundary conditions, 
we consider the analytical boundary condition for the magnetic potential $A$ in the constrained transport framework.
Along the boundary $\xi=\xi_0$, the magnetic boundary condition~\eqref{eqn:analyticalPECB} gives
\begin{align*}
 \bar B_1 =  \partial_x\xi\,\partial_yA-\partial_y\xi \,\partial_x A = (\partial_x\xi\,\partial_y\eta-\partial_x\eta \,\partial_y\xi)\partial_\eta A = 0,
\end{align*}
which further leads to
\begin{align*}
    \partial_\eta A = 0.
\end{align*}
Therefore, the potential at the boundary satisfies
\begin{align*}
    \partial_t A + (u \, \partial_x \xi + v \,\partial_y\xi) \partial_\xi A = 0. \qquad 
\end{align*}
Due to $u \, \partial_x \xi + v \,\partial_y\xi = \nabla_\xv \xi \cdot \uv$ and the velocity boundary condition~\eqref{eqn:analyticalPECV}, the analytical boundary condition for the potential therefore is
\begin{align}
    A(\xi_0, \eta, t) = A^0,
    \label{eqn:analyticalPECA}
\end{align}
with $A^0$ being some constant given in the initial condition.

Based on the above derivations for the PEC boundary, we propose a numerical compatibility boundary condition for the PEC boundary as follows.
First,   the velocities on the boundary are projected such that
\begin{align*}
    \nv \cdot \uv_\iv = 0, \quad \nv \cdot \Bv_\iv = 0, \quad \text{on } i=0.
\end{align*}
A Dirichlet boundary condition is applied to the magnetic potential 
\begin{align*}
    A_\iv = A^0, \quad \text{on } i=0.
\end{align*}
Note that the divergence in a general curvilinear coordinate can be written in a conservative form as 
\begin{align}
    \nabla \cdot \Bv = J \sum_{m=1}^d \partial_{r_m} (\av_{r_m} \cdot \Bv), \qquad \av_{r_m}  = J^{-1} \nabla_\xv r_m.
    \label{}
\end{align}
To implement the boundary condition numerically, 2nd-order discrete operators are defined as
\begin{align*}
    D_{0\xi} q_\iv := \frac{q_{i+1,j}-q_{i-1,j}}{2 \Delta \xi}, \qquad
    D_{0\eta} q_\iv := \frac{q_{i,j+1}-q_{i,j-1}}{2 \Delta \eta}.
\end{align*}
Its normal derivative is defined as
\begin{align*}
    D_{0n} q_\iv := \nv \cdot \left( \partial_x \xi \,D_{0\xi}+\partial_x \eta \,D_{0\eta}, \partial_y\xi\, D_{0\xi} +\partial_y\eta \, D_{0\eta} \right) q_\iv.
\end{align*}
The divergence-free condition can be used to determine the magnetic field in the normal direction as
\begin{align}
    D_{0\xi} (\av_\xi \cdot \Bv_\iv) = - D_{0\eta} (\av_\eta \cdot \Bv_\iv), \quad \text{on } i = 0.
    \label{eqn:divergenceBoundary}
\end{align}
The compatibility condition is used to determine the total pressure at the ghost points
\begin{align}
    D_{0n} p_{\text{tot},\iv} = -\rho_\iv \bar u_{2,\iv}\, \nv_\iv \cdot D_{0\eta} \uv_\iv + \bar B_{2,\iv}\, \nv_\iv \cdot D_{0\eta} \Bv_\iv, \quad \text{on } i = 0.
    \label{eqn:compatibilityBoundary}
\end{align}
Note that the conditions~\eqref{eqn:divergenceBoundary} and~\eqref{eqn:compatibilityBoundary} are used to determine $p_\text{tot}$ and the magnetic field $\av_\xi \cdot \Bv$ in the normal direction at $i = -1$. 
The same procedure can be used to determine those values at $i=-2$ and $-3$ if centered differences of wider stencils are used on the left-hand side of~\eqref{eqn:divergenceBoundary} and~\eqref{eqn:compatibilityBoundary}.
In our implementation those approximations are all implemented in second-order for simplicity.
The rest quantities $\{\rho, \uv, B_t, B_3, A\}$ at the ghost points
$i = -1, -2, -3$ are determined by extrapolations.
Here $B_t$ stands for the magnetic field in the transpose direction of the interface.
The WENO extrapolation in~\cite{tan2010inverse} is used here as the limited extrapolation.
For the quantity $q \in \{\rho, \uv, B_t, B_3, A\}$, let 
\begin{align*}
    p_0(\xi) & = q_{0,j}, \\
    p_1(\xi) & = \frac{q_{1,j}-q_{0,j}}{\Delta \xi} \xi +  q_{0,j}, \\
    p_2(\xi) & = \frac{q_{0,j}-2 q_{1,j}+ q_{2,j}}{2\,\Delta \xi^2} \xi^2 + \frac{-3 q_{0,j}+4q_{1,j}-q_{2,j}}{2\, \Delta \xi} \xi +  q_{0,j},
\end{align*}
define the first, second and third-order extrapolations in the negative $\xi$ direction, respectively.
The limited extrapolation is defined by
\begin{align*}
    \bar{p}(\xi) = \omegaBar_0 p_0(\xi) + \omegaBar_1 p_1(\xi) + \omegaBar_2 p_2(\xi),
\end{align*}
where
\begin{align*}
    \omegaBar_{r}=\frac{\alphaBar_{r}}{\alphaBar_{0}+\alphaBar_{1}+\alphaBar_{2}}, \quad \alphaBar_{r}=\frac{\dBar_{r}}{(\betaBar_{r}+\epsilon)^2}, \quad r=0,1,2,
\end{align*}
$\dBar_0 = \Delta \xi^2$, $\dBar_1 = \Delta \xi$ and $\dBar_2 =1 - \Delta \xi -  \Delta \xi^2$.
The smoothness indicators here are given by
\begin{equation*}
\begin{aligned}
\betaBar_{0} = \Delta \xi^2,  \quad
\betaBar_{1} =  (q_{1,j}-q_{0,j})^2, \quad
\betaBar_{2} = \frac{13}{12}( q_{0,j} - 2q_{1,j} + q_{2,j} )^2 + ( 2q_{0,j} - 3q_{1,j} +q_{2,j} )^2.
\end{aligned}
\end{equation*}
More details can be found in~\cite{tan2010inverse}  and note that $\beta_2$ therein contains a typo, which has been fixed here.
We find this WENO extrapolation is slightly more robust than the limited extrapolation (a weighted average of first and third-order extrapolations) used in~\cite{Henshaw2006moving}
for the current work.
	
Here the extrapolations are used to determine the transpose magnetic field $B_t$ and $B_3$.
Note in~\cite{max2006b} the boundary conditions for the electric fields are implemented differently by taking another time derivative of the evolution equations.
However, for the MHD system we consider, it is easy to show that taking another time derivative of the magnetic induction equation does not provide a useful constraint. 
The primary reason is that the electric field is fully determined by an ideal Ohm's law $\Ev = \Bv \times \uv$ in the MHD system. 

The idea of the compatibility boundary condition is similar to the idea of the so-called inverse Lax-Wendroff method~\cite{tan2010inverse}.
Both methods convert the normal derivatives of certain variables to the tangential derivatives at the PEC or slip-wall boundary, see~\eqref{eqn:divergenceBoundary} and~\eqref{eqn:compatibilityBoundary} for instance. 
To extend it to a boundary condition of higher order (higher than the second-order version presented here), 
the compatibility condition relies on high-order discrete operators of a wider stencil to approximate the conditions 
such as~\eqref{eqn:divergenceBoundary} and~\eqref{eqn:compatibilityBoundary},
which typically results into global coupling of all the ghost points, see~\cite{max2006b} for instance.
On the other hand, the inverse Lax-Wendroff approach relies on high-order time-derivatives to avoid global coupling,
but for a PEC boundary condition considered here, such a procedure leads to very complicated algebraic relations.
In practice, extrapolations are used to approximate those high-order time-derivatives, which may cause a new issue of numerical stabilities.
Therefore, the extension to a higher-order boundary condition is challenging and remains an interesting line of future research.

In practice, a reflective boundary condition is typically used for the PEC boundary, 
see~\cite{de1999numerical, de2001} for instance.
The boundary condition is simply implemented by copying the solution at the grid points $(i,j)$ to the ghost points $(-i,j)$ and
changing the sign of the velocity and magnetic field that are perpendicular to the interface.
For the compressible Euler case, this boundary condition is only valid on straight walls and 
introduces a low-order error on curved walls, see~\cite{Henshaw2006moving} for instance.
For the ideal MHD equations considered here, besides the same issue from hydrodynamics, another  potential issue is that 
such a reflective magnetic field will also affect the divergence of the magnetic field and may lead to numerical instabilities
for certain methods. In Section~\ref{sec:bowShock}, we will use a bow shock benchmark to discuss 
those issues for high-order methods  in further details.

%% file: texSISC/numerical.tex
\newcommand\subcap[1]{#1}
\section{Numerical results}
  \label{sec:numericalResults}

In this section, numerical results are presented to demonstrate the accuracy and performance of the WENO scheme.
For ease of reference, 
the WENO scheme based on the alternative flux formulation is simply referred to as the WENO scheme here.
Since the schemes using different Riemann solvers produce similar solutions, 
we only present the numerical solutions using the Lax-Friedirchs and HLLD fluxes.
Throughout the simulations, a third-order TVD-RK method is used as the time integrator.
A CFL number of $0.5$ is typically used, with the largest wave speed estimated in the curvilinear coordinates.
Unless otherwise stated, the constrained transport and positivity-preserving limiter are turned on in numerical simulations, although those steps are not required for some easy problems.


  \input texSISC/alfvenWave

  \input texSISC/shockTube
\input texSISC/fieldloop

  \input texSISC/orszagTang

  \input texSISC/cloudShock

\input texSISC/rotor

\input texSISC/blastwave

  \input texSISC/bowShock

%% file: texSISC/alfvenWave.tex
  \subsection{2D smooth \Alfven~wave problem}
\label{sec:alfvenWave}
  We first consider the smooth \Alfven~wave problem on a curvilinear mesh. 
  This problem is used to verify the accuracy of the numerical schemes on general curvilinear meshes.
  The initial conditions is
  \begin{equation*}
      (\rho, u, v, w, p, B_1, B_2, B_3)(0, x, y)  = (1, 0, 0.1\sin(2\pi x), 0.1\cos(2\pi x), 0.1, 1,
      0.1\sin(2\pi x), 0.1\cos(2\pi x)),
  \end{equation*}
  and its initial magnetic potential is 
  \begin{equation*}
      A(0,x,y)  = y + 0.1 \frac{\cos(2\pi x)}{2 \pi} .
  \end{equation*}
  The exact solution is an \Alfven~wave propagating along $x$-direction with a wave speed of one.

{
\newcommand{\figWidth}{10cm}
\newcommand{\trimfig}[2]{\trimFig{#1}{#2}{0}{0}{0.1}{.8}}
\begin{figure}[htb]
\begin{center}
\begin{tikzpicture}[scale=1]
  \useasboundingbox (0,0) rectangle (8,4.);  
  \draw(-1,-1) node[anchor=south west,xshift=0pt,yshift=+0pt] {\trimfig{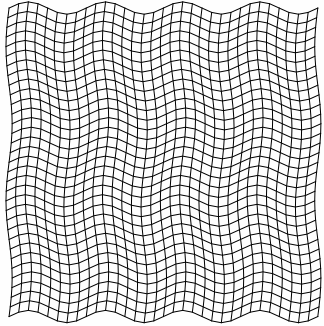}{\figWidth}};
%
\end{tikzpicture}
\end{center}
    \caption{Computational grid for the 2D \Alfven~wave problem.
      The (coarse) grid is of size $32 \times 32$  in the computational domain of $(\xi, \eta) \in {[0,1]}^2$ with
      $x = \xi + \epsilon_x \sin(2 \pi \, \eta  \, a_x)$ and
      $y = \eta +  \epsilon_y\sin(2 \pi \, \xi \, a_y)$.}\label{fig:2DAlfvenPerturbedMesh}
\end{figure}
}

In previous work such as~\cite{Christlieb2016,Christlieb2014} the direction of
 the \Alfven~wave was rotated so that it is not parallel to any grid lines in a Cartesian grid.  
 In the current test, the same goal is achieved by keeping the direction of the wave
  parallel to $x$-direction but perturbing the Cartesian grid.
 The computational domain is set to be $(\xi, \eta) \in {[0,1]}^2$, with the grid lines perturbed according to the mapping
  \begin{equation*}
    \begin{aligned}
      x & = \xi + \epsilon_x \sin ( 2 \pi \, \eta \, a_x), \\
      y & = \eta + \epsilon_y \sin ( 2 \pi\, \xi\, a_y),
    \end{aligned}
  \end{equation*}
  where $\epsilon_x$ and $\epsilon_y$ are the magnitude of
  perturbation and $a_x$ and $a_y$ are the wave numbers of the
  perturbation.  In the results presented below,  the parameters are taking by $\epsilon_x = 0.01$, $\epsilon_y = 0.02$, $a_x = 2$, and $a_y = 4$.
  As an illustration, a coarse grid of size $32 \times 32$ is presented in Figure~\ref{fig:2DAlfvenPerturbedMesh}.
  The boundary condition are all periodic for this smooth test.
{
  \begin{table}[htb]
    \centering
    \caption{$L^{\infty}$-errors of the 2D smooth \Alfven~wave problem. 
 }
 \label{tab:2DAlfven}
    \tableFont
    \begin{tabular}{|c|c|c|c|c|c|c|}
        \hline
        \multicolumn{7}{|c|}{WENO with Lax-Friedrichs flux} \\ \hline
      Mesh                  & Error in $\uv$        &   Order &   Error in $\Bv$        &   Order  &   Error in $A$         &  Order \\
        \hline
        \hline
        $32 \times 32$	  &  \num{3.324}{-3} &   ---	          &   \num{5.131}{-3} &   ---	          &   \num{1.560}{-4}  &   ---	\\
        \hline                                                      
        $64 \times 64$	  &  \num{8.234}{-5}   &   5.34	  &   \num{4.090}{-4}   &   3.65	  &   \num{8.394}{-6}  &   4.22	\\
        \hline                                                      
      $128 \times 128$	  &  \num{6.713}{-6}   &   3.62	  &   \num{2.658}{-5}   &   3.94	  &   \num{5.266}{-7}  &   3.99	\\
        \hline                                                      
      $256 \times 256$	  &  \num{4.544}{-7}   &   3.88	  &   \num{1.677}{-6}   &   3.99	  &   \num{3.322}{-8}  &   3.99	\\
        \hline
    \end{tabular}

\medskip

    \begin{tabular}{|c|c|c|c|c|c|c|}
        \hline
        \multicolumn{7}{|c|}{WENO with HLLD flux} \\ \hline
      Mesh                  &  Error in $\uv$        &   Order & Error in $\Bv$        &   Order  &   Error in $A$         &  Order \\
        \hline
        \hline
        $32 \times 32$	  &   \num{3.415}{-3} &   ---	 &  \num{5.257}{-3}  &   ---	  &   \num{1.604}{-4}   &   ---	\\
        \hline                                              
      $64 \times 64$	  &   \num{8.260}{-5}   &   5.37	  &\num{4.090}{-4}  &   3.68	  &   \num{8.388}{-6}   &   4.26	\\
        \hline                                              
      $128 \times 128$	  &   \num{6.716}{-6}   &   3.62	  &\num{2.658}{-5}  &   3.94	  &   \num{5.264}{-7}   &   3.99	\\
        \hline                                              
      $256 \times 256$	  &   \num{4.542}{-7}   &   3.89	  &\num{1.677}{-6}  &   3.99	  &   \num{3.322}{-8}   &   3.99	\\
        \hline
    \end{tabular}
 
  \end{table}
}

  A refinement study is conducted on a sequence of grids of increasing resolution to verify the accuracy
  of the WENO methods with two Riemann solvers.
  The numerical solutions are compared to the exact solutions at $t = 1$.
  Throughout the refinement study, a fixed CFL number of 0.6 is used to determine the time step. 
  Table~\ref{tab:2DAlfven} presents the $L^{\infty}$-errors of $\uv$, $\Bv$ and $A$ and the estimated convergence rates.
 The error of the vector is the maximum taken over the Euclidean norm of the vector.
  The results are obtained using the WENO methods with the Lax-Friedrichs flux and the HLLD flux.  
  The results confirm that the numerical schemes are both fourth-order accurate.  
  Recall that a fourth-order constrained transport method is used 
  and the resulting scheme in this work is fourth-order (in space).
  Note that the difference between the results of two fluxes is very small for this smooth problem.
  In the simulations, the constrained transport and positivity-preserving limiter are both turned on.

%% file: texSISC/shockTube.tex
 \subsection{Brio-Wu shock tube}
\label{sec:shockTube}

The second problem we consider is a commonly tested Riemann problem of the Brio-Wu shock tube test.  
The initial conditions in 1D are
  \begin{equation*}
    (\rho, u, v, w, p, B_1, B_2, B_3)
    =
    \begin{cases}
      (1, 0, 0, 0, 1, 0.75, 1, 0) &
      \text{if $x<0$,} \\
      (0.125, 0, 0, 0, 0.1, 0.75, -1, 0) & \text{if $x \geq 0$.}
    \end{cases}
  \end{equation*}
The schemes are first tested in 1D  on both uniform and non-uniform meshes,
and then they are tested on a 2D uniform mesh with the initial conditions rotated.

\subsubsection{1D shock tube}
\label{sec:1DShockTube}

{
\newcommand{\figWidth}{8.cm}
\newcommand{\trimfig}[2]{\trimFig{#1}{#2}{.0}{.0}{.0}{.0}}
\newcommand{\figWidthz}{3.8cm}
\newcommand{\trimfigz}[2]{\trimFig{#1}{#2}{.0}{.0}{.0}{.0}}
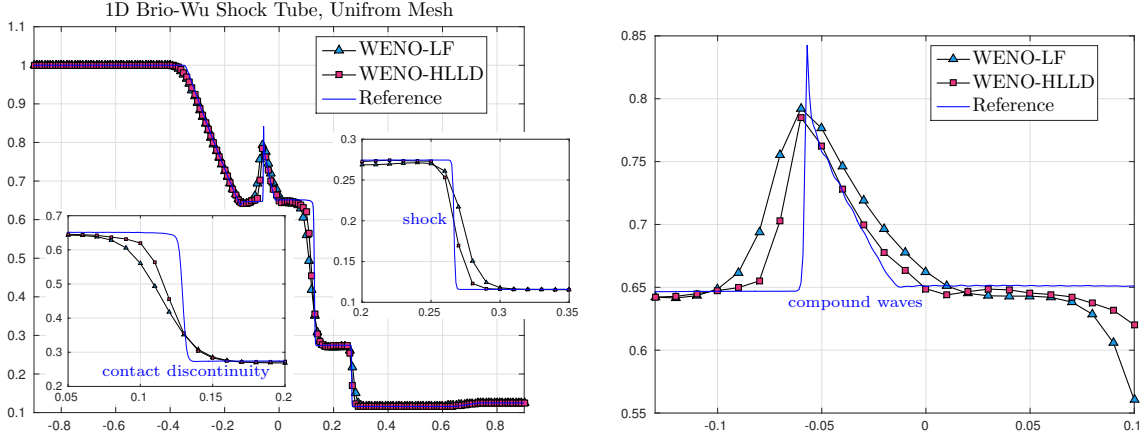
\begin{figure}[htb]
\begin{center}
\resizebox{14cm}{!}{
\begin{tikzpicture}[scale=1]
  \useasboundingbox (0.0,0) rectangle (16,6.2);  
\begin{scope}[yshift=.4cm]
  \draw(-1,-.75) node[anchor=south west,xshift=0pt,yshift=+0pt] {\trimfig{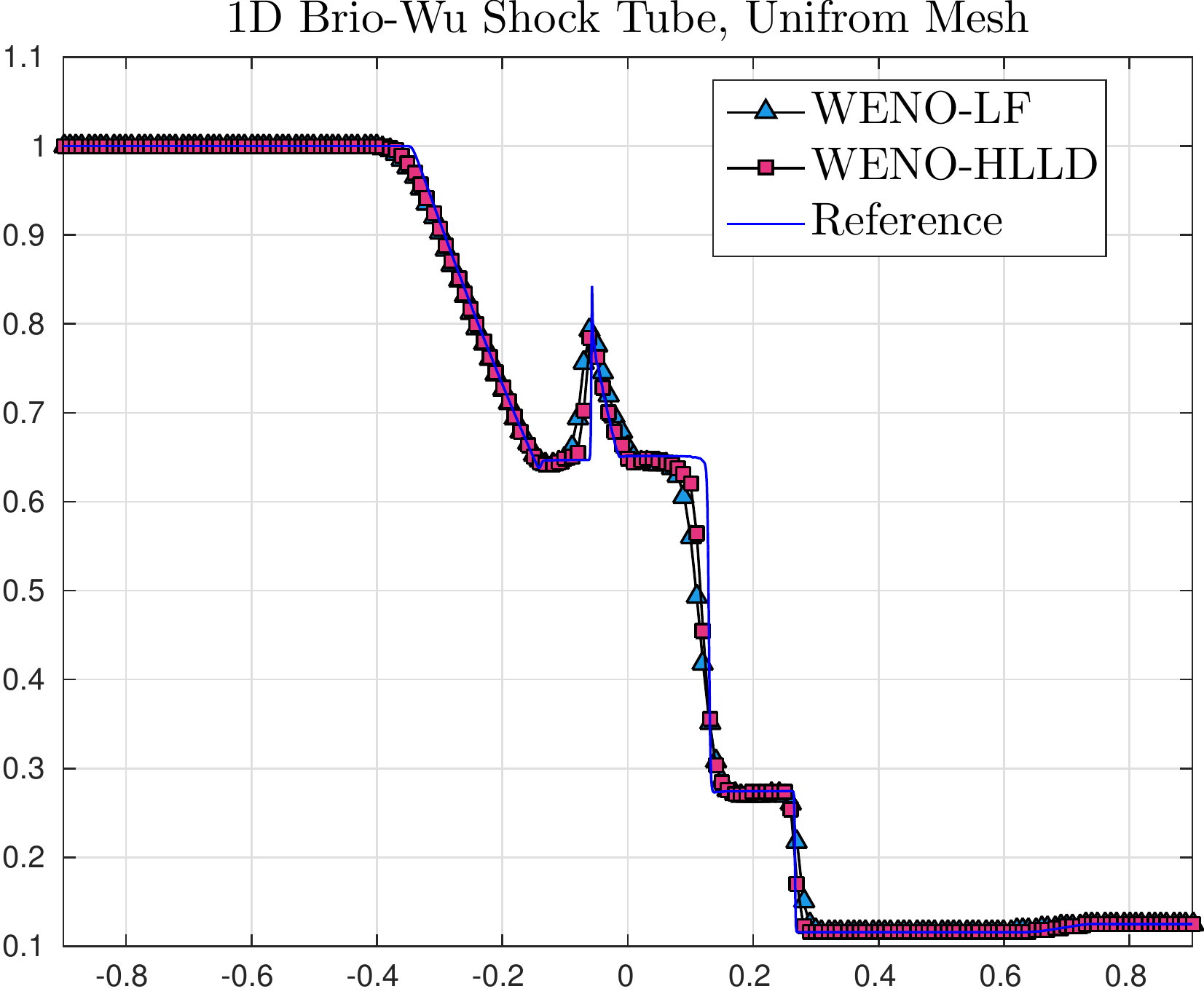}{\figWidth}};
  \draw(-.4, -.3) node[anchor=south west,xshift=0pt,yshift=+0pt] {\trimfigz{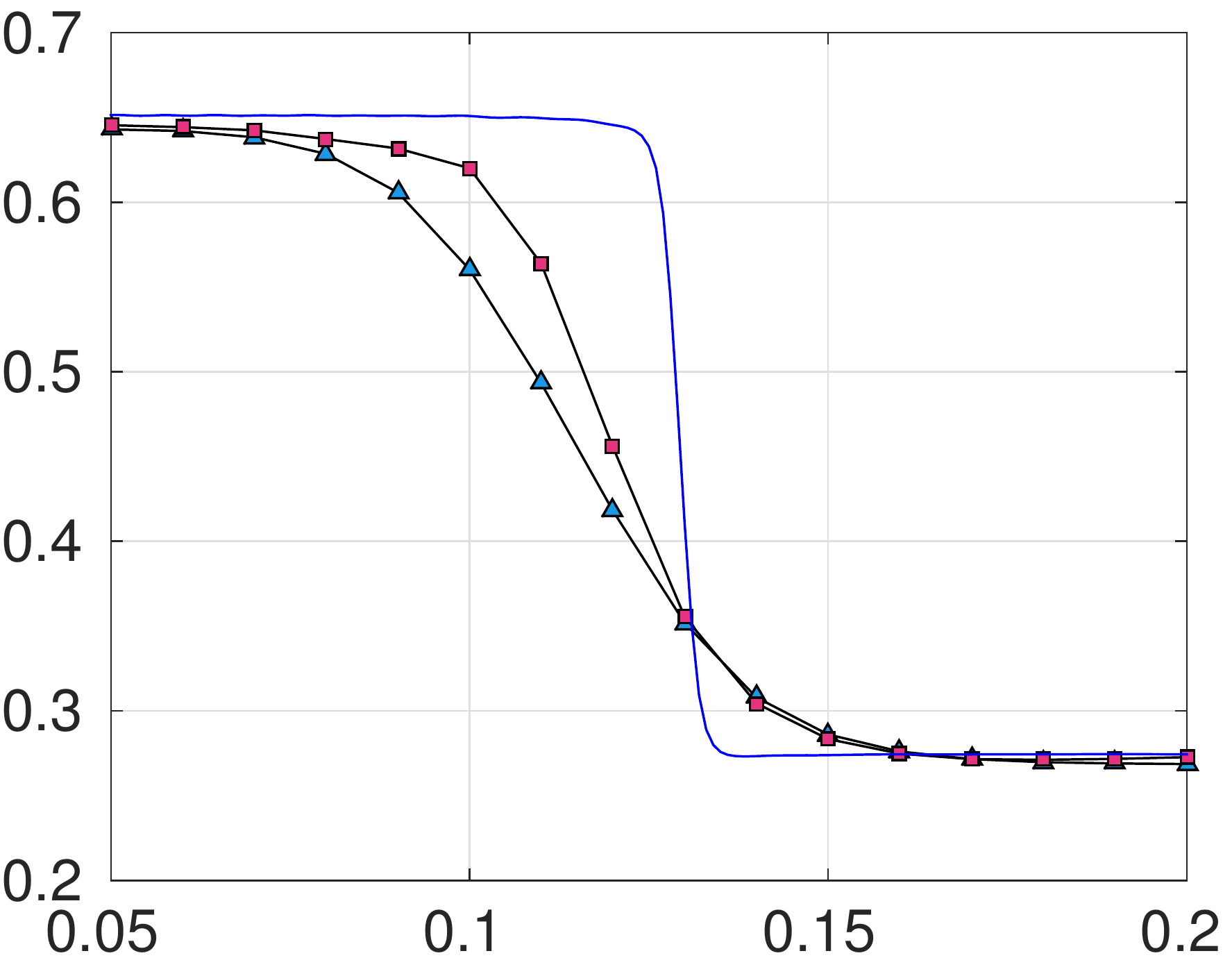}{\figWidthz}};
  \draw(4, 1.) node[anchor=south west,xshift=0pt,yshift=+0pt] {\trimfigz{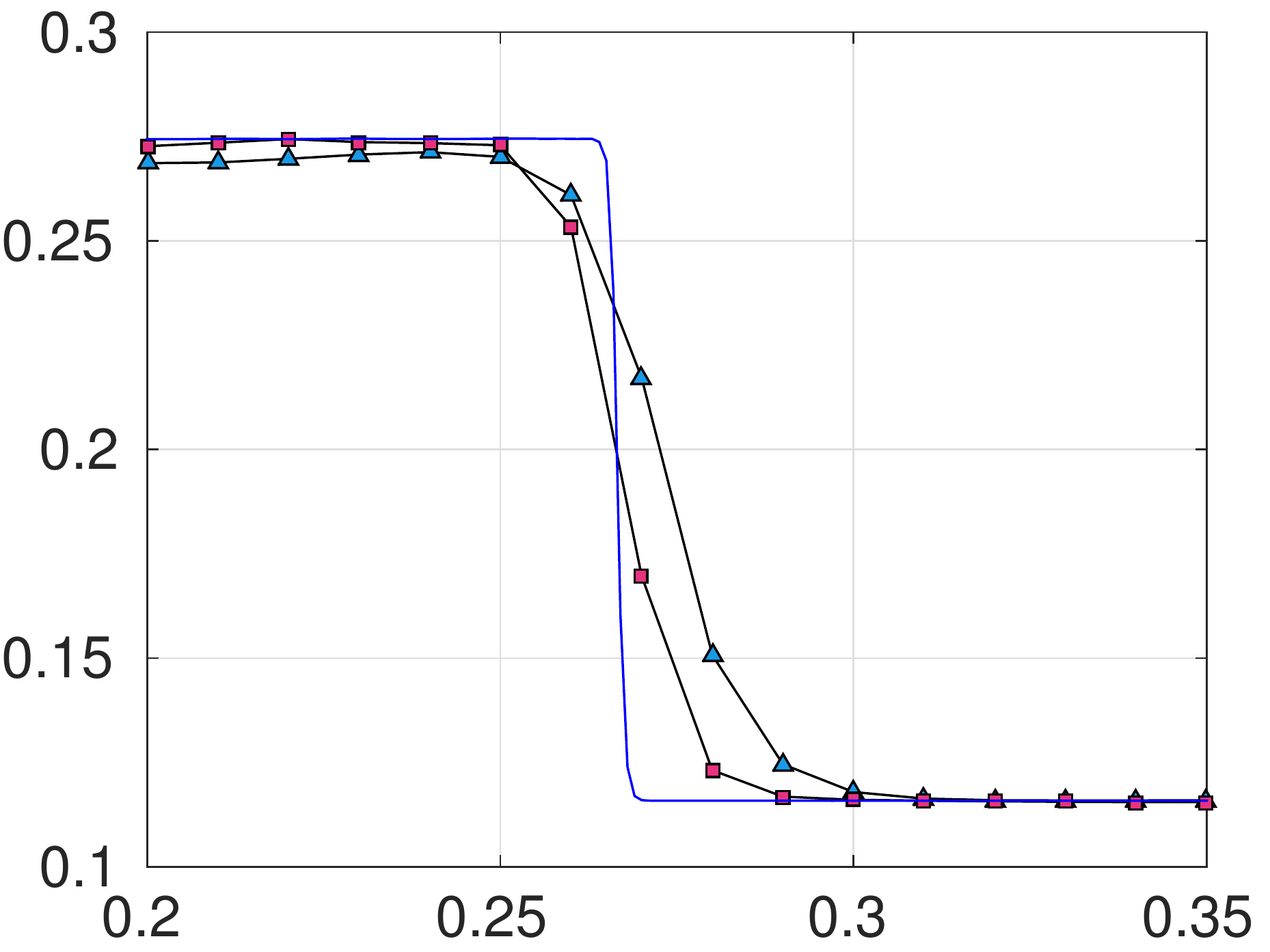}{\figWidthz}};
  \draw(8.4,-.75) node[anchor=south west,xshift=0pt,yshift=+0pt] {\trimfig{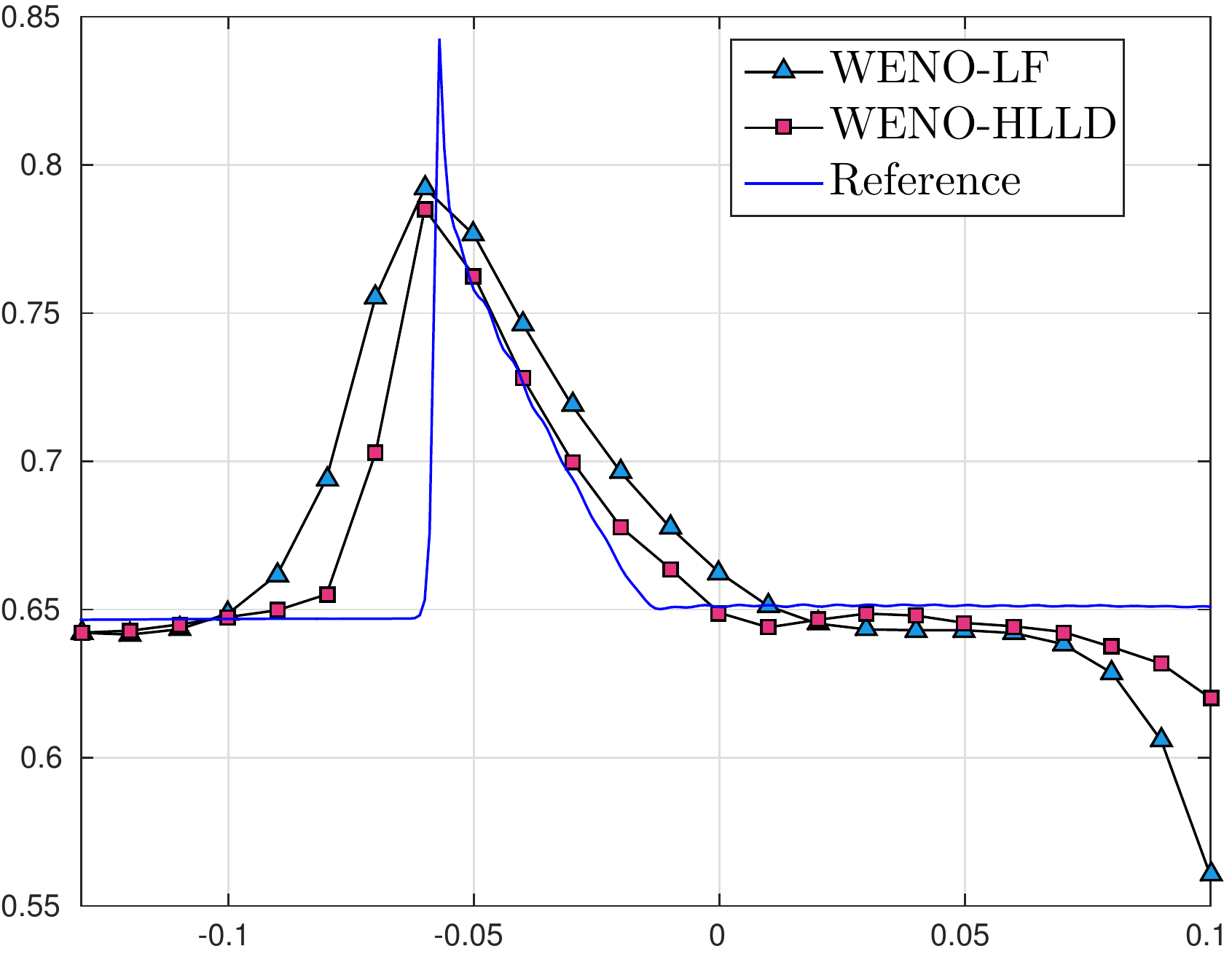}{\figWidth}};
  \draw(.5,.35) node[blue, anchor=west] {\scriptsize contact discontinuity};
  \draw(5.1,2.7) node[blue,anchor=west] {\scriptsize shock};
  \draw(11,1.4) node[blue, anchor=west] {\scriptsize compound waves};
  \end{scope}
\end{tikzpicture}
} 
\end{center}
    \caption{
1D Brio-Wu shock tube.
Density solved using the WENO schemes with two fluxes at $t = 0.2$ and its three zoomed views around the contact discontinuity, shock and compound waves.
The solutions are computed on a uniform mesh of 200 points.   
The reference solution is a numerical solution on a fine mesh of 2000 points.
}
\label{fig:1DBrioWuAltUniform}
\end{figure}
}

Figure~\ref{fig:1DBrioWuAltUniform}  present the density on a uniform mesh of 200 grid points.  
The solutions of the WENO schemes with the Lax-Friedrichs (WENO-LF) and HLLD (WENO-HLLD) fluxes are presented.
The results are compared to a reference solution on a very fine mesh solved using the scheme in~\cite{Christlieb2014}.
The numerical results of two fluxes match well with the reference solution as well as other numerical results in the literature.
As shown from the zoomed views in Figure~\ref{fig:1DBrioWuAltUniform}, 
the numerical solutions of the HLLD flux show less smeared structures around the contact discontinuity, shock and compound waves.
The solutions of the HLLD flux around rarefaction (not presented in the zoomed views) also show some improvements over the solutions of the Lax-Friedrichs flux.
This 1D results show the low-dissipative solver performs better in problems involving shocks or contact discontinuities.

{
\newcommand{\figWidth}{8.cm}
\newcommand{\trimfig}[2]{\trimFig{#1}{#2}{.0}{.0}{.0}{.0}}
\newcommand{\figWidthz}{3.8cm}
\newcommand{\trimfigz}[2]{\trimFig{#1}{#2}{.0}{.0}{.0}{.0}}
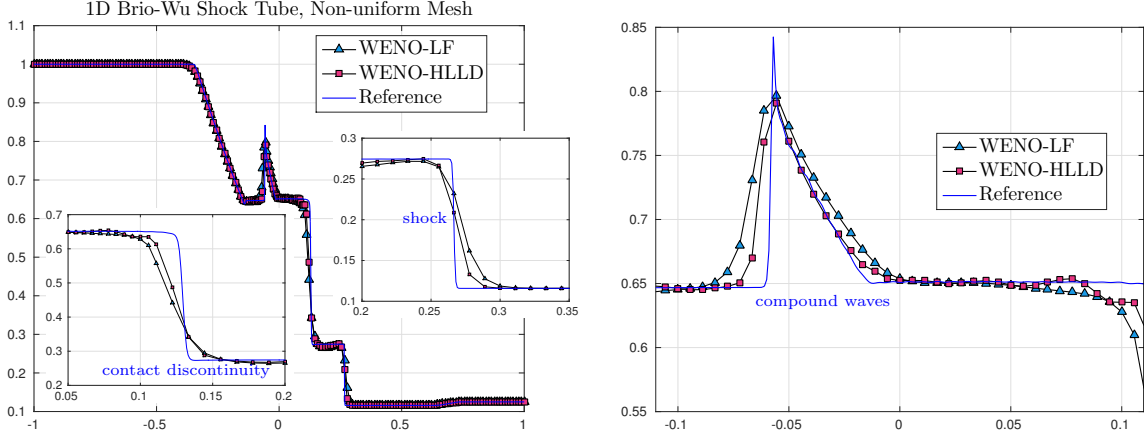
\begin{figure}[htb]
\begin{center}
\resizebox{14cm}{!}{
\begin{tikzpicture}[scale=1]
  \useasboundingbox (0.0,0) rectangle (16,6.2);  
\begin{scope}[yshift=.4cm]
  \draw(-1,-.75) node[anchor=south west,xshift=0pt,yshift=+0pt] {\trimfig{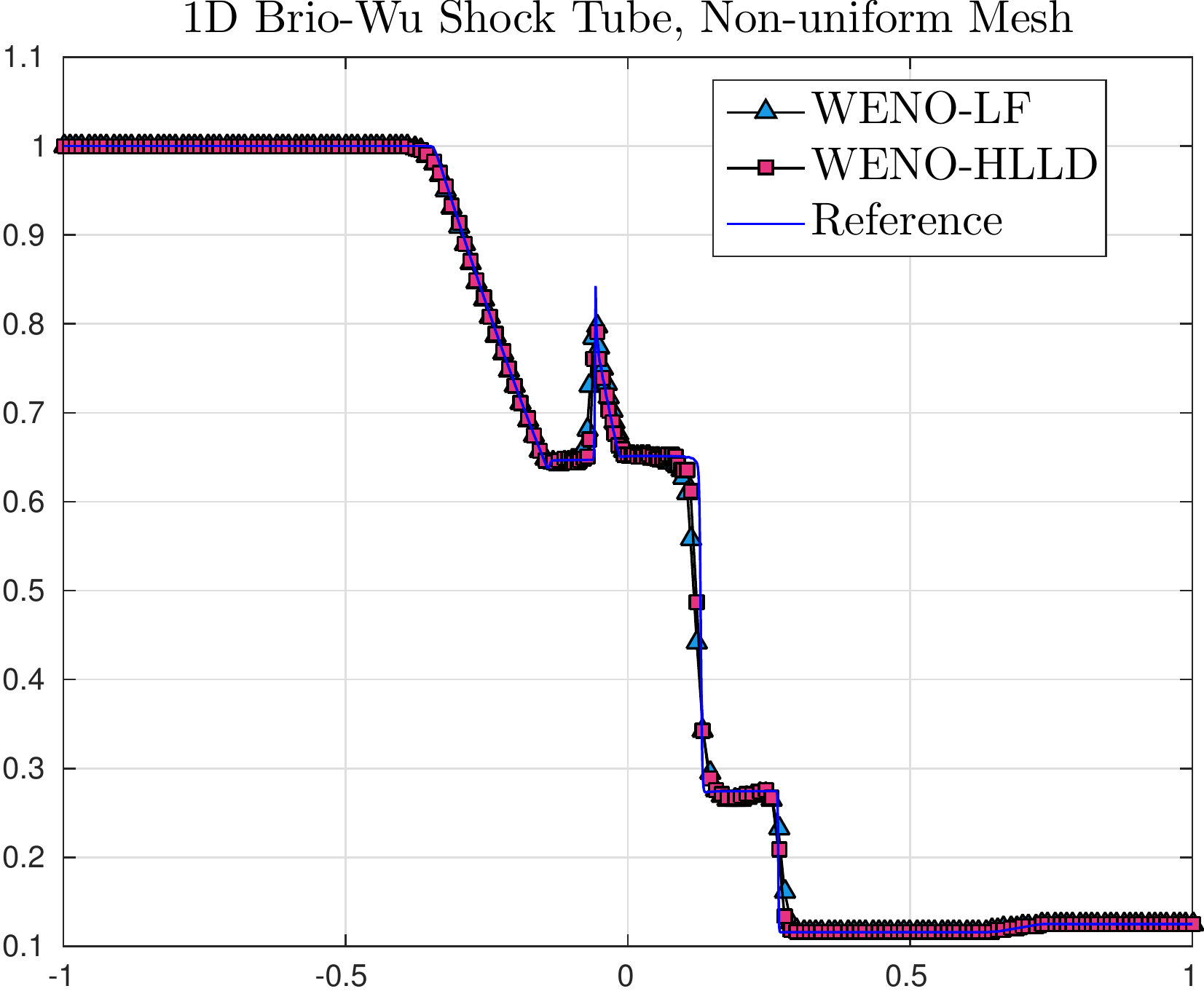}{\figWidth}};
  \draw(-.4, -.3) node[anchor=south west,xshift=0pt,yshift=+0pt] {\trimfigz{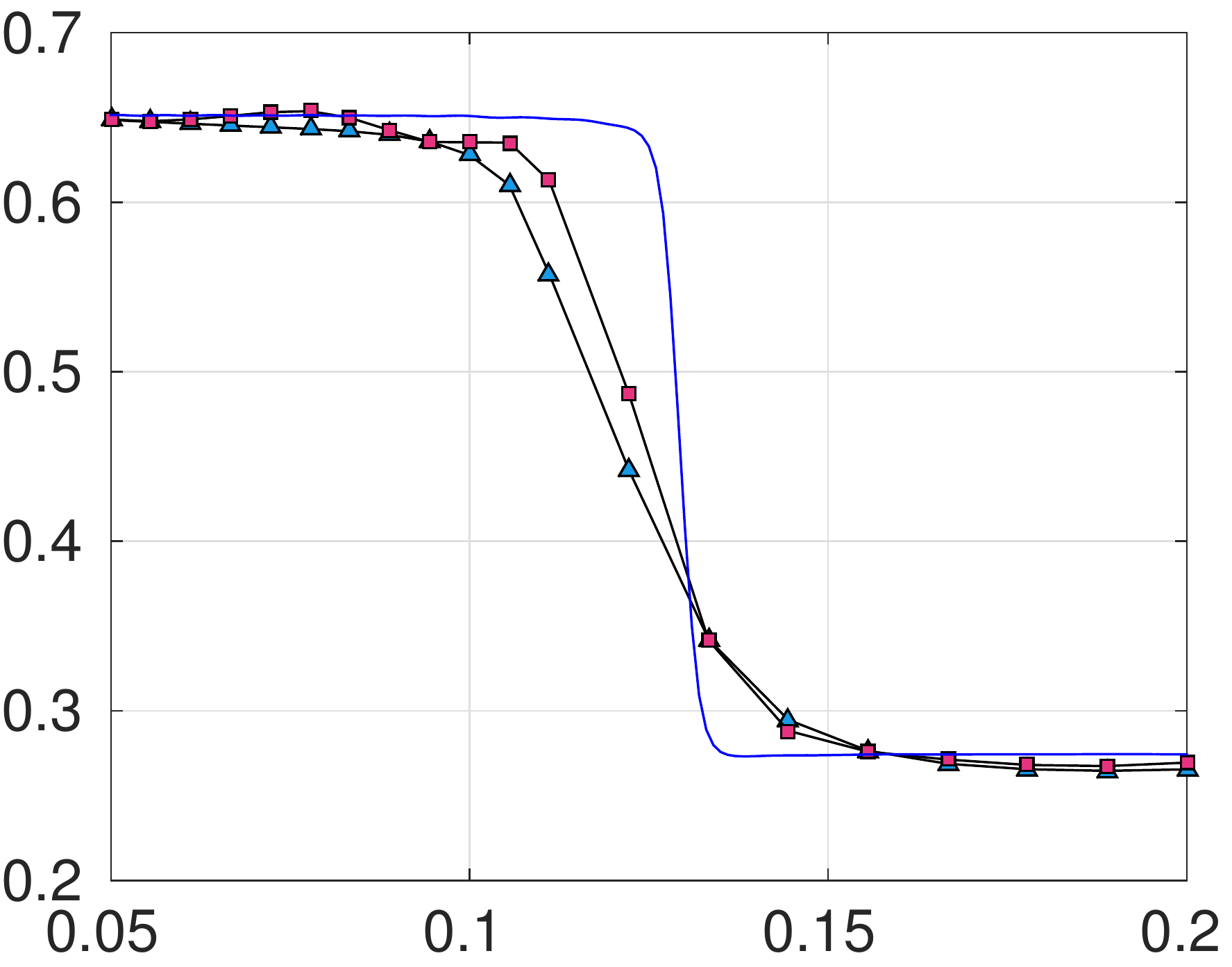}{\figWidthz}};
  \draw(4, 1.) node[anchor=south west,xshift=0pt,yshift=+0pt] {\trimfigz{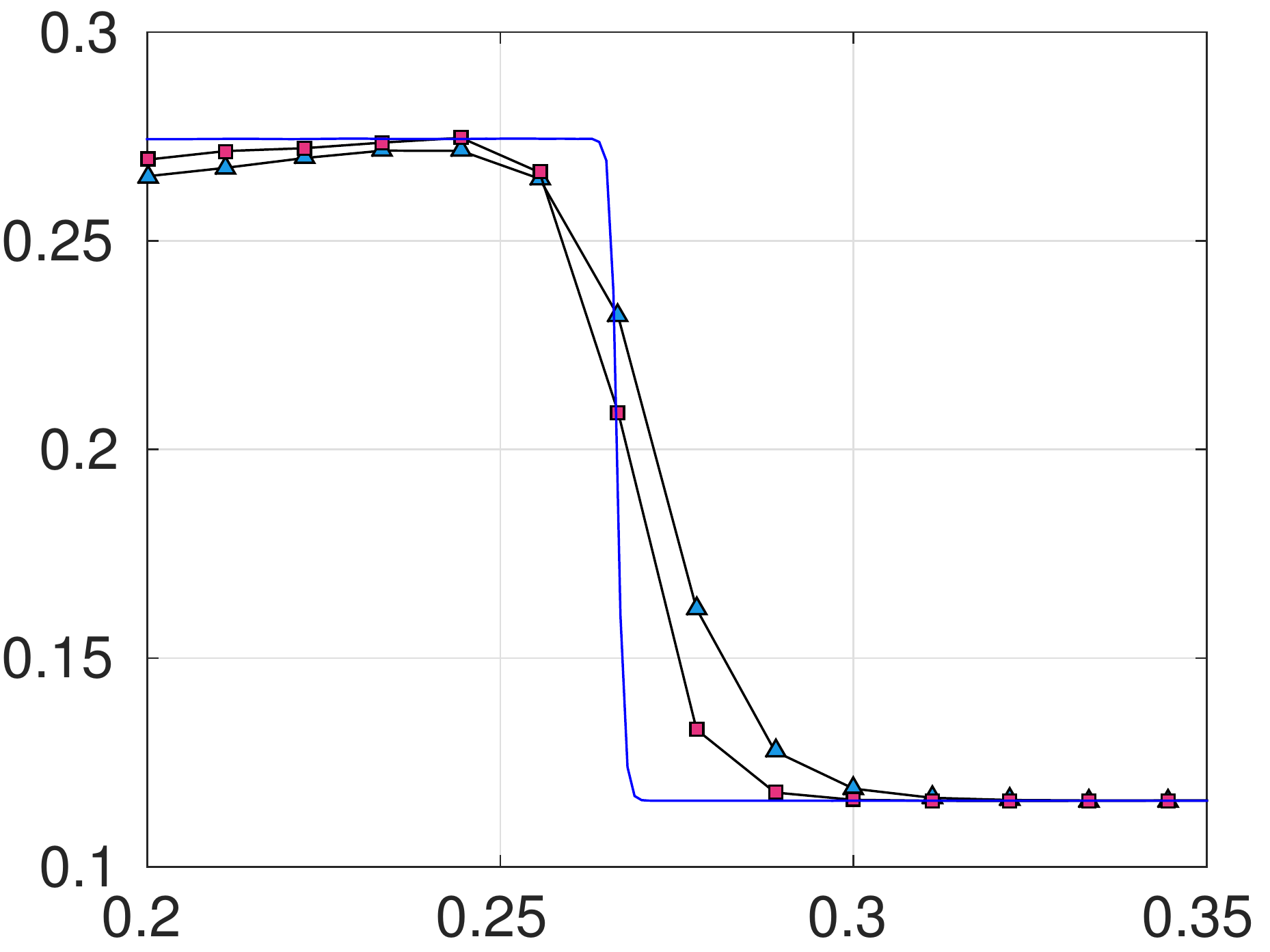}{\figWidthz}};
  \draw(8.4,-.75) node[anchor=south west,xshift=0pt,yshift=+0pt] {\trimfig{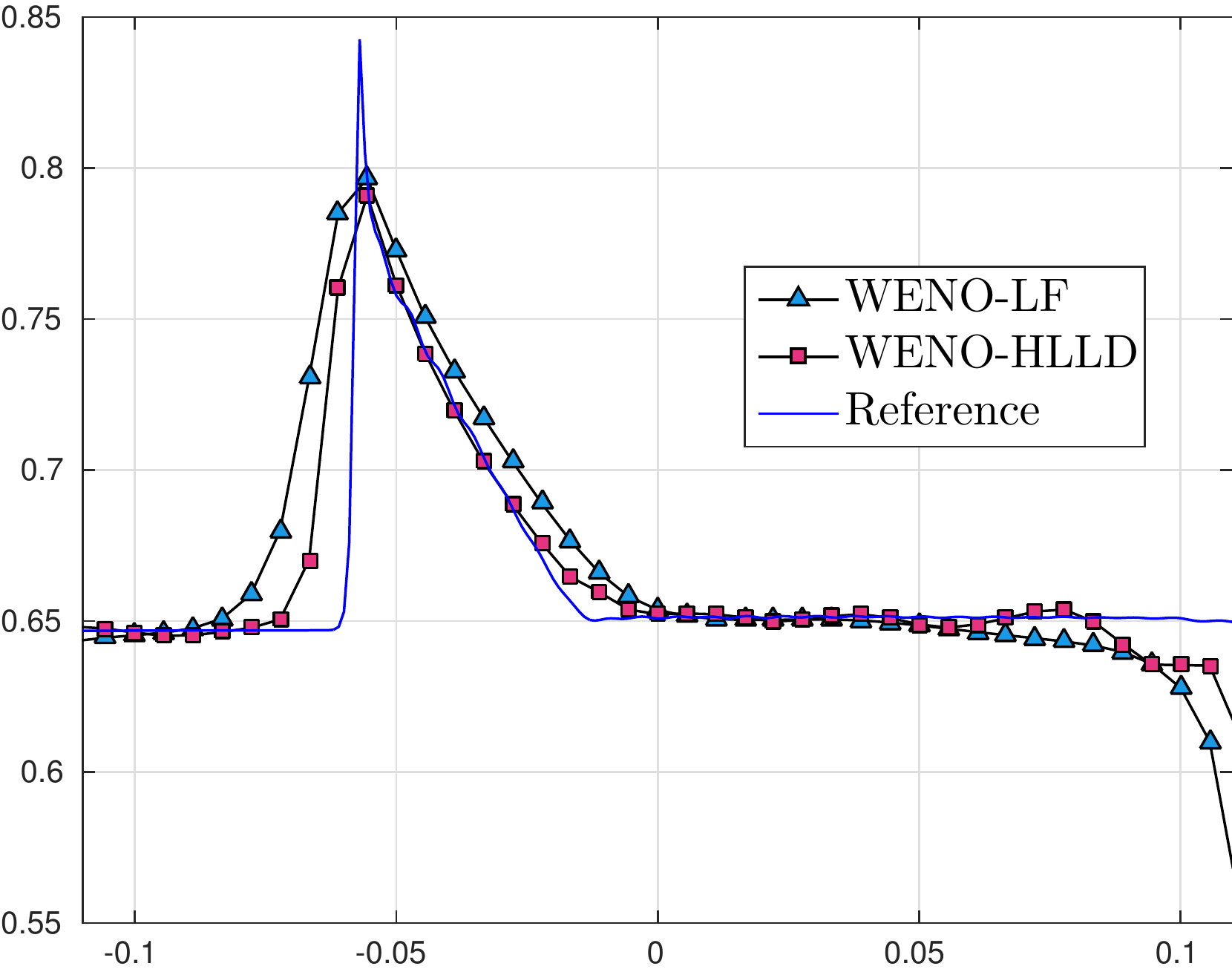}{\figWidth}};
  \draw(.5,.35)   node[blue, anchor=west] {\scriptsize contact discontinuity};
  \draw(5.1,2.7)  node[blue, anchor=west] {\scriptsize shock};
  \draw(10.5,1.4) node[blue, anchor=west] {\scriptsize compound waves};
  \end{scope}
\end{tikzpicture}
} 
\end{center}
    \caption{
1D Brio-Wu shock tube.
Density solved using the WENO schemes with two fluxes at $t = 0.2$ and its three zoomed views around the contact discontinuity, shock and compound waves.
The solutions are computed on a non-uniform mesh of 200 points.   
The reference solution is a numerical solution on a fine mesh of 2000 points.
}
\label{fig:1DBrioWuAltNonUniform}
\end{figure}
}

Next we examine the scheme on a non-uniform mesh given by the mapping
\begin{equation*}
  x = \begin{cases}
    \frac{5}{9} \xi & \text{if ${\left| \xi \right|} \leq 0.2$,} \\
    \sign(\xi) \left( \frac{1}{9} + \frac{10}{9} \left( \left| \xi
        \right| - 0.2 \right) \right) & \text{otherwise,}
  \end{cases}
\end{equation*}
with $-1 \leq \xi \leq 1$.
Note that the mesh is clustered around the region $[-0.11, 0.11]$ by a factor of $9/5$ and coarsened by a factor of 0.9 at the remaining region.
Figure~\ref{fig:1DBrioWuAltNonUniform} presents the density on a non-uniform mesh of 200 grid points that are solved using the WENO schemes with two fluxes.
The results are similar to the results on uniform meshes. 
The solutions also resolve the contact discontinuity and compound waves better due to the clustered grids points around those regions.
We also note that both the WENO schemes can handle the abrupt change in the grid spacing 
of the non-uniform meshes.

\subsubsection{2D rotated shock tube}
  \label{sec:Alt2DRotatedShockTube}
{
\newcommand{\figWidth}{8cm}
\newcommand{\trimfig}[2]{\trimFig{#1}{#2}{.3}{.0}{.15}{.0}}
\begin{figure}[htb]
\begin{center}
\begin{tikzpicture}[scale=1]
  \useasboundingbox (0,0) rectangle (16,3.5);  
  \draw(-.5,-.4) node[anchor=south west,xshift=0pt,yshift=+0pt] {\trimfig{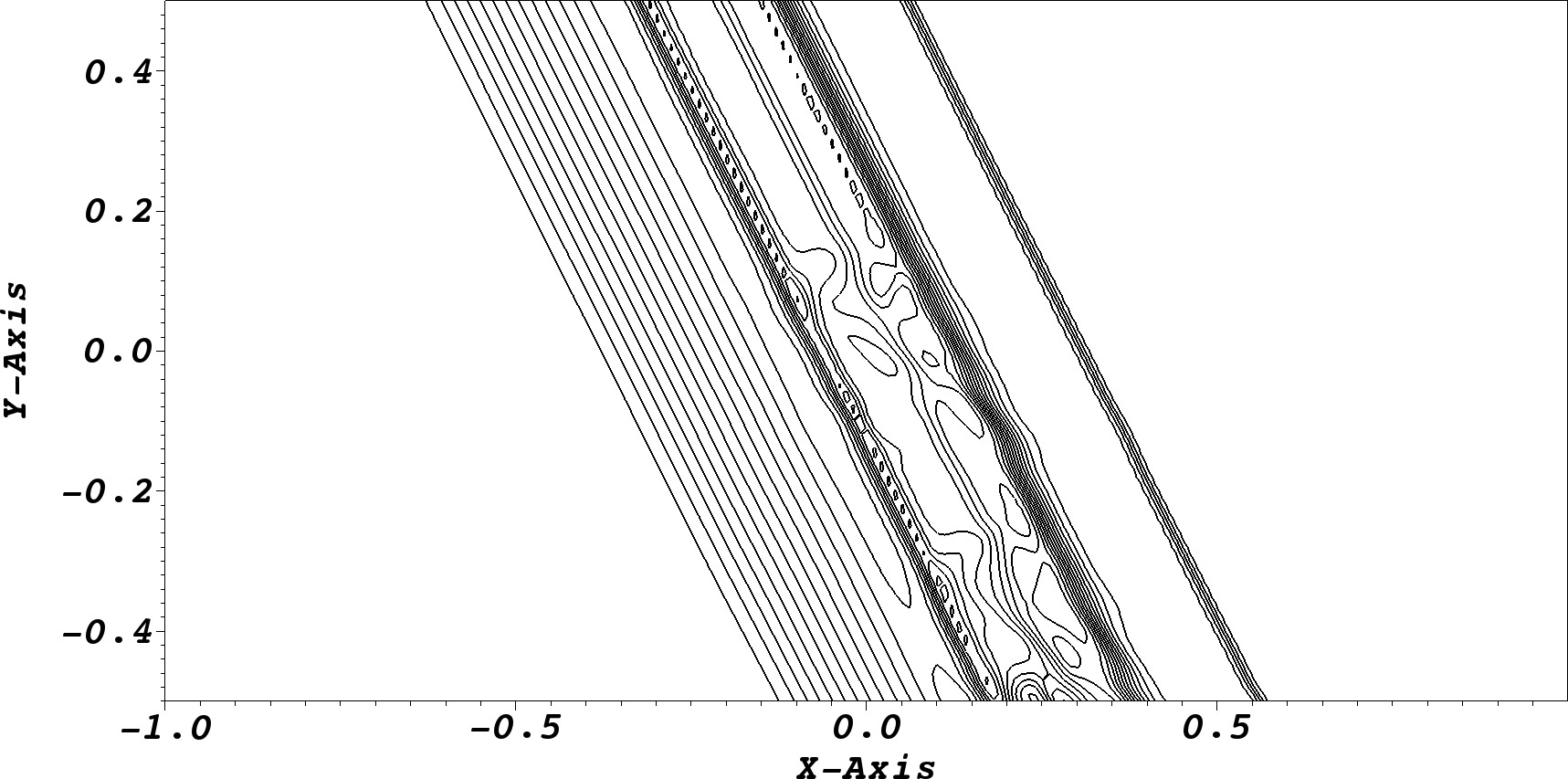}{\figWidth}};
  \draw(8.,-.4) node[anchor=south west,xshift=0pt,yshift=+0pt] {\trimfig{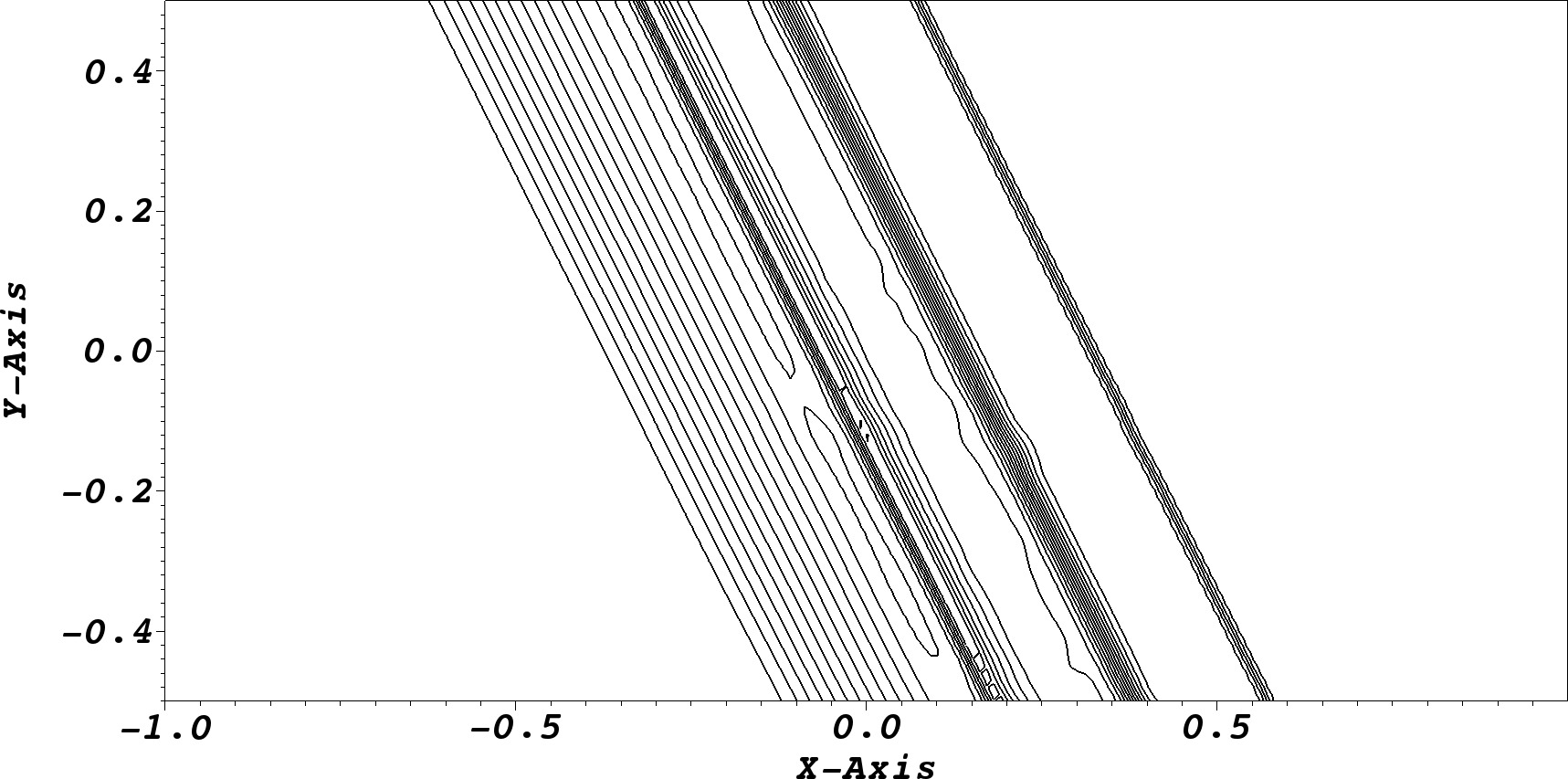}{\figWidth}};
  \draw(0.2,3.4) node[draw,fill=white,anchor=west,xshift=2pt,yshift=-1pt] {\small $\rho$};
  \draw(8.7,3.4) node[draw,fill=white,anchor=west,xshift=2pt,yshift=-1pt] {\small $\rho$};
%
\end{tikzpicture}
\end{center}
    \caption{
    2D rotated Brio-Wu shock tube.  Density at $t=0.2$. 
    Left: constrained transport and positivity-preserving limiter turned off.
    Right: constrained transport and positivity-preserving limiter turned on.
The WENO scheme with the Lax-Friedrichs flux is used on a uniform grid of size $200 \times 100$.
}
\label{fig:2DRotatedBrioWuContour}
\end{figure}
}

{
\newcommand{\figWidth}{8.cm}
\newcommand{\trimfig}[2]{\trimFig{#1}{#2}{.0}{.0}{.0}{.0}}
\newcommand{\figWidthz}{3.6cm}
\newcommand{\trimfigz}[2]{\trimFig{#1}{#2}{.0}{.0}{.0}{.0}}
\begin{figure}[htb]
\begin{center}
\resizebox{14cm}{!}{
\begin{tikzpicture}[scale=1]
  \useasboundingbox (0.0,0) rectangle (16,6.2);  
\begin{scope}[yshift=.4cm]
  \draw(-1,-.75) node[anchor=south west,xshift=0pt,yshift=+0pt] {\trimfig{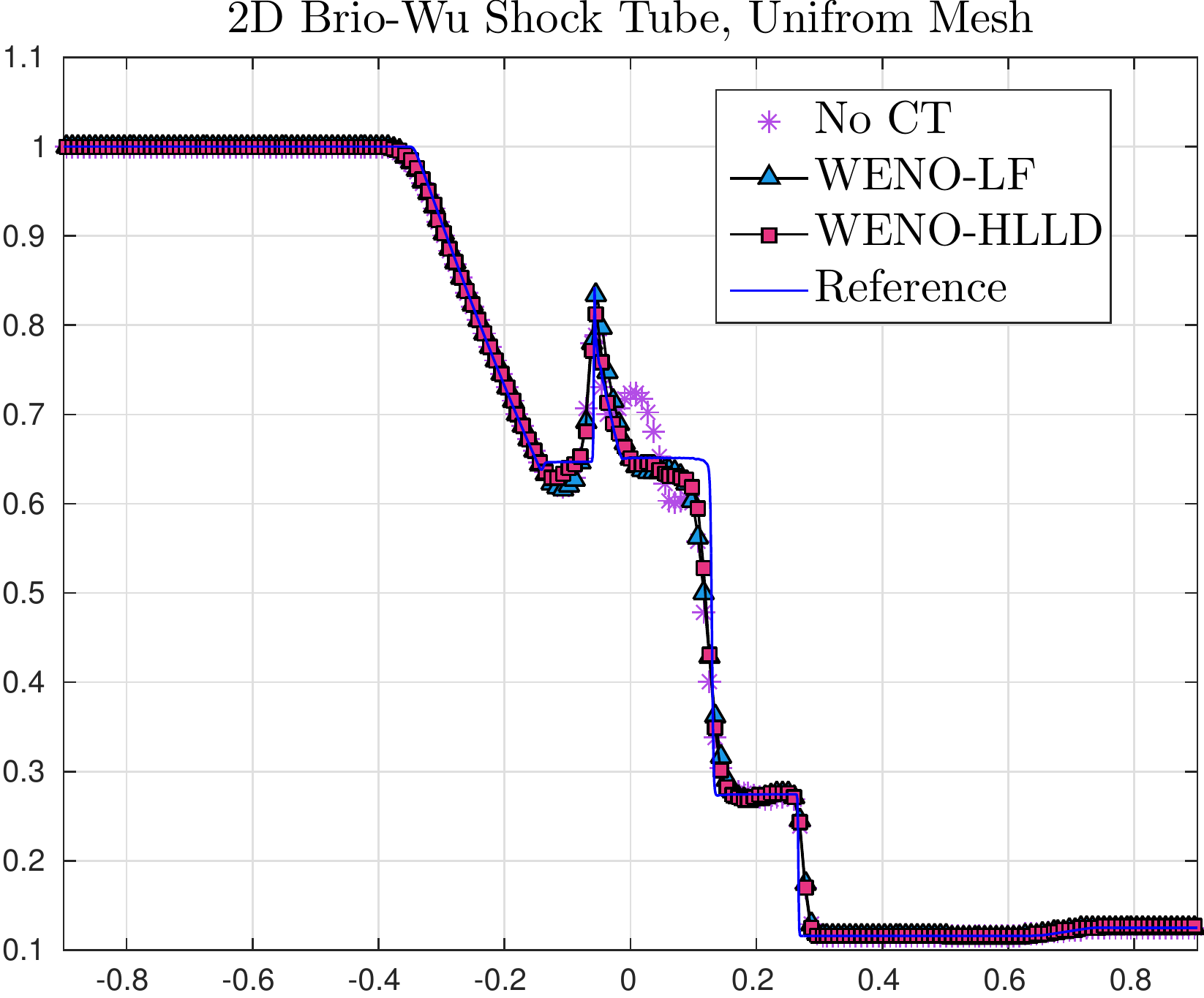}{\figWidth}};
  \draw(-.4, -.3) node[anchor=south west,xshift=0pt,yshift=+0pt] {\trimfigz{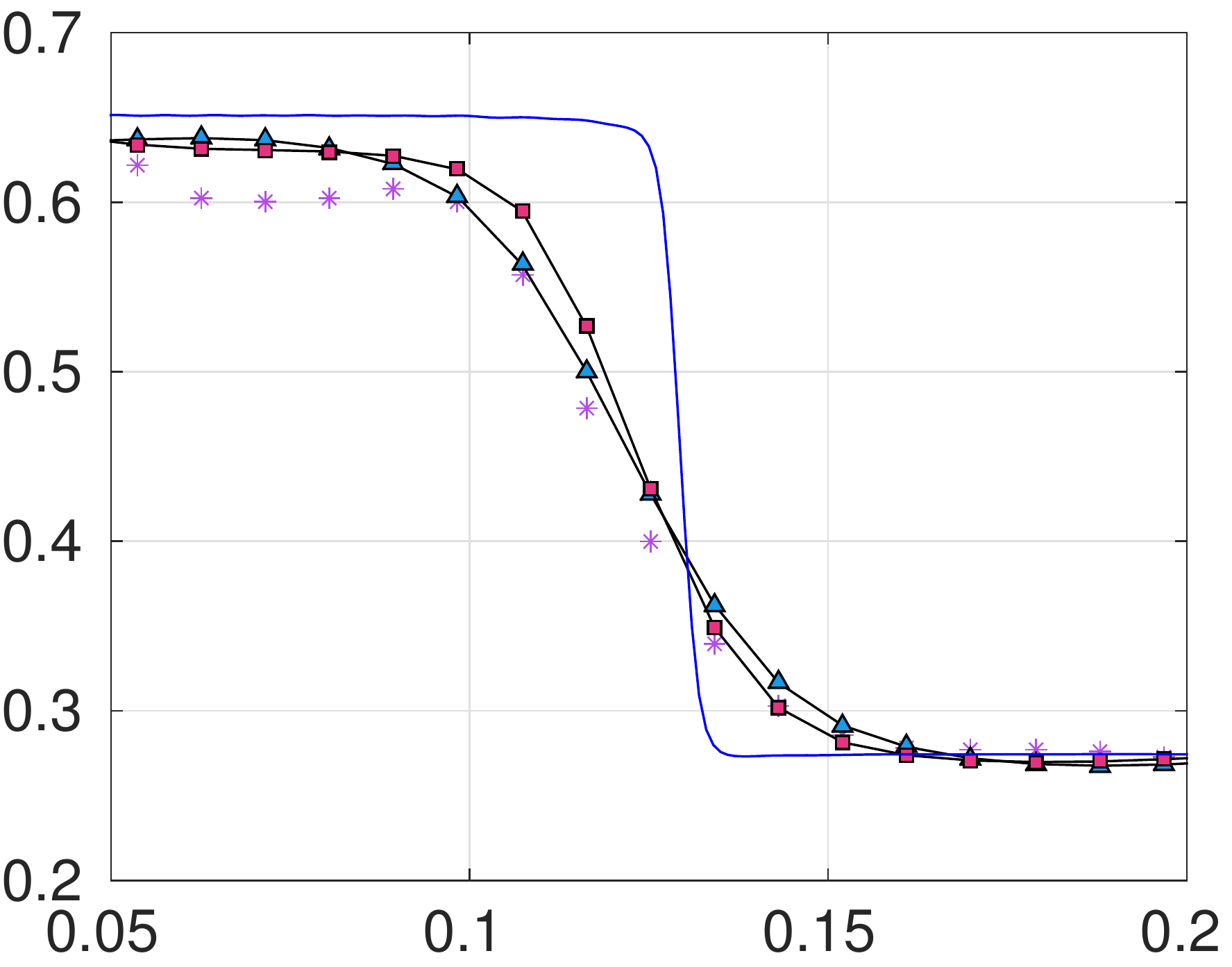}{\figWidthz}};
  \draw(4.2, .9) node[anchor=south west,xshift=0pt,yshift=+0pt] {\trimfigz{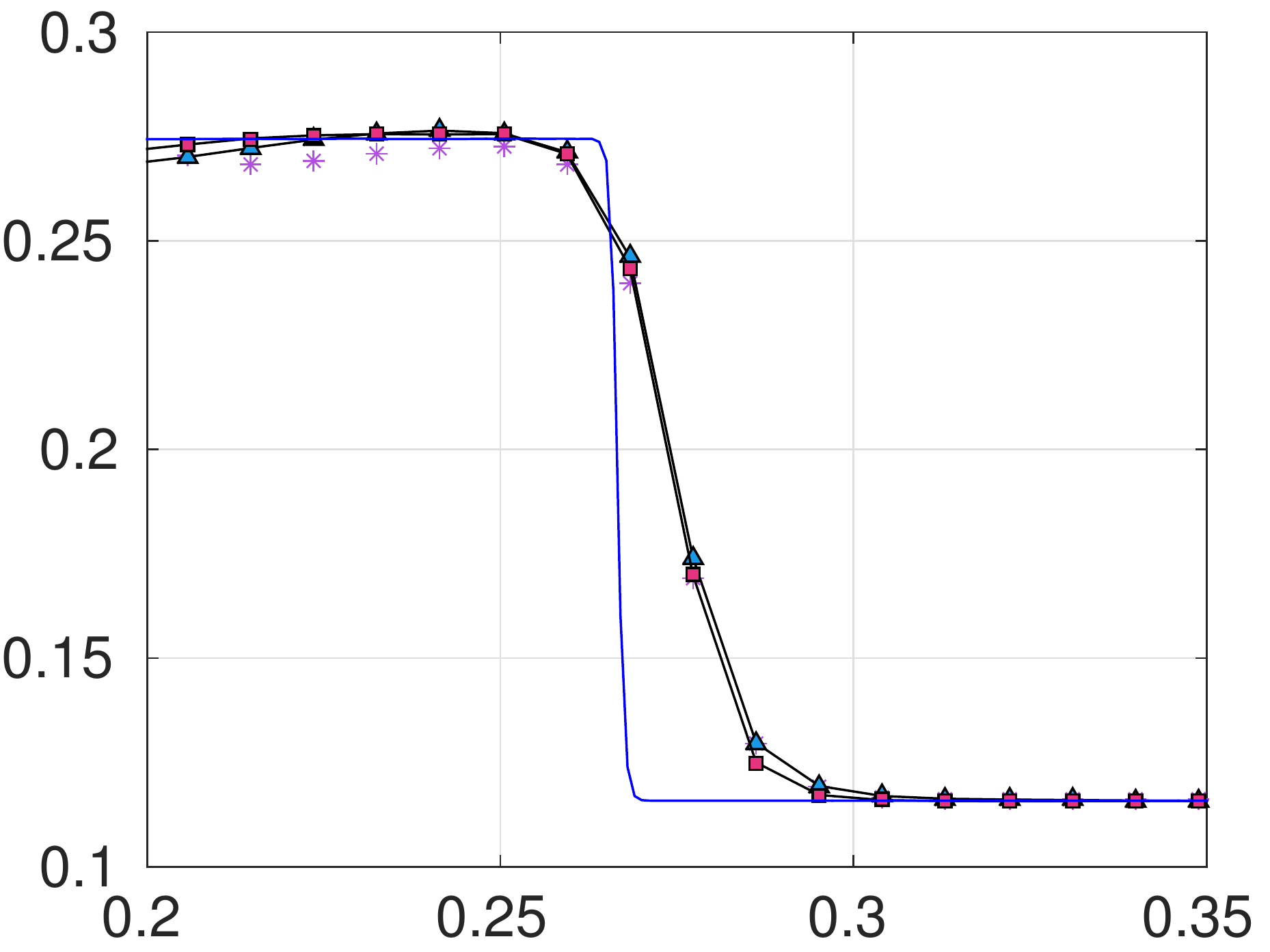}{\figWidthz}};
  \draw(8.4,-.75) node[anchor=south west,xshift=0pt,yshift=+0pt] {\trimfig{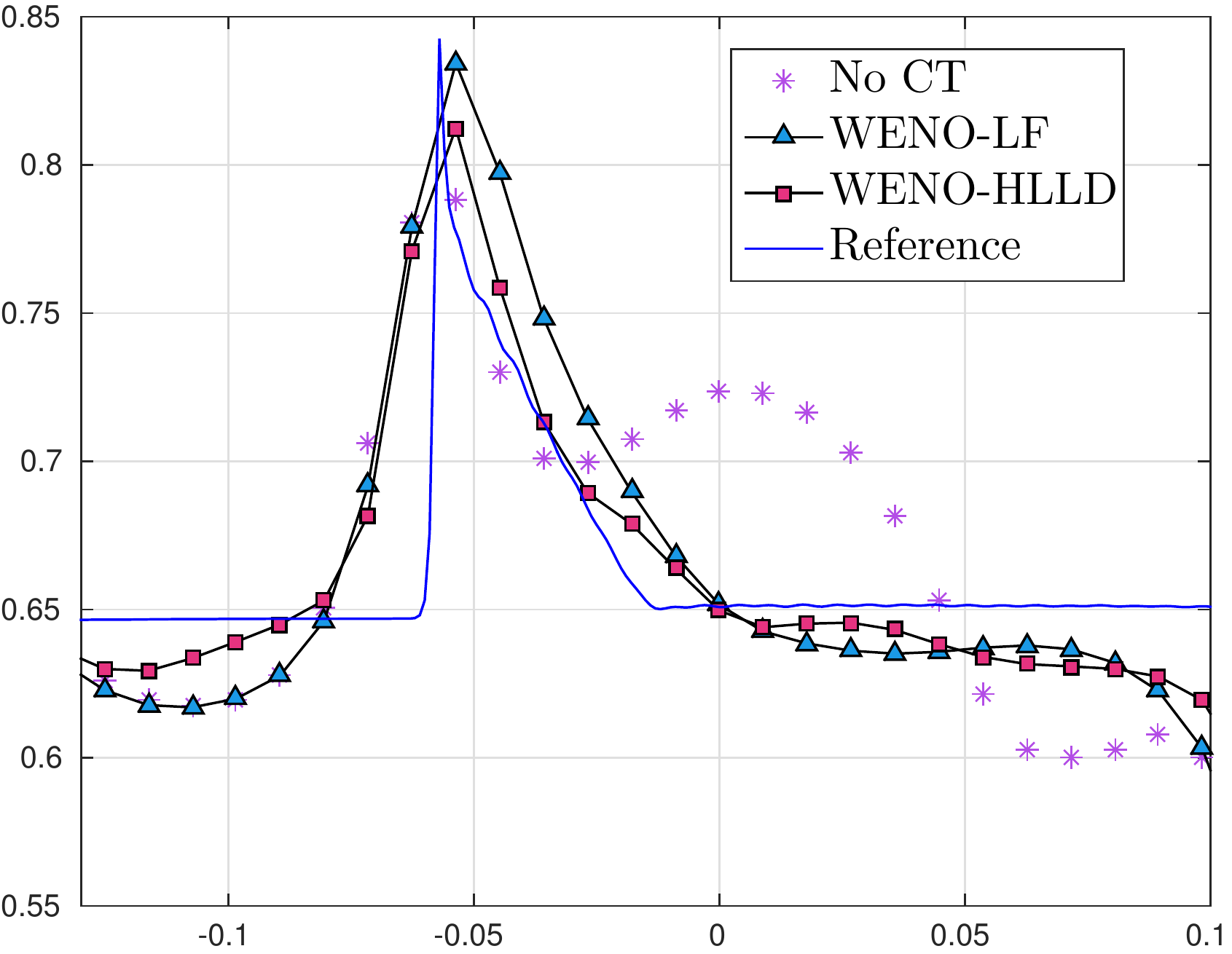}{\figWidth}};
  \draw(.4,.32)   node[blue, anchor=west] {\scriptsize contact discontinuity};
  \draw(5.1,2.7)  node[blue, anchor=west] {\scriptsize shock};
  \draw(10.5,1.4) node[blue, anchor=west] {\scriptsize compound waves};
  \end{scope}
\end{tikzpicture}
} 
\end{center}
    \caption{
2D rotated Brio-Wu shock tube.
Density solved using the WENO schemes at $t = 0.2$ and its three zoomed views around the contact discontinuity, shock and compound waves.
The solutions are computed on a uniform mesh of size $200 \times 100$.   
The reference solution is a numerical solution on a fine mesh of 2000 points.
}
\label{fig:2DRotatedBrioWuCurves}
\end{figure}
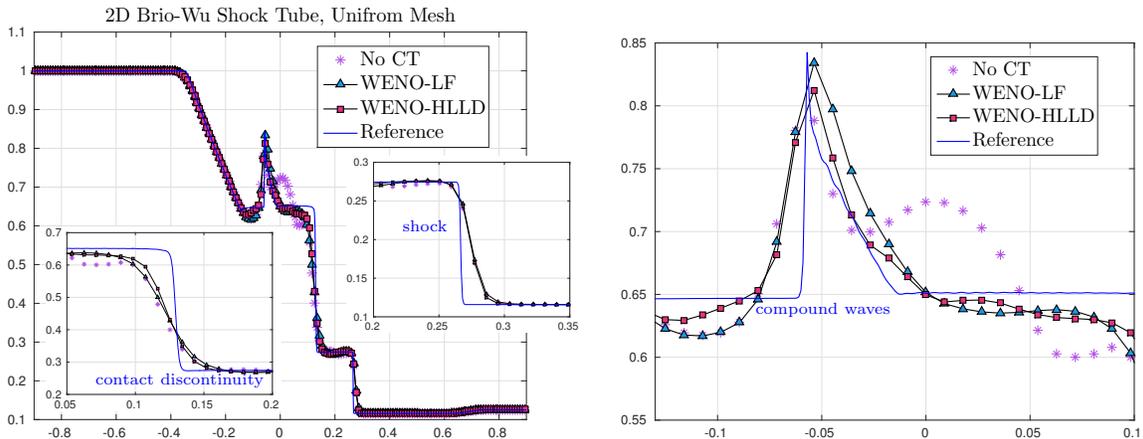
}

The Brio-Wu shock tube is then solved in 2D.  
As illustrated in Figure~\ref{fig:2DRotatedBrioWuContour}, the Riemann problem is rotated by an angle of $\tan^{-1}(0.5)$ with respect to $x$-direction.
The computational domain is $[-1, 1]\times [-0.5, 0.5]$, an inflow boundary condition is used on the left and an outflow boundary condition is used on the right.
The top and bottom boundary conditions are zero-order extrapolations along the tangential direction of the wave propagation for the conserved quantities
and a linear extrapolation along the same direction is used for the potential.
The results presented in this section use a uniform mesh of size $200 \times 100$.  
Since the direction of the wave propagation is not parallel to the coordinate axes any more,
this 2D Riemann problem requires the divergence-free condition to be handled properly.
Figure~\ref{fig:2DRotatedBrioWuContour} presents the contour plot of density solved using the WENO schemes with the Lax-Friedrichs flux. 
Note that there are spurious oscillations around the region of compound waves and contact discontinuity in the solutions using the scheme with the constrained transport step turned off.
For a low-dissipative scheme,  it is found that controlling the divergence error is even more important, 
since the scheme with the HLLD flux becomes unstable before $t$ reaches the final time when the constrained transport step is turned off.
Figure~\ref{fig:2DRotatedBrioWuCurves} shows the computed solutions along $y=0$ that are projected to the direction of the wave propagation.
The solution without the constrained transport step has some oscillations which are not observed in the solutions with the constrained transport step.
The solutions of two schemes with the constrained transport step show a good agreement with the reference solution.
It is observed that the solution of the HLLD flux is still better than the solution of the Lax-Friedrichs flux around the component waves and contact discontinuity,
while the improvement are not so obvious around the shock region.
The plausible reason is that the constrained transport step introduces extra dissipations around the shock.
To reduce the dissipations from the constrained transport step in those shock regions can be a challenging task, since the dissipations may be important to maintain the stability of the schemes. 

%% file: texSISC/fieldloop.tex
\subsection{2D field loop}
\label{sec:2DFieldLoop}
	
{
	
In this section we consider a 2D advection of a weakly magnetized field loop 
from~\cite{gardiner2005unsplit}. 
The initial conditions are	
\begin{align*}
	\left(\rho, u, v, w, p \right) (0,x,y)=\left( 1, \sqrt{5}\cos(\theta), \sqrt{5}\sin(\theta), 0, 1\right)
\end{align*}
with the advection angle of $\theta=\tan^{-1}(0.5)$. Magnetic field components are initialized by taking the curl of the magnetic potential $A$
\begin{align*}
	A (0,x,y) = \left\{\begin{array}{ll}
	0.001 (R-r), & \text{if} \ r\leq R,\\
	0, & \text{otherwise}.\\
	\end{array}
	\right.
\end{align*}
with $r=\sqrt{x^2+y^2}$ and $R=0.3$.
The example is solved on a stationary curve grid and a randomized grid.
The curve grid is mapped from the computational domain $(\xi, \eta)\in[-1,1]\times[-0.5,0.5]$:
\begin{align*}
x=& \xi  +  \epsilon_{x} \sin( 2\pi \eta ),\\
y=& \eta + \epsilon_{y} \sin( 2\pi \xi  ).
\end{align*}
where $\epsilon_{x}=-0.03$ and $\epsilon_{y}=-0.05$.
The randomized grid is formed by randomizing the uniform computation domain $(\xi,\eta)$ with 10\% magnitude grid spacing $\Delta \xi$ or $\Delta\eta$ in a random direction.
Periodic boundary conditions are used in both directions.  
In Figure~\ref{fig:fieldloop}, we present gray-scale images of  $B_1^2+B_2^2$ and contour plots of the potential $A$ at $t=2$, with $200\times 100$ grid points. 
Note that the magnetic field maintains the circular symmetry of the loop as expected. 
The numerical dissipations are observed around the center and edge of the loop, which is similar to the results in~\cite{gardiner2005unsplit}.
Note that for this problem the schemes of different numerical fluxes produce almost identical results, since the solutions are essentially determined by the constrained transport step.
Therefore, the results of the Lax-Friedrichs flux are not presented in Figure~\ref{fig:fieldloop}.

{
\newcommand{\figWidth}{8.5cm}
\newcommand{\figWidthb}{8.3cm}
\newcommand{\trimfig}[2]{\trimFig{#1}{#2}{.77}{.0}{.35}{.9}}
\newcommand{\trimfigb}[2]{\trimFig{#1}{#2}{.78}{.0}{.355}{.9}}
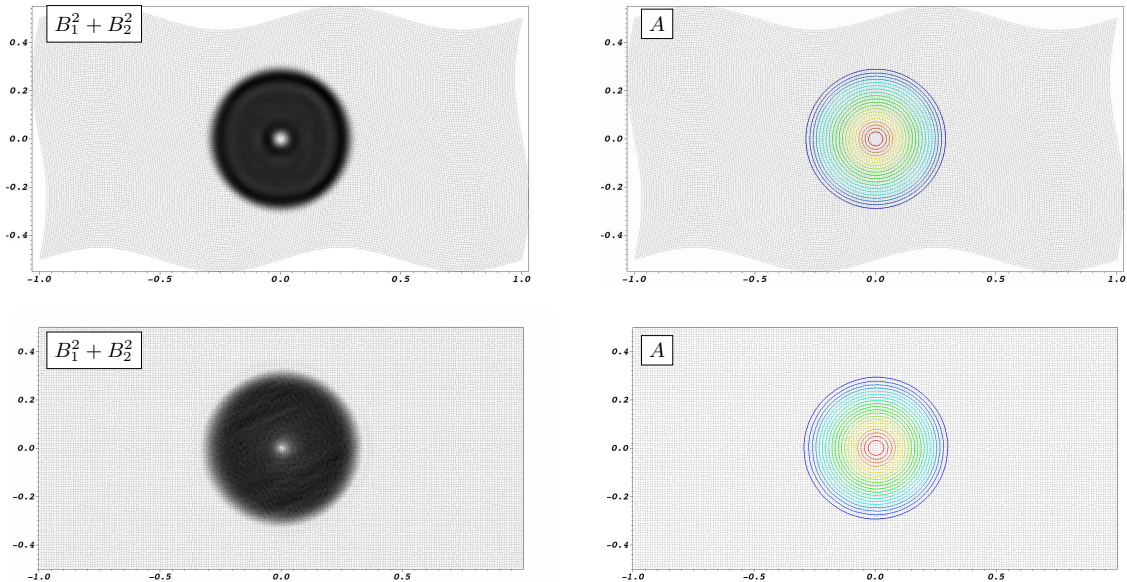
\begin{figure}[htb]
\begin{center}
\resizebox{14cm}{!}{
\begin{tikzpicture}[scale=1]
  \useasboundingbox (0.0,0) rectangle (16,8.6);  
\begin{scope}[yshift=.4cm]
  \draw(-0.9,-.75) node[anchor=south west,xshift=0pt,yshift=+0pt] {\trimfigb{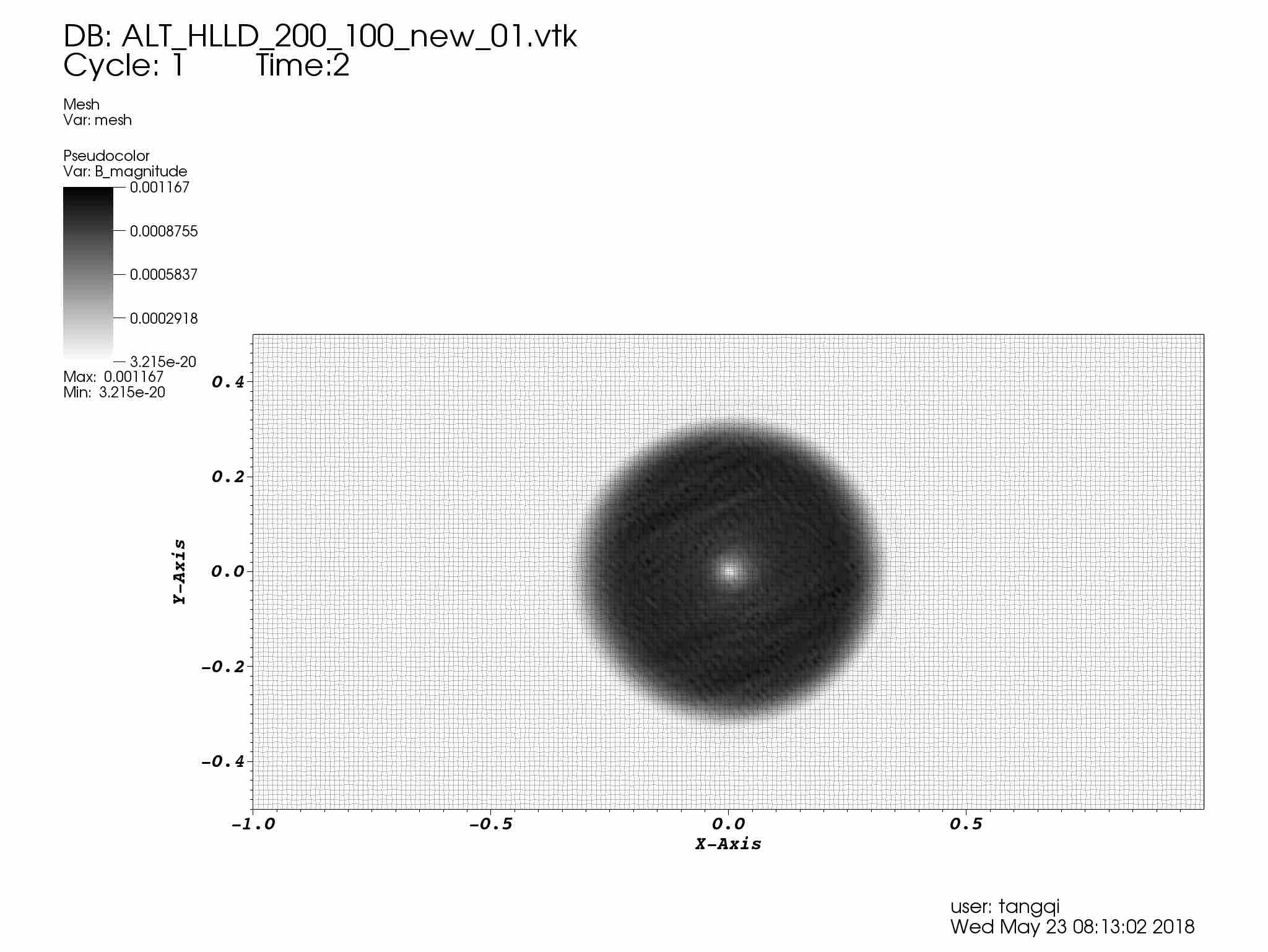}{\figWidthb}};
  \draw(8.2, -.75) node[anchor=south west,xshift=0pt,yshift=+0pt] {\trimfigb{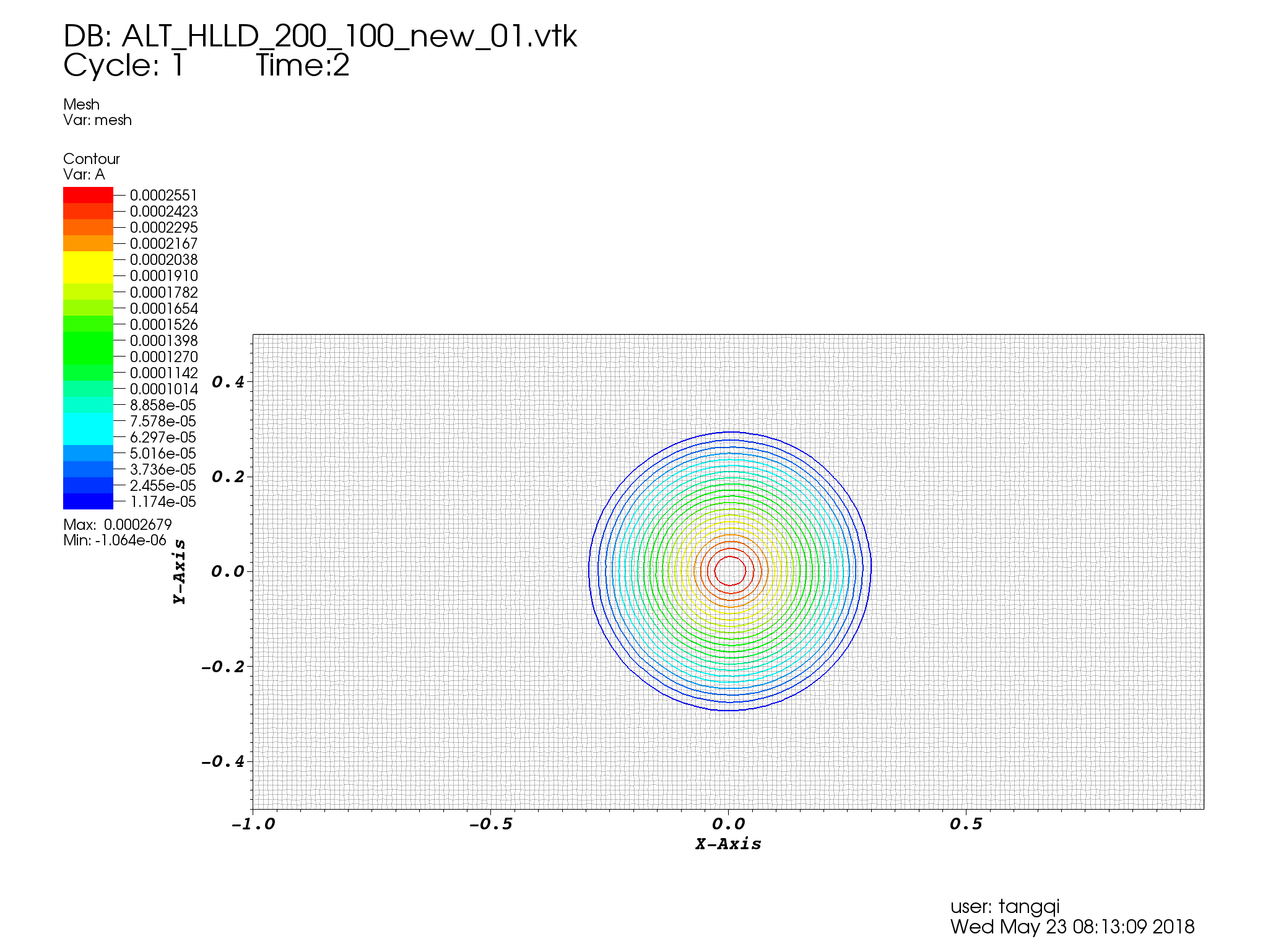}{\figWidthb}};
  \draw(-1,3.8) node[anchor=south west,xshift=0pt,yshift=+0pt] {\trimfig{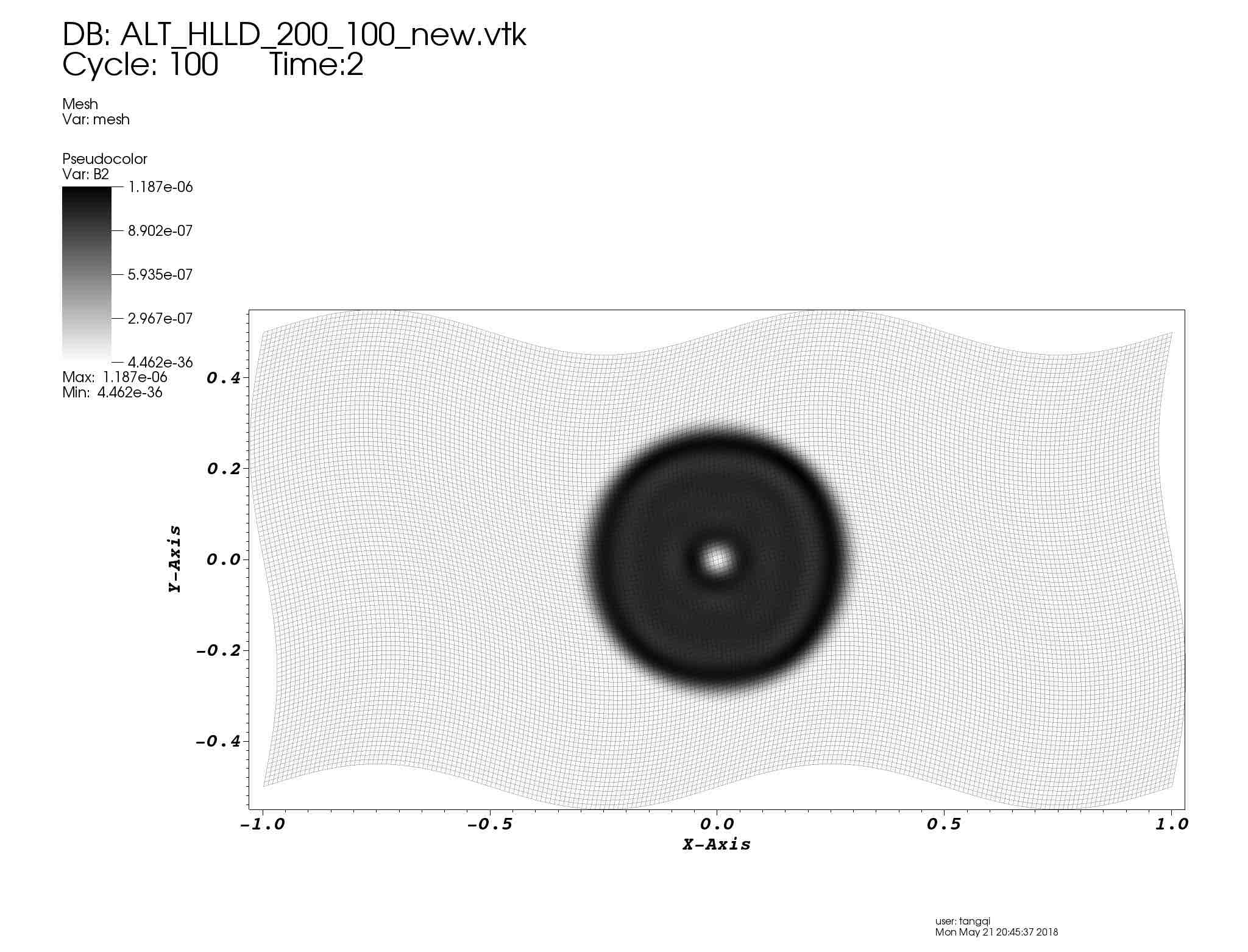}{\figWidth}};
  \draw(8.1, 3.8) node[anchor=south west,xshift=0pt,yshift=+0pt] {\trimfig{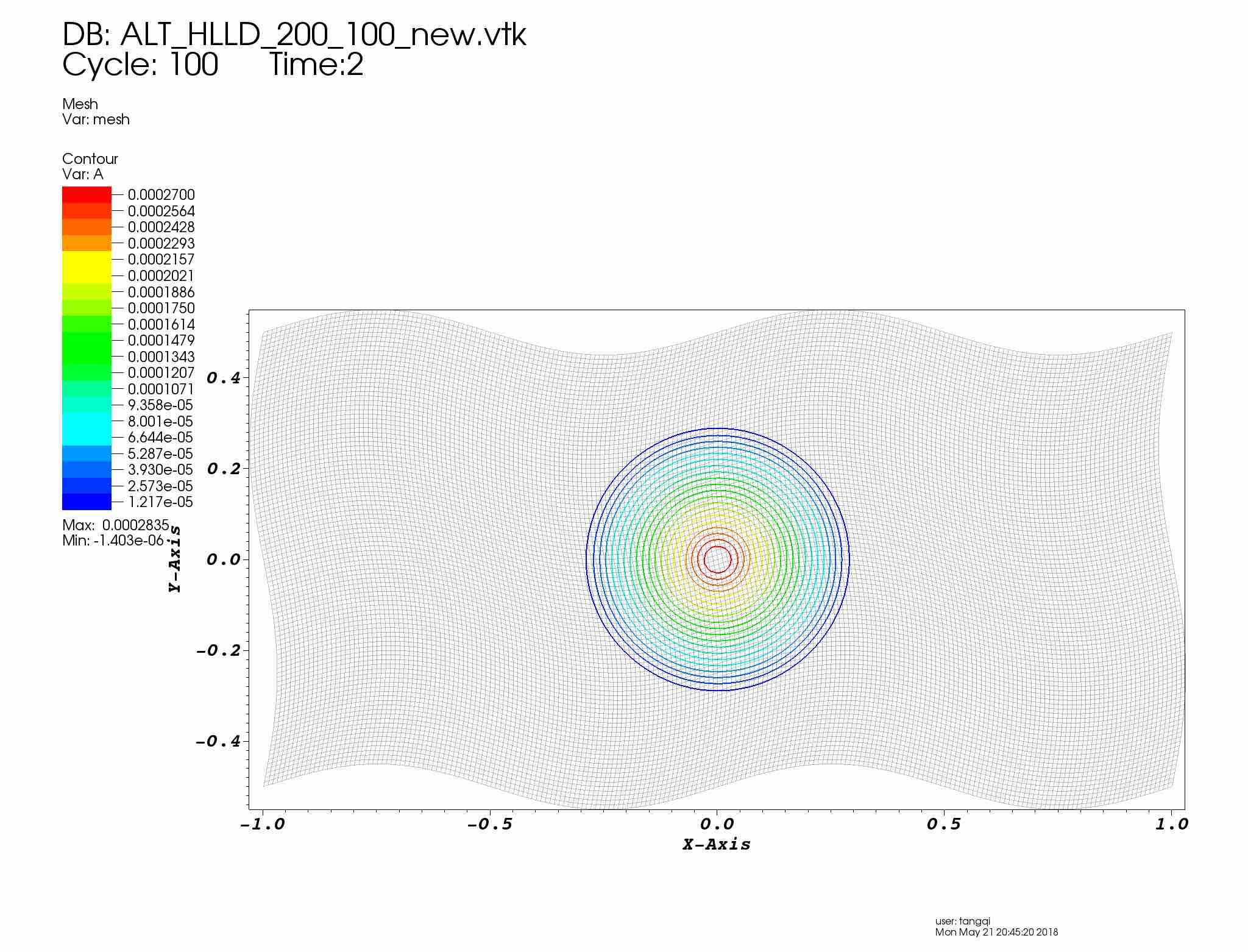}{\figWidth}};
  \draw(-.3,3.) node[draw,fill=white,anchor=west,xshift=2pt,yshift=-1pt] {\small $B_1^2+B_2^2$};
  \draw(8.8,3.) node[draw,fill=white,anchor=west,xshift=2pt,yshift=-1pt] {\small $A$};
  \draw(-.3,8.) node[draw,fill=white,anchor=west,xshift=2pt,yshift=-1pt] {\small $B_1^2+B_2^2$};
  \draw(8.8,8.) node[draw,fill=white,anchor=west,xshift=2pt,yshift=-1pt] {\small $A$};
  \end{scope}
 grid:
\end{tikzpicture}
} 
\end{center}
    \caption{
2D field loop. Top row: the solutions on the curve grid at t = 2. 
Bottom row:  the solutions on the randomized grid at t = 2.
The computational grids of size $200\times100$ are plotted (the light solid lines).
}
\label{fig:fieldloop}
\end{figure}
}

%
}

%% file: texSISC/orszagTang.tex
  \subsection{2D Orszag-Tang vortex}
  \label{sec:orszagTang}
  We next consider a common benchmark problem of the 2D Orszag-Tang vortex problem.  
  The initial conditions are
  \begin{equation*}
      (\rho, u, v, w, p, B_1, B_2, B_3)(0, x, y)  = (\gamma^2, -\sin (y), \sin (x),0, \gamma, -\sin (y), \sin (2x), 0),
  \end{equation*}
  and its initial magnetic potential is 
  \begin{equation*}
    A(0, x, y) = 0.5\cos(2x) + \cos(y).
  \end{equation*}
  To examine the performance of the schemes on general curvilinear meshes, 
  a mesh similar to the one in Section~\ref{sec:alfvenWave} is used.
  In particular, the computational domain is $(\xi, \eta) \in [0, 2\pi] \times [0, 2\pi]$ 
  and the curvilinear grid is given by the mapping
  \begin{equation*}
    \begin{aligned}
      x & = \xi + \epsilon_x \sin (\eta \, a_x), \\
      y & = \eta + \epsilon_y \sin (\xi \, a_y),
    \end{aligned}
  \end{equation*}
  where $\epsilon_x = 0.03$, $\epsilon_y = 0.05$, $a_x = 2$ and
  $a_y = 4$. The boundary conditions are all periodic. 
  
{
\newcommand{\figWidth}{7cm}
\newcommand{\figWidthb}{6.95cm}
\newcommand{\trimfig}[2]{\trimFig{#1}{#2}{.58}{.05}{.37}{.0}}
\newcommand{\trimfigb}[2]{\trimFig{#1}{#2}{.55}{.05}{.34}{.0}}
\begin{figure}[htb]
\begin{center}
\resizebox{12cm}{!}{
\begin{tikzpicture}[scale=1]
  \useasboundingbox (0,0) rectangle (16,6.6);  
  \draw(-.2,-.5) node[anchor=south west,xshift=0pt,yshift=+0pt] {\trimfig{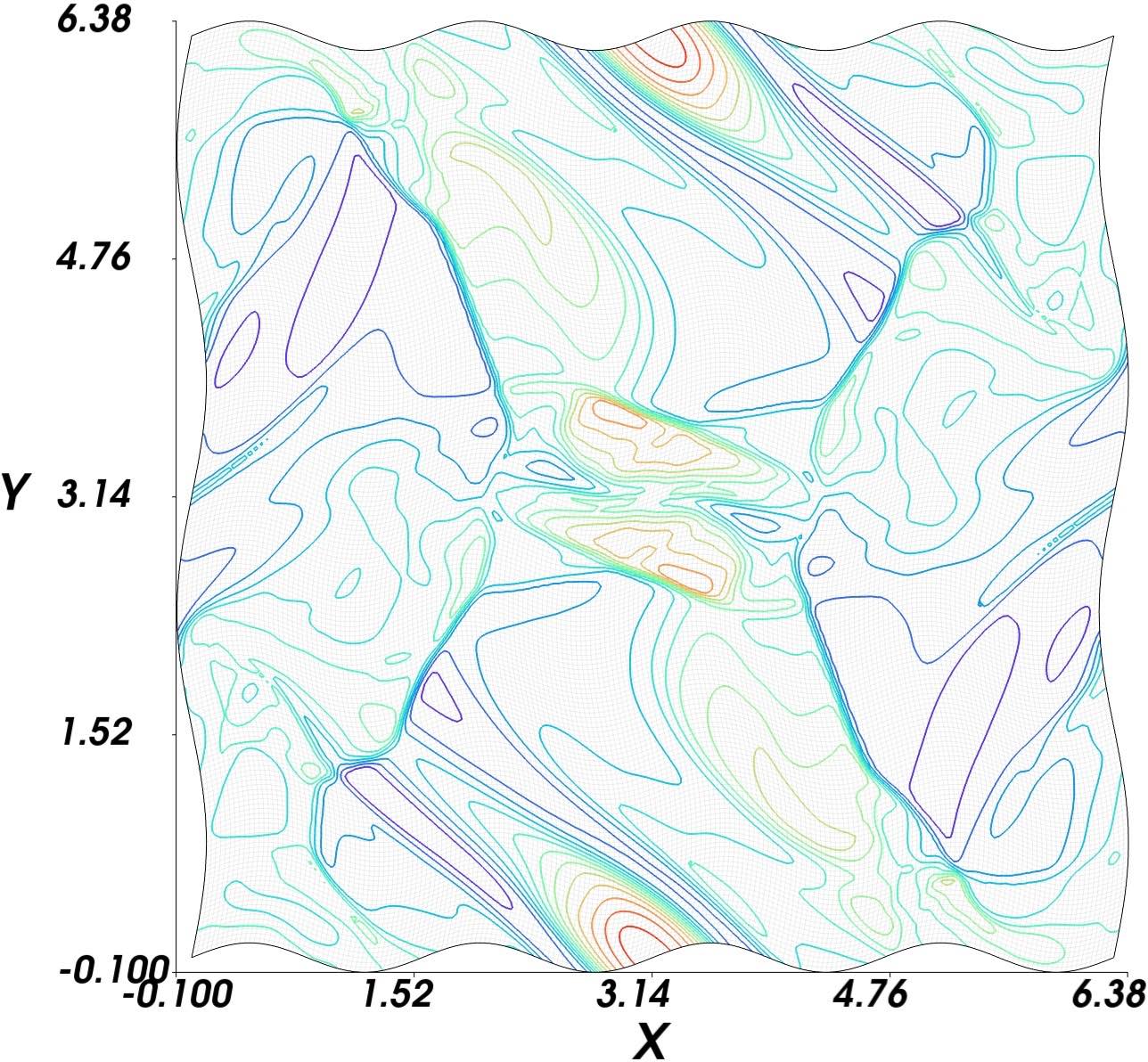}{\figWidth}};
  \draw(9,-.5) node[anchor=south west,xshift=0pt,yshift=+0pt] {\trimfigb{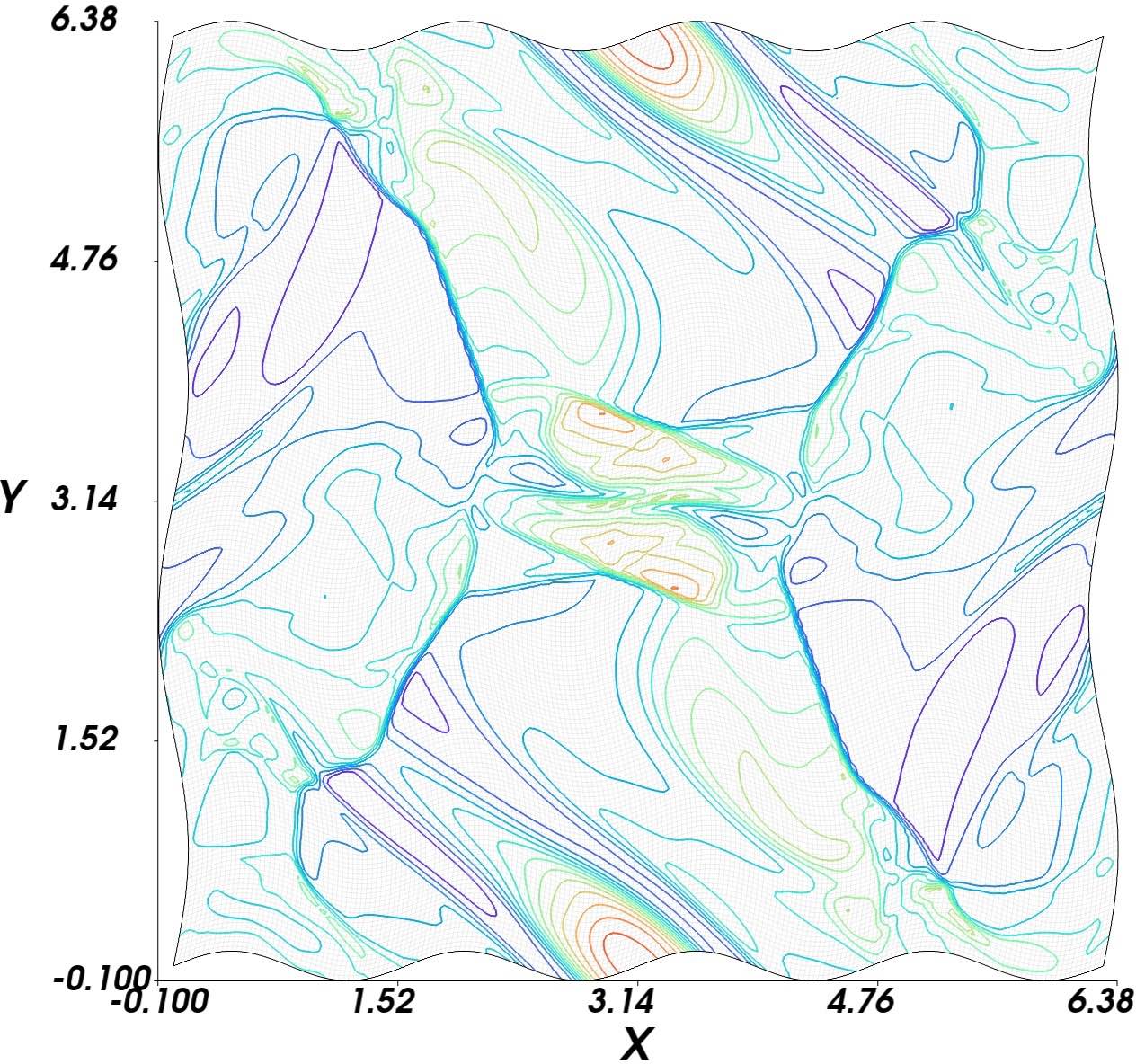}{\figWidthb}};
%
\end{tikzpicture}
}
\end{center}
    \caption{Orszag-Tang vortex problem.
Contour plots of density at $t=3$ are presented with 15 equally spaced contour lines.  A perturbed mesh of size $192 \times 192$ is used.
Left: WENO with Lax-Friedrichs flux. Right: WENO with HLLD flux.
The constrained transport step and positivity-preserving limiter are turned on.}
\label{fig:orszagTang}
\end{figure}
}
  Figure~\ref{fig:orszagTang} presents the contour plots of the density at time $t=3$. 
  The solutions of two schemes are presented with 15 equally spaced contour lines.
  The problem starts from an smooth initial condition and develops an vortex and several MHD shock waves.
  Those waves interacts with each other (see Figure~\ref{fig:orszagTang}) and eventually result in turbulence.
  For such a problem, a low-dissipative scheme is preferred to capture small structures.
  It is observed from Figure~\ref{fig:orszagTang} that the HLLD flux produces the less dissipative solutions, for instance, around the shock region.
  Although the previous Riemann problem test in Section~\ref{sec:shockTube} shows the dissipation from the constrained transport problem may smear those shocks,
  the current results show that the improvement of a low-disspative scheme is still significant for practical problems, such as MHD turbulence simulations.
  The computed solutions match well with those found in the literature~\cite{Christlieb2014,Dai1998a,Rossmanith2006,Toth2000,Zachary1994}.
We note that the simulations run successfully to a much later time of $t=10$,
which indicates that the divergence-free condition is handled 
properly on a curvilinear mesh by the constrained transport approach.
Without the constrained transport step, the simulations becomes 
unstable as soon as discontinuities develop in the solutions.

%

%% file: texSISC/cloudShock.tex
  \subsection{2D cloud-shock interaction}

  \label{sec:cloudShock}

  In this section we consider the 2D cloud-shock interaction problem.
  The initial conditions are
  \begin{align*}
    &(\rho, u, v, w, p, B_1, B_2, B_3) (0,x,y) \\ \nonumber
    &\quad =
      \begin{cases}
        (3.86859, 11.2536, 0, 0, 167.345, 0, 2.1826182, -2.1826182) &
        \text{if $x < 0.05$}, \\
        (10, 0, 0, 0, 1, 0, 0.56418958, 0.56418958) &
        \text{if ${x > 0.05}$, ${r < 0.15}$,} \\
        (1, 0, 0, 0, 1, 0, 0.56418958, 0.56418958) & \text{otherwise,}
      \end{cases}
  \end{align*}
  where $r = \sqrt{{(x-0.25)}^2 + {(y-0.5)}^2}$.  The initial magnetic
  potential is,
  \begin{equation*}
    A (0,x,y)=
    \begin{cases}
      -2.1826182\, x + 0.080921431 &
      \text{if $x \leq 0.05$}, \\
      -0.56418958 \,x & \text{if $x > 0.05$}.
    \end{cases}
  \end{equation*}
The problem models an MHD shock propagating toward a dense bubble, 
resulting into very complex structures as the shock passes through the bubble.
Those structures around the bubble are very sensitive to the numerical dissipations 
and  low dissipative schemes are advantageous to obtain less smeared structures.
Here we use this problem to study the effects of the Riemann solvers on both the Cartesian and curvilinear grids.

{
\newcommand{\figWidth}{5cm}
\newcommand{\trimfig}[2]{\trimFig{#1}{#2}{.0}{.0}{.0}{.0}}
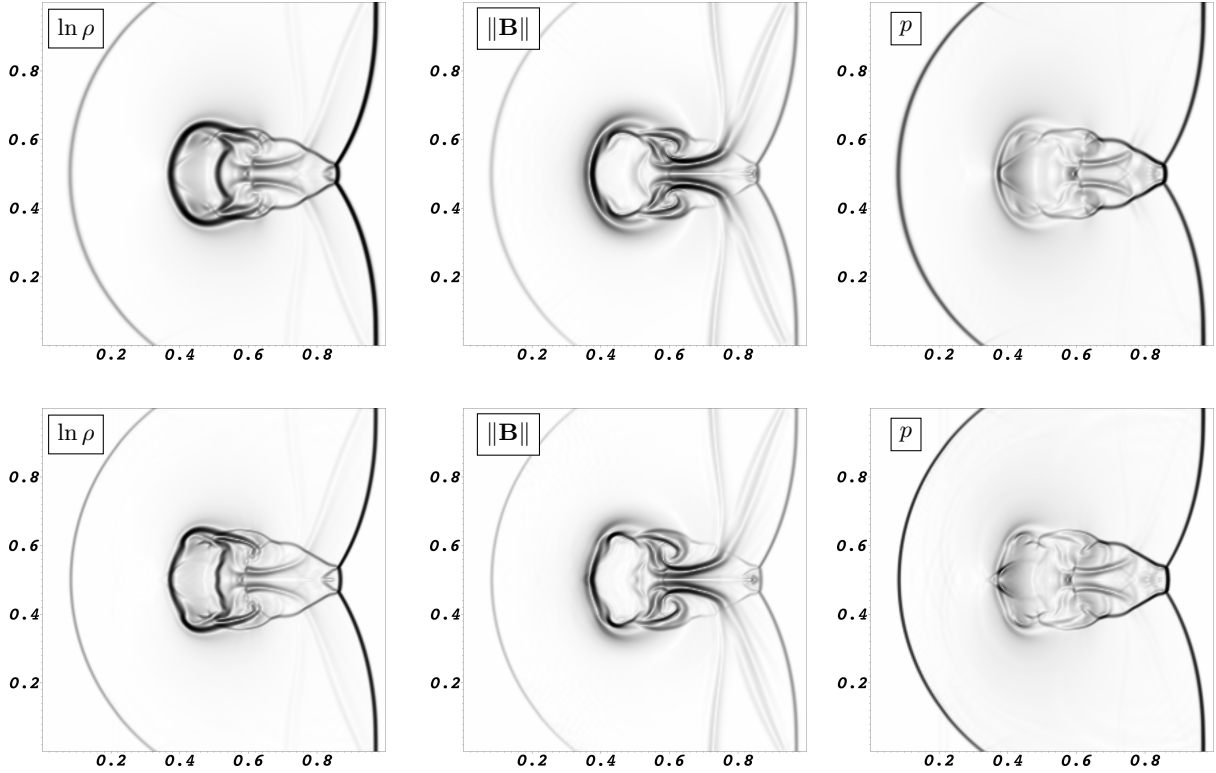
\begin{figure}[htb]
\begin{center}
\begin{tikzpicture}[scale=1]
  \useasboundingbox (0,0) rectangle (15,10);  
  \draw(-0.8,5.) node[anchor=south west,xshift=0pt,yshift=+0pt] {\trimfig{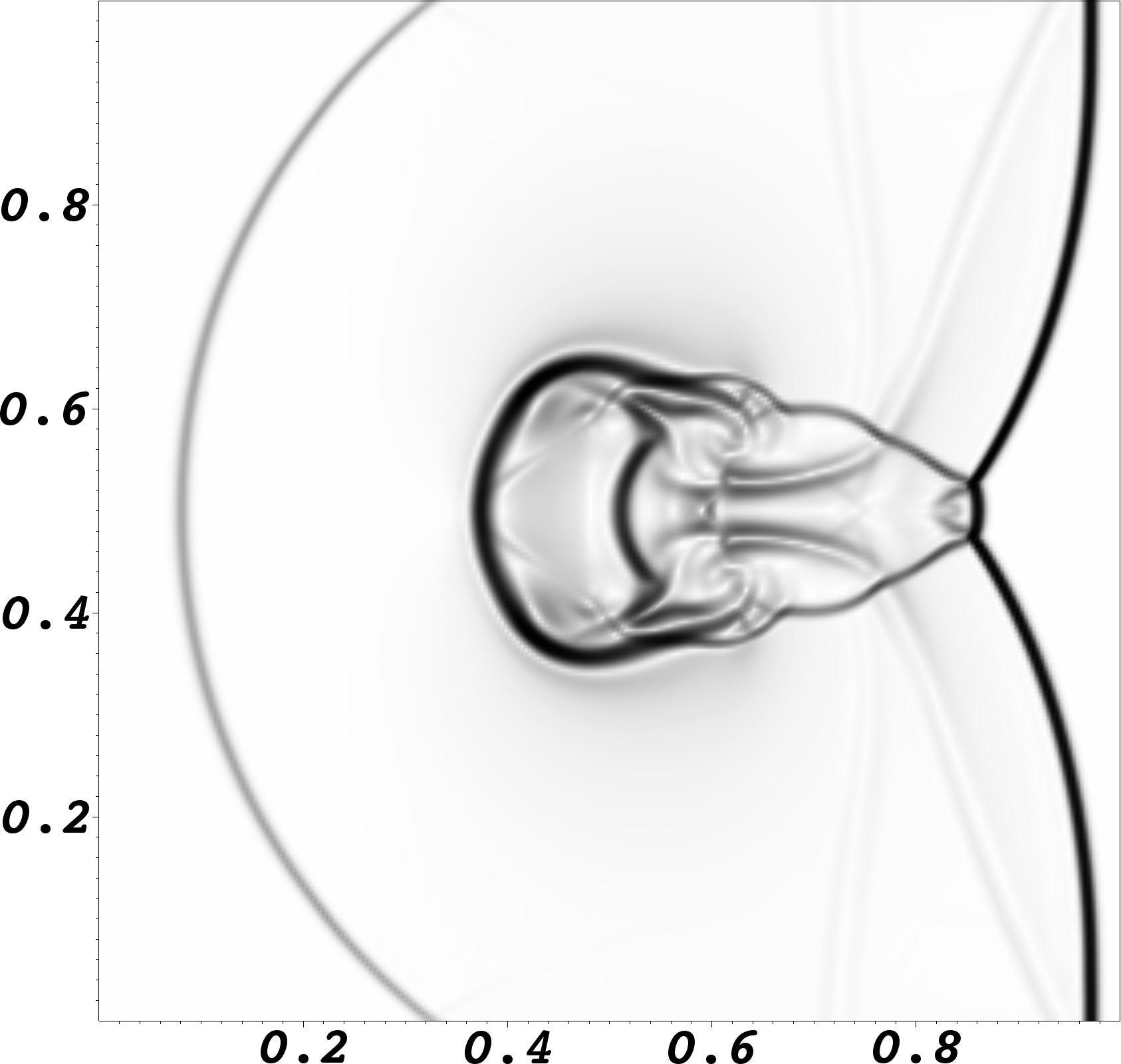}{\figWidth}};
  \draw(4.8 ,5.) node[anchor=south west,xshift=0pt,yshift=+0pt] {\trimfig{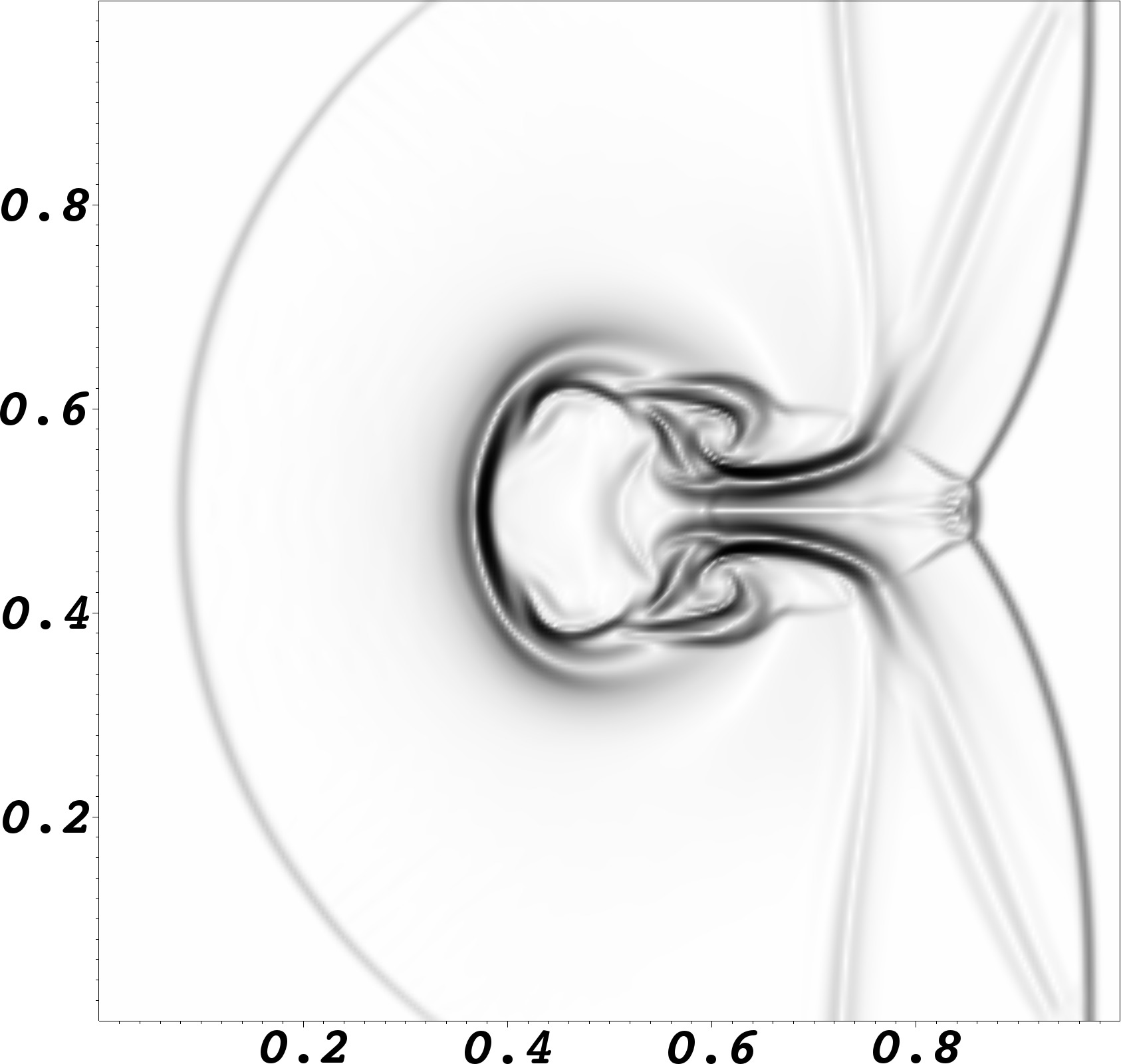}{\figWidth}};
  \draw(10.2,5.) node[anchor=south west,xshift=0pt,yshift=+0pt] {\trimfig{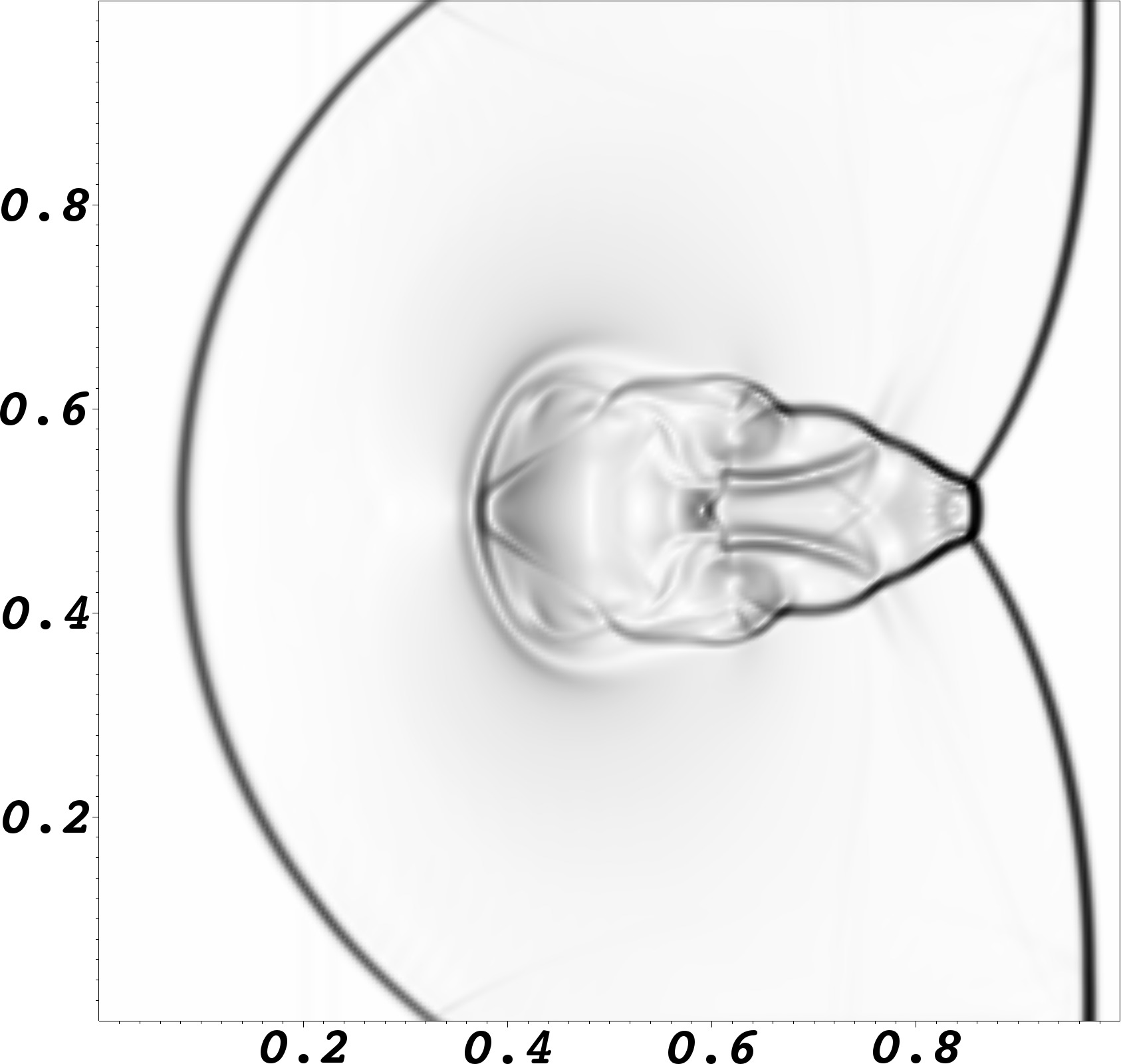}{\figWidth}};
  \draw(-.2,9.6) node[draw,fill=white,anchor=west,xshift=2pt,yshift=-1pt] {\small $\ln \rho$};
  \draw(5.5,9.6) node[draw,fill=white,anchor=west,xshift=2pt,yshift=-1pt] {\small $\|\Bv\|$};
  \draw(11,9.6) node[draw,fill=white,anchor=west,xshift=2pt,yshift=-1pt] {\small $p$};
  \draw(-.8 ,-.4) node[anchor=south west,xshift=0pt,yshift=+0pt] {\trimfig{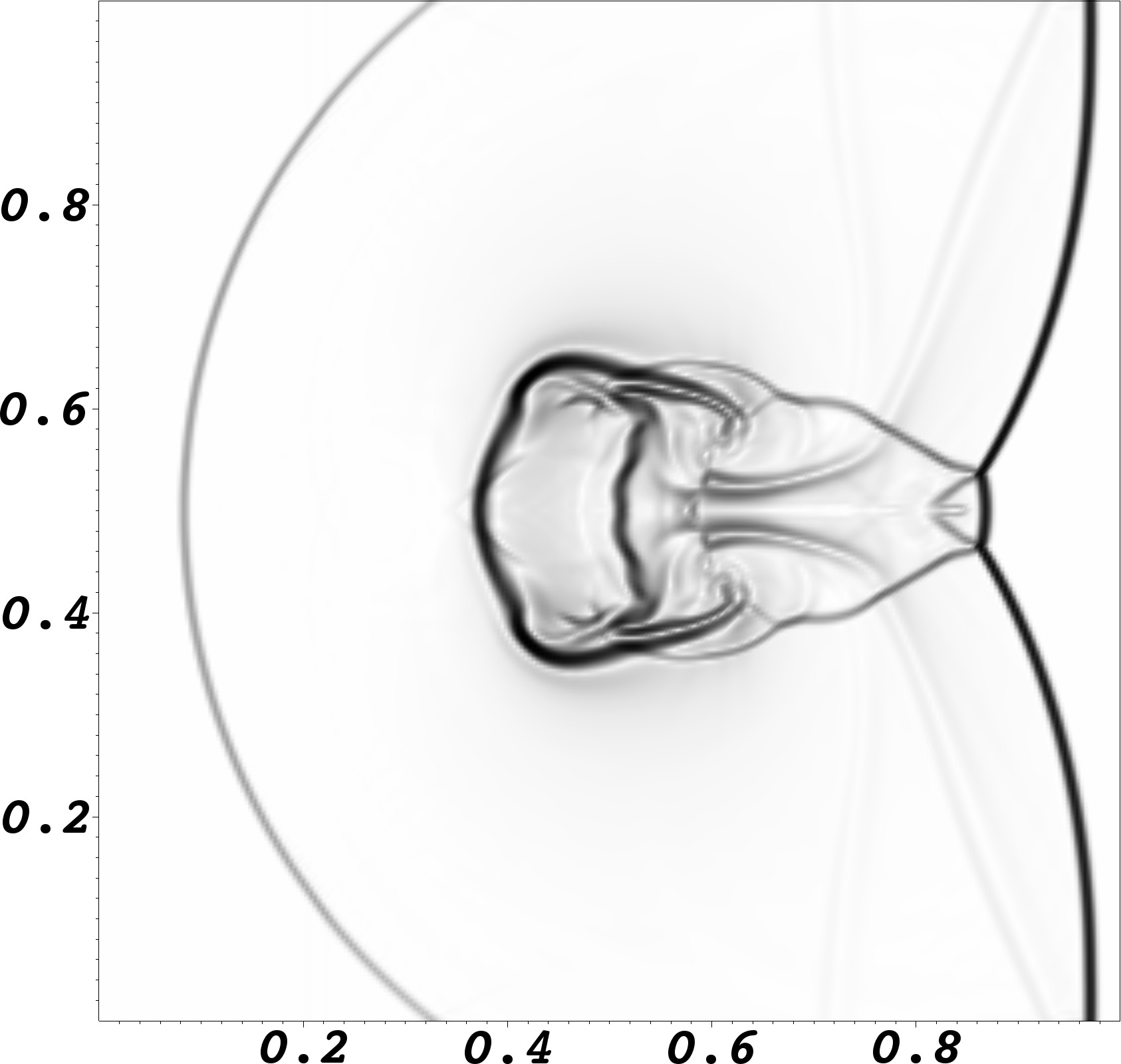}{\figWidth}};
  \draw(4.8 ,-.4) node[anchor=south west,xshift=0pt,yshift=+0pt] {\trimfig{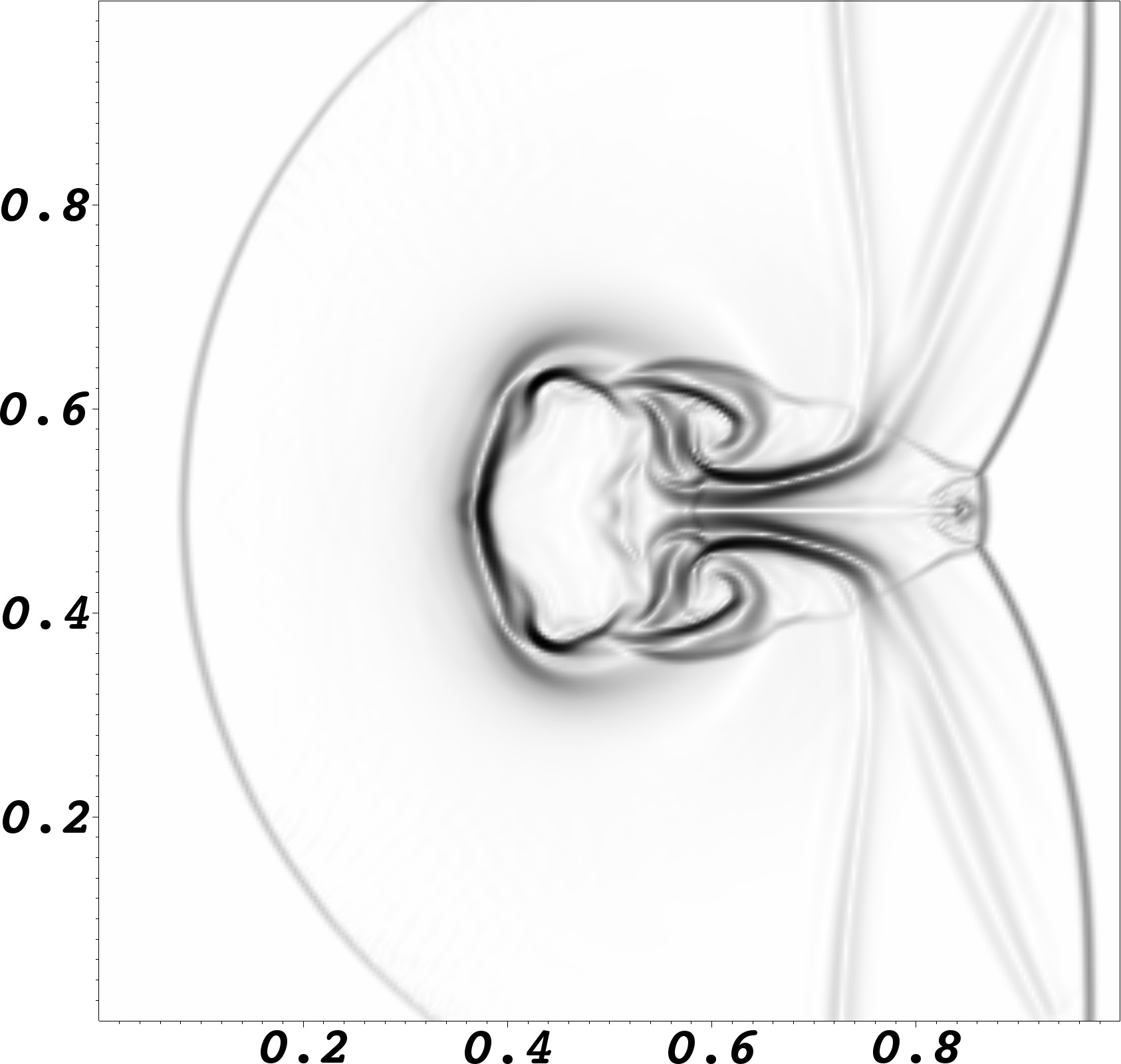}{\figWidth}};
  \draw(10.2,-.4) node[anchor=south west,xshift=0pt,yshift=+0pt] {\trimfig{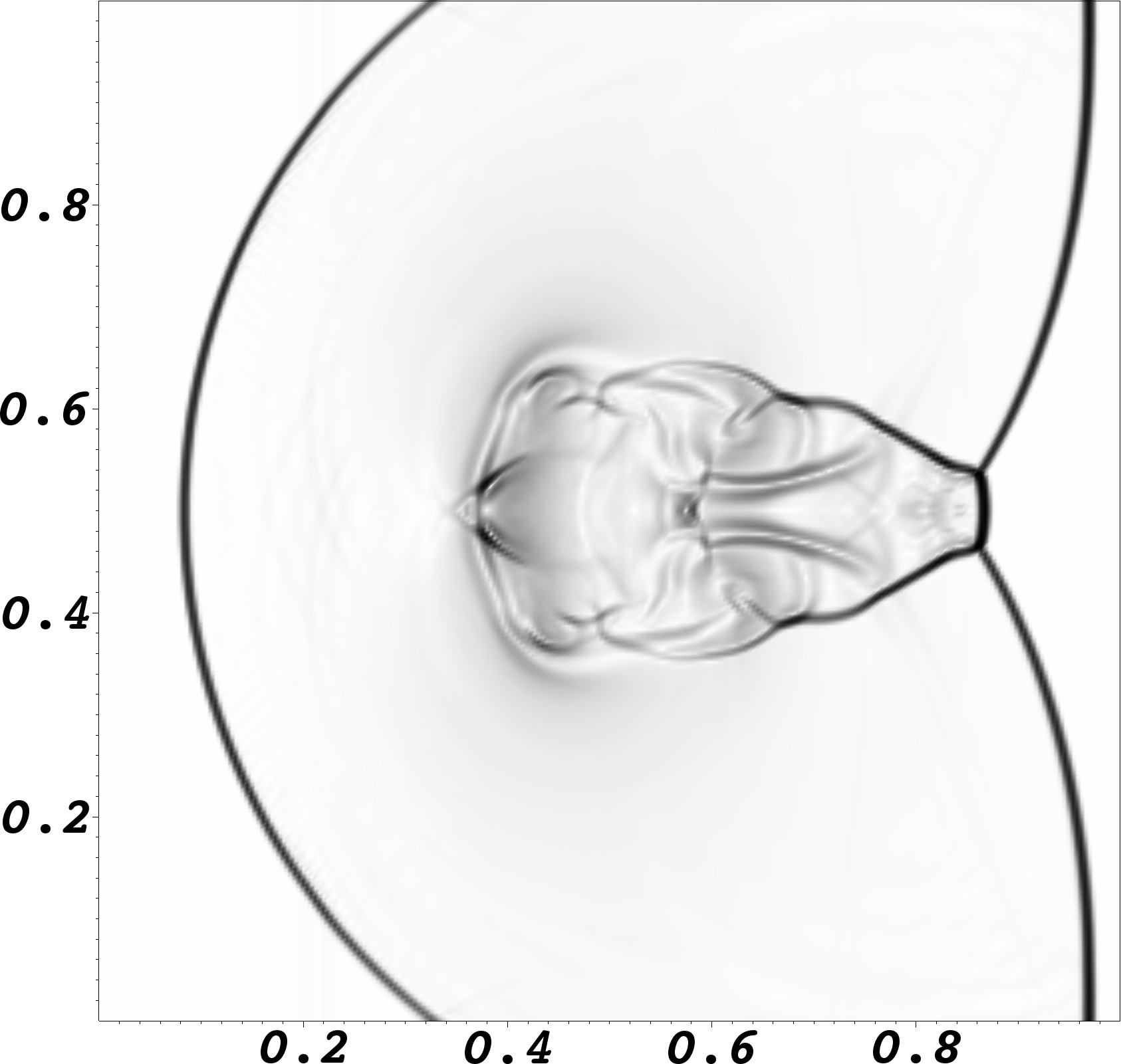}{\figWidth}};
  \draw(-.2,4.2) node[draw,fill=white,anchor=west,xshift=2pt,yshift=-1pt] {\small $\ln \rho$};
  \draw(5.5,4.2) node[draw,fill=white,anchor=west,xshift=2pt,yshift=-1pt] {\small $\|\Bv\|$};
  \draw(11.,4.2) node[draw,fill=white,anchor=west,xshift=2pt,yshift=-1pt] {\small $p$};
%
\end{tikzpicture}
\end{center}
    \caption{\label{fig:2DCloudShockAltWENOUniform}
    2D cloud-shock interaction. 
    Schlieren plots of the logarithm of the density, norm of the magnetic field  and pressure  at $t=0.06$.
Top row: WENO with the Lax-Friedrichs flux.
Bottom row: WENO with the HLLD flux.
A uniform Cartesian grid of size $256 \times 256$ is used.  
 The constrained transport step and positivity-preserving limiter are turned on.
}
\end{figure}
}

The problem is first solved on a square domain of $(x,y) \in [0, 1] \times [0, 1]$.
A uniform Cartesian grid of size $256 \times 256$ is used with
an inflow boundary condition applied at the left boundary and
the outflow boundary condition applied at the other three boundaries.
Figure~\ref{fig:2DCloudShockAltWENOUniform} presents Schlieren plots of the logarithm of the density, norm of the magnetic field and pressure at $t=0.06$.
  The solution matches well with the results in the literature, 
  such as those in~\cite{Christlieb2016,Christlieb2014,Dai1998a,Rossmanith2006}.
  It is observed that the HLLD flux resolves shocks and other complex features much better, although both schemes are fourth-order accurate 
  and produce similar results in the smooth~\Alfven~wave test in Section~\ref{sec:alfvenWave}.  
  In particular, the complex structures around the initial bubble locations of the HLLD flux are less smeared than those in the Lax-Friedrichs flux.

{
\newcommand{\figWidth}{4cm}
\newcommand{\figWidthb}{3.3cm}
\newcommand{\trimfig}[2]{\trimFig{#1}{#2}{.0}{.0}{.0}{.0}}
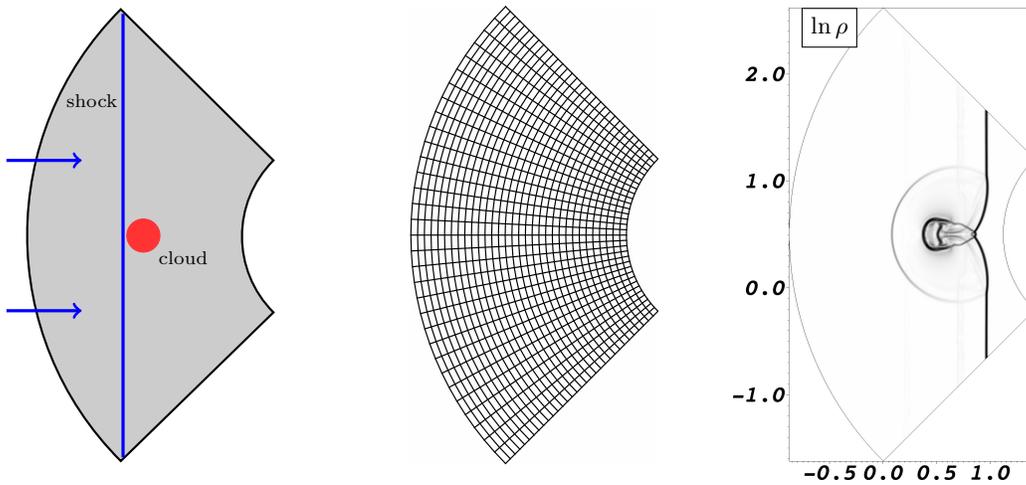
\begin{figure}[htb]
\begin{center}
\begin{tikzpicture}[scale=1]
  \useasboundingbox (0,0) rectangle (12,6);  
  \draw(4.2,-.3) node[anchor=south west,xshift=0pt,yshift=+0pt] {\trimfig{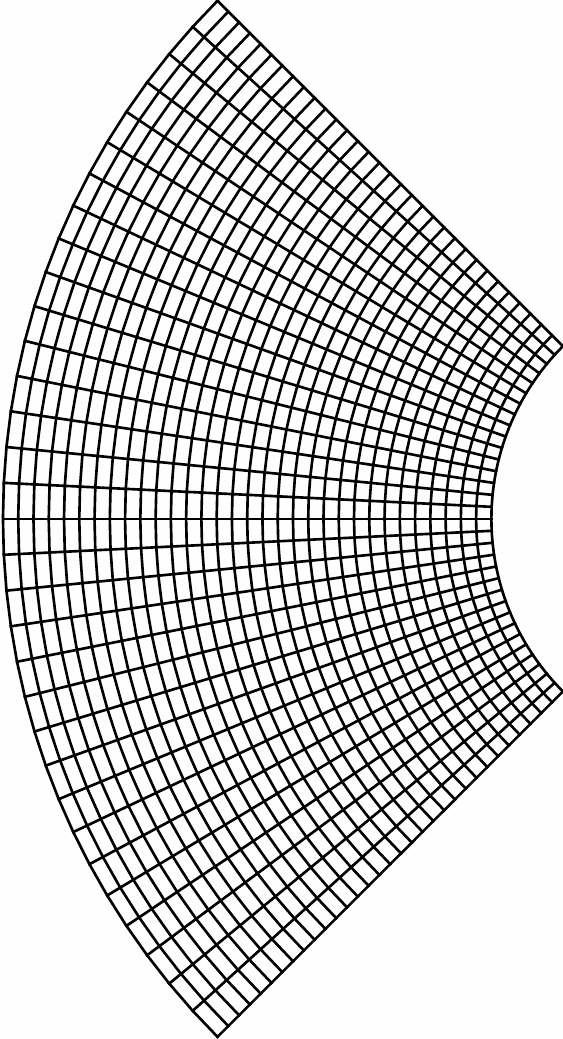}{\figWidthb}};
  \draw(8.5,-.5) node[anchor=south west,xshift=0pt,yshift=+0pt] {\trimfig{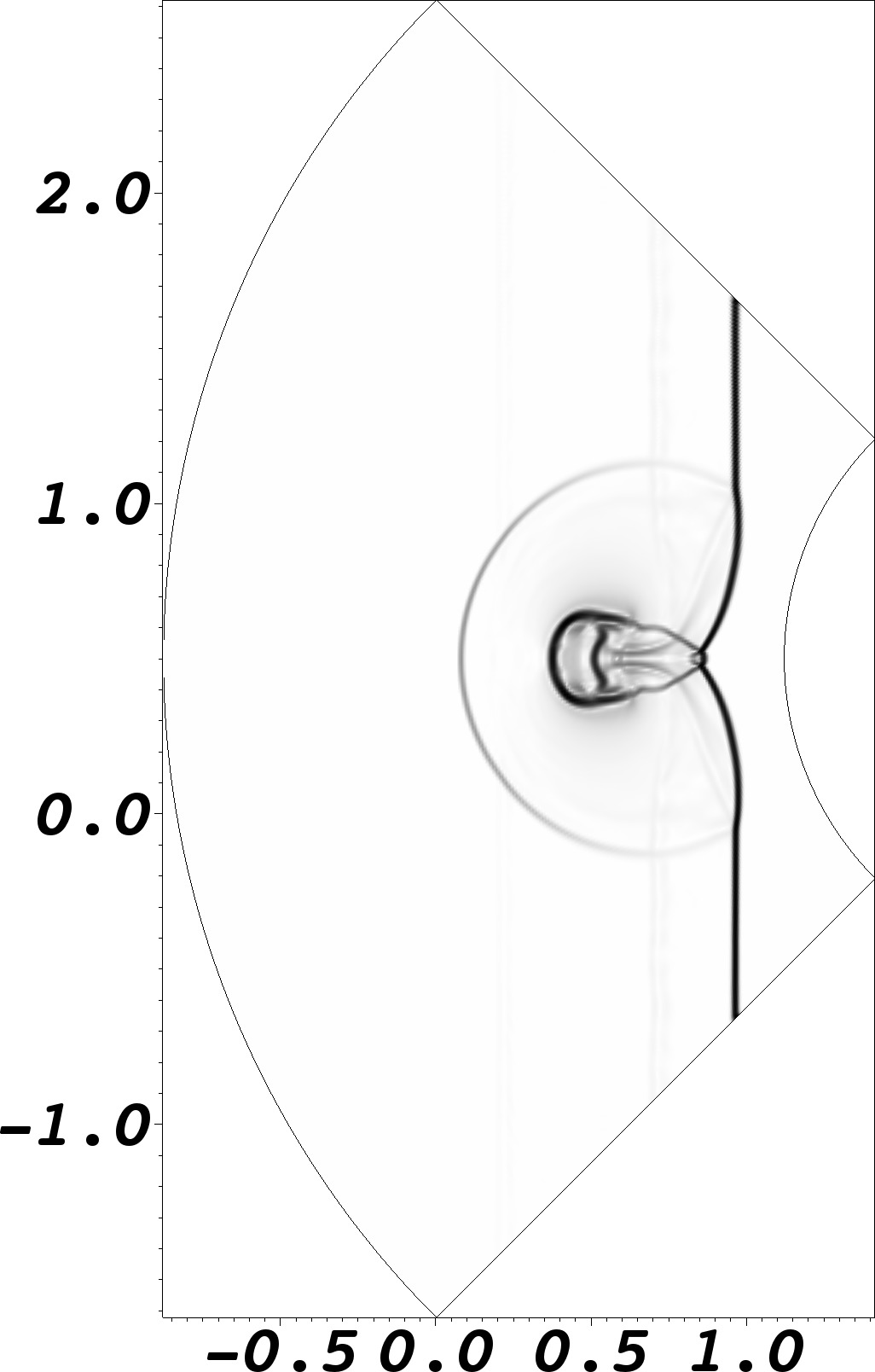}{\figWidth}};
\begin{scope}[xshift=-7.5cm,yshift=-.1cm]
  \filldraw[fill=gray!40,thick] (10.05,4) arc [radius=1.43, start angle=135, delta angle=90]
                    -- (8.02,0) arc [radius=4.25, start angle=225, delta angle=-90]
                    -- cycle;
      \draw[->,very thick,blue] (6.5,4) -- (7.5,4); 
      \draw[->,very thick,blue] (6.5,2) -- (7.5,2);
  \draw[red!80,fill=red!80] (8.32,3) circle (.22);
  \draw[blue,very thick](8.05,0.05) -- (8.05,5.95);
  \draw(8.4,2.7) node[anchor=west] {\scriptsize cloud};
  \draw(8.1,4.8) node[anchor=east] {\scriptsize shock};
  \end{scope}
  \draw(9.5,5.7) node[draw,fill=white,anchor=west,xshift=2pt,yshift=-1pt] {\small $\ln\rho$};
%
\end{tikzpicture}
\end{center}
    \caption{\label{fig:2DCloudShockGrid}
    Left: a diagram of the cloud-shock interaction in a sector domain.
    Middle: a coarse grid of size $32 \times 32$.
    Right: the density at $t = 0.06$ on a fine grid of size $256 \times 256$.
    }
\end{figure}
}

We next consider the same problem but change the physical domain to a sector region.
As illustrated in Figure~\ref{fig:2DCloudShockGrid}, the domain is determined by the following mapping 
  \begin{equation*}
    \begin{aligned}
      x &= (3 - 2 \xi) \cos(\pi + (1-2\,\eta) \pi / 4) + 3 \cos(\pi / 4), \\
      y &= (3 - 2 \xi) \sin(\pi + (1-2\,\eta) \pi / 4) + 0.5,
    \end{aligned}
  \end{equation*}
with $(\xi, \eta) \in [0, 1] \times [0, 1]$.
A uniform grid in the computational domain $(\xi, \eta)$ is used for simulations as illustrated in Figure~\ref{fig:2DCloudShockGrid}.
The same initial condition is used but the initial cloud becomes relatively smaller compared to the computational domain,
which can be also seen in the plot of density in Figure~\ref{fig:2DCloudShockGrid}.
An inflow boundary condition applied at the left edge of the sector domain and the remaining boundary conditions are outflow.

{
\newcommand{\figWidth}{5cm}
\newcommand{\figWidtha}{4.5cm}
\newcommand{\trimfig}[2]{\trimFig{#1}{#2}{.0}{.0}{.0}{.0}}
\newcommand{\trimfiga}[2]{\trimFig{#1}{#2}{.8}{.0}{.34}{.0}}
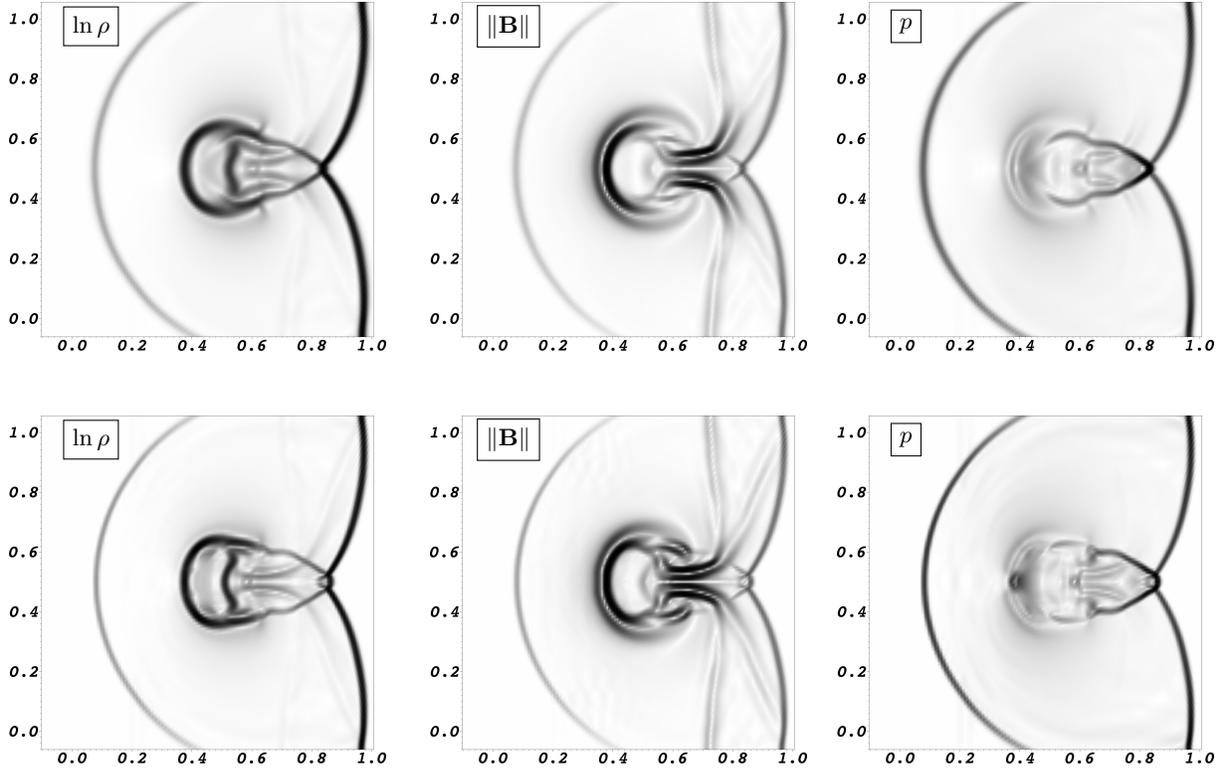
\begin{figure}[htb]
\begin{center}
\begin{tikzpicture}[scale=1]
  \useasboundingbox (0,0) rectangle (15,10);  
  \draw(-.8 ,-.4) node[anchor=south west,xshift=0pt,yshift=+0pt] {\trimfig{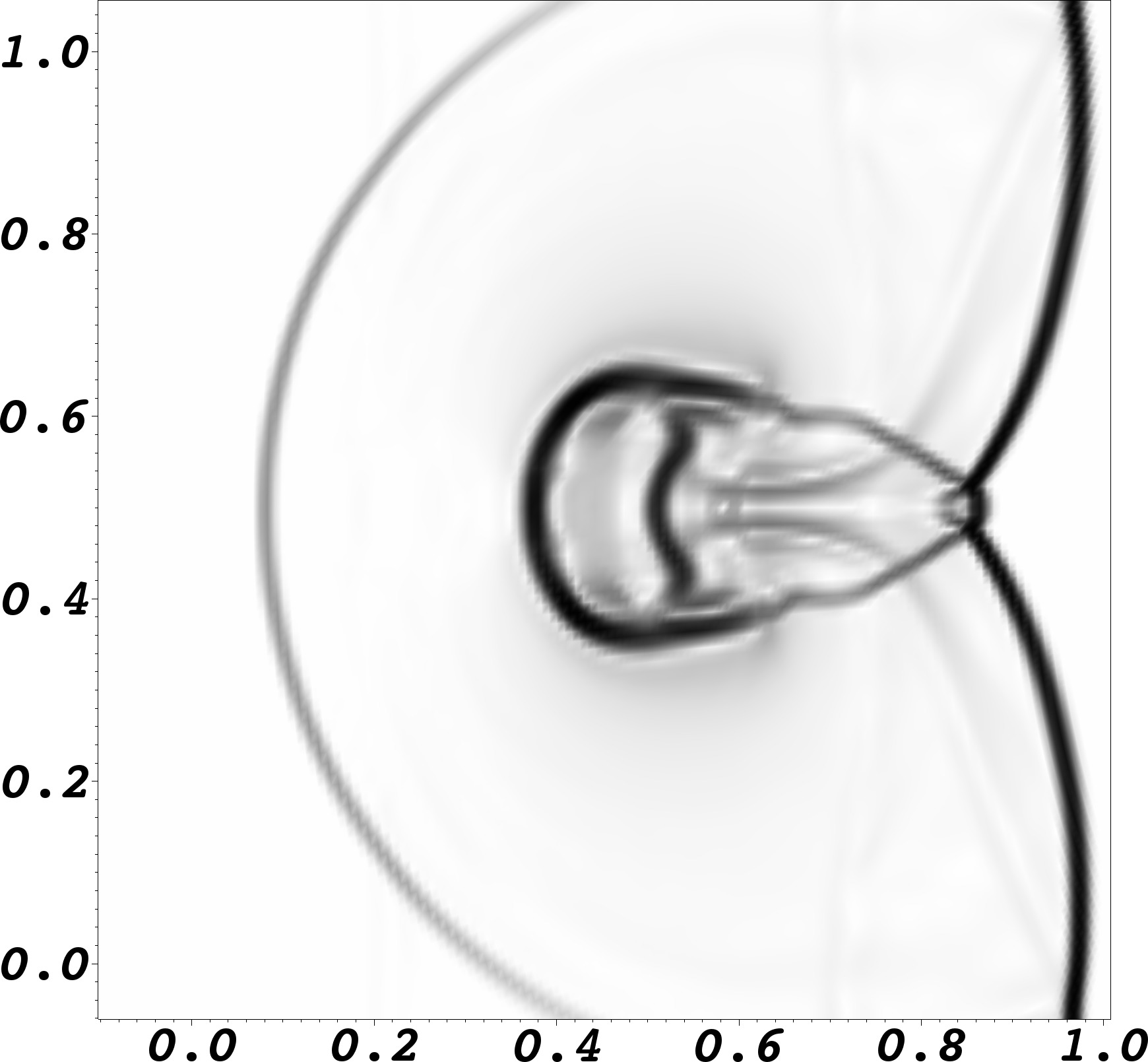}{\figWidth}};
  \draw(4.8 ,-.4) node[anchor=south west,xshift=0pt,yshift=+0pt] {\trimfig{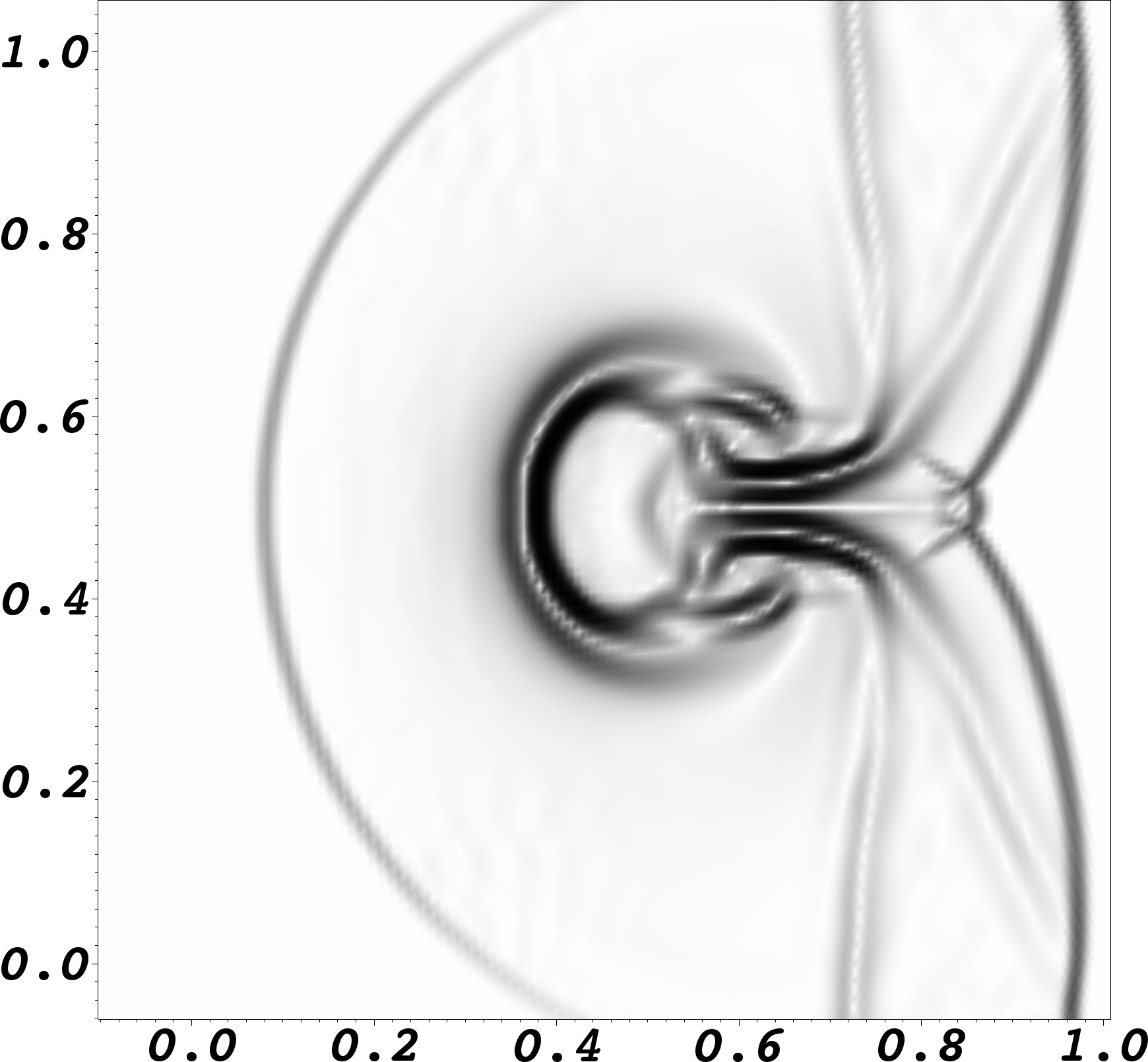}{\figWidth}};
  \draw(10.2,-.4) node[anchor=south west,xshift=0pt,yshift=+0pt] {\trimfig{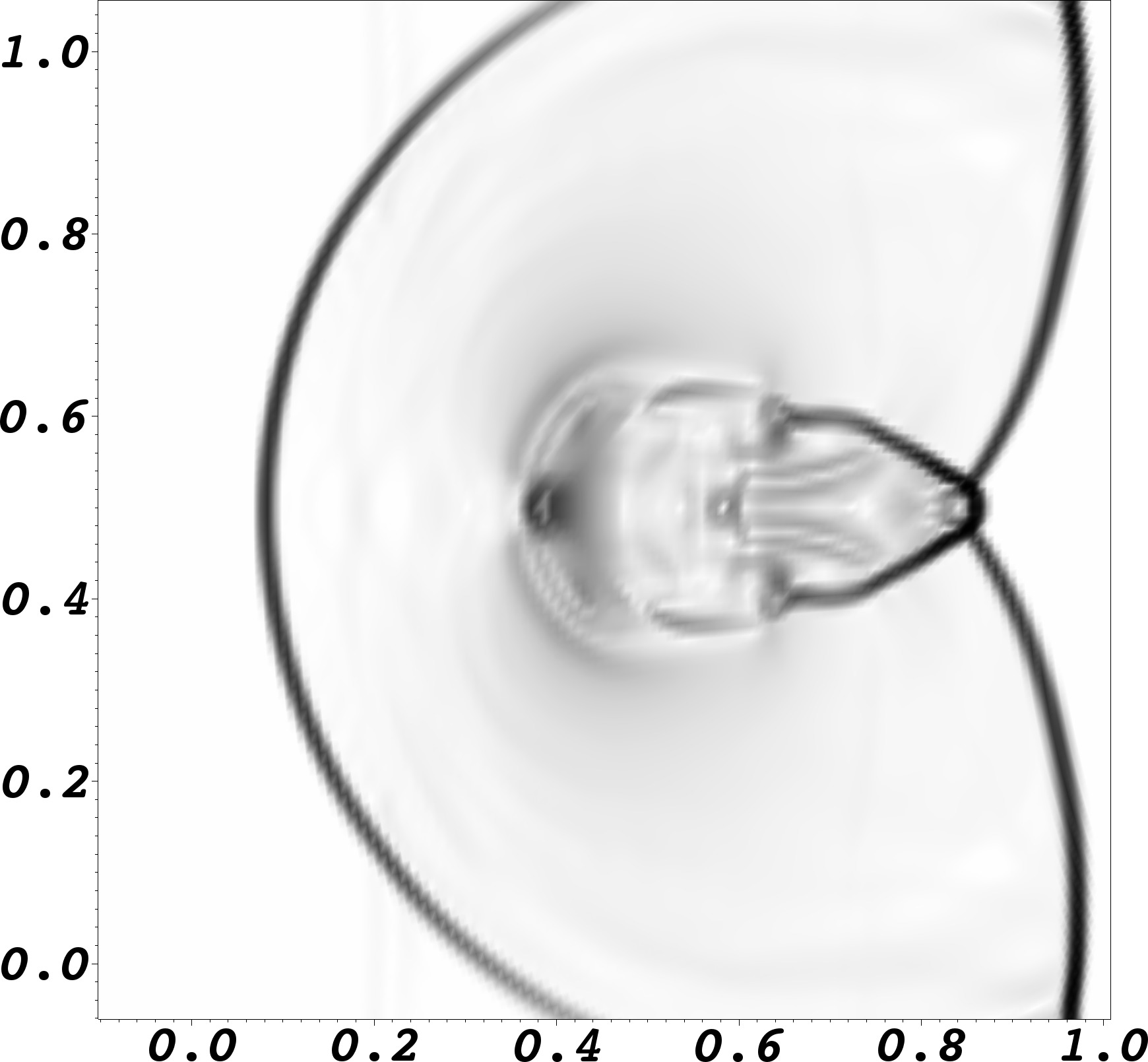}{\figWidth}};
  \draw(0,4.1) node[draw,fill=white,anchor=west,xshift=2pt,yshift=-1pt] {\small $\ln \rho$};
  \draw(5.5,4.1) node[draw,fill=white,anchor=west,xshift=2pt,yshift=-1pt] {\small $\|\Bv\|$};
  \draw(11,4.1) node[draw,fill=white,anchor=west,xshift=2pt,yshift=-1pt] {\small $p$};
  \begin{scope}[yshift=5.5cm]
  \draw(-.8 ,-.4) node[anchor=south west,xshift=0pt,yshift=+0pt] {\trimfig{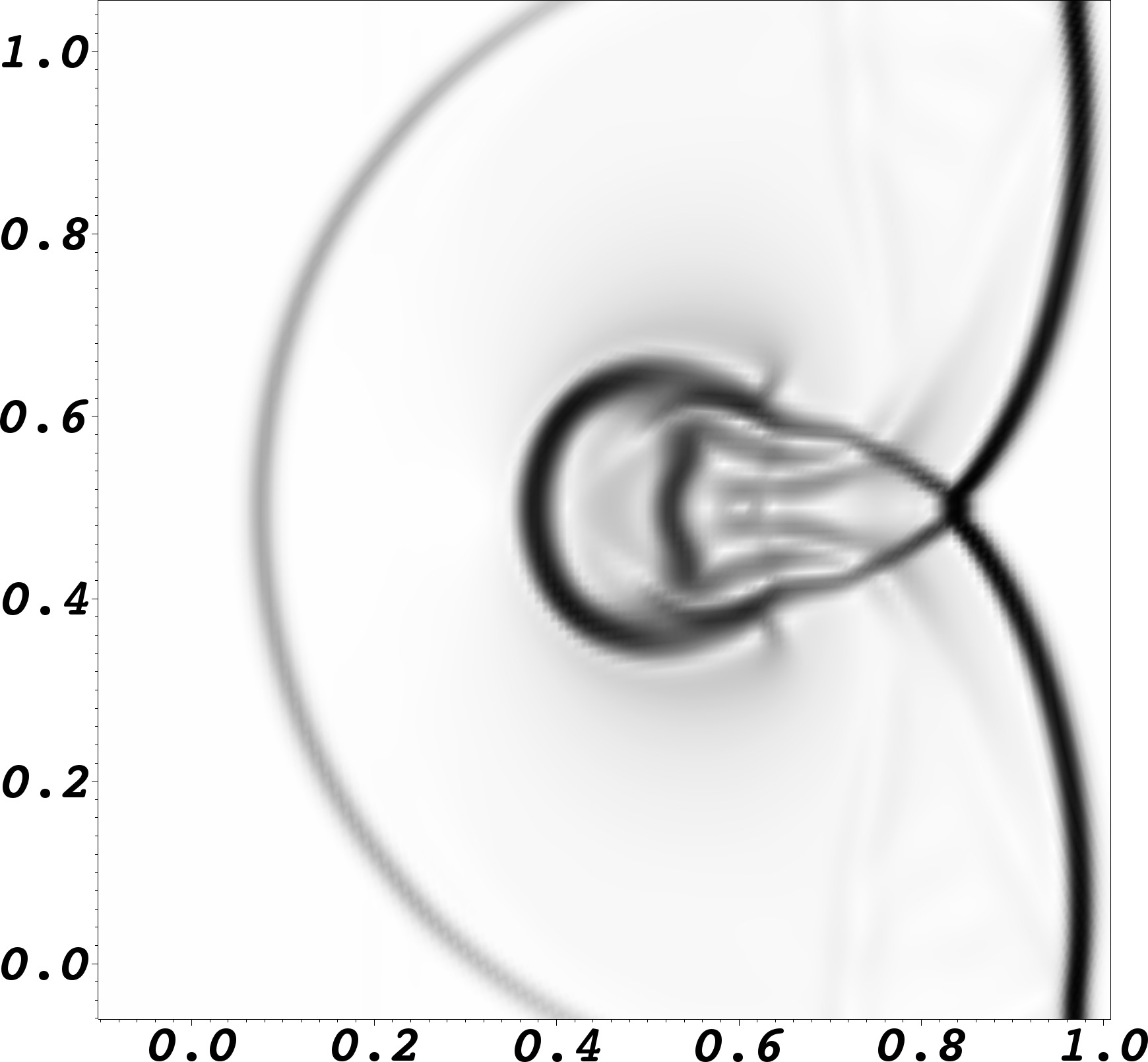}{\figWidth}};
  \draw(4.8 ,-.4) node[anchor=south west,xshift=0pt,yshift=+0pt] {\trimfig{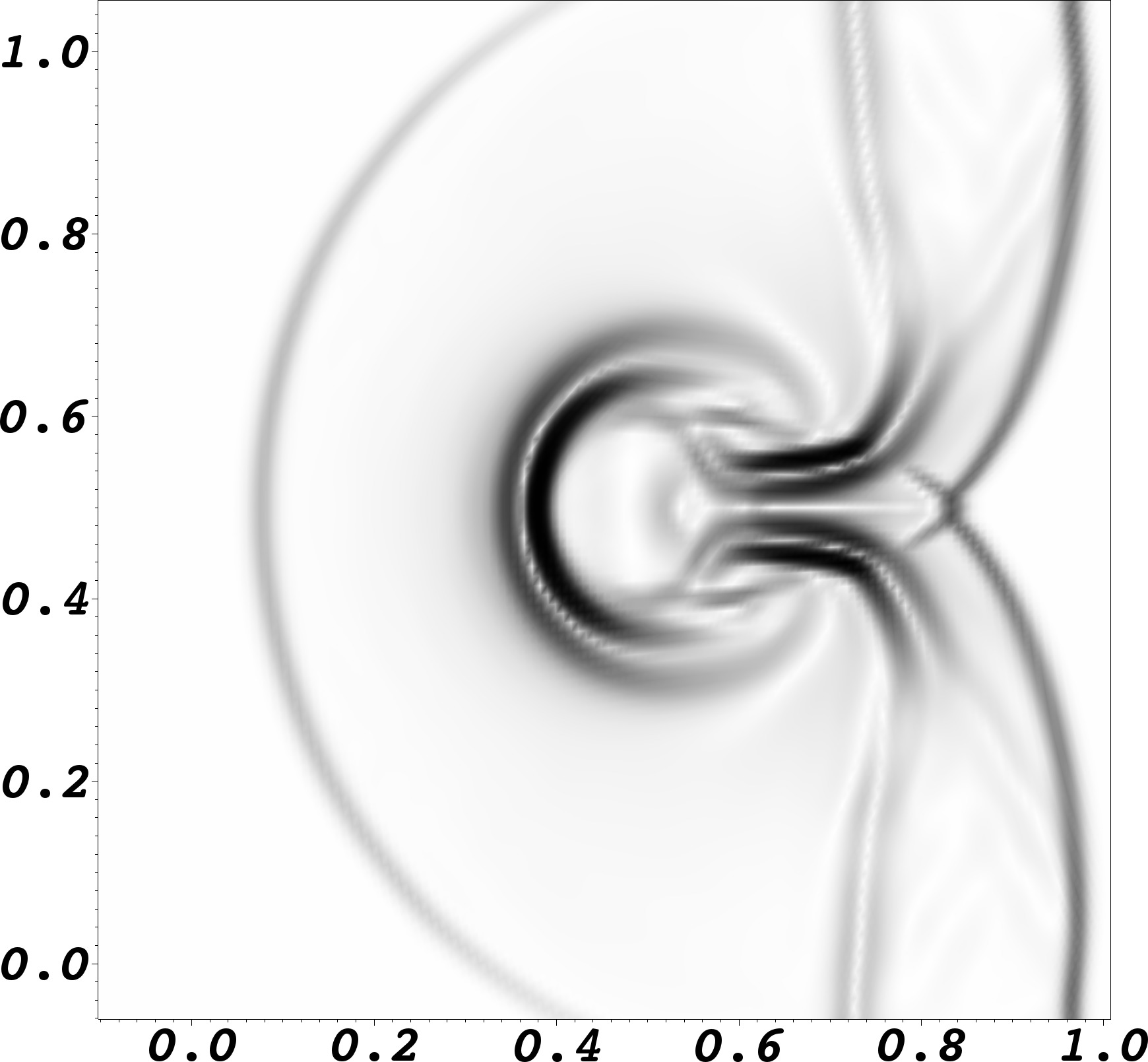}{\figWidth}};
  \draw(10.2,-.4) node[anchor=south west,xshift=0pt,yshift=+0pt] {\trimfig{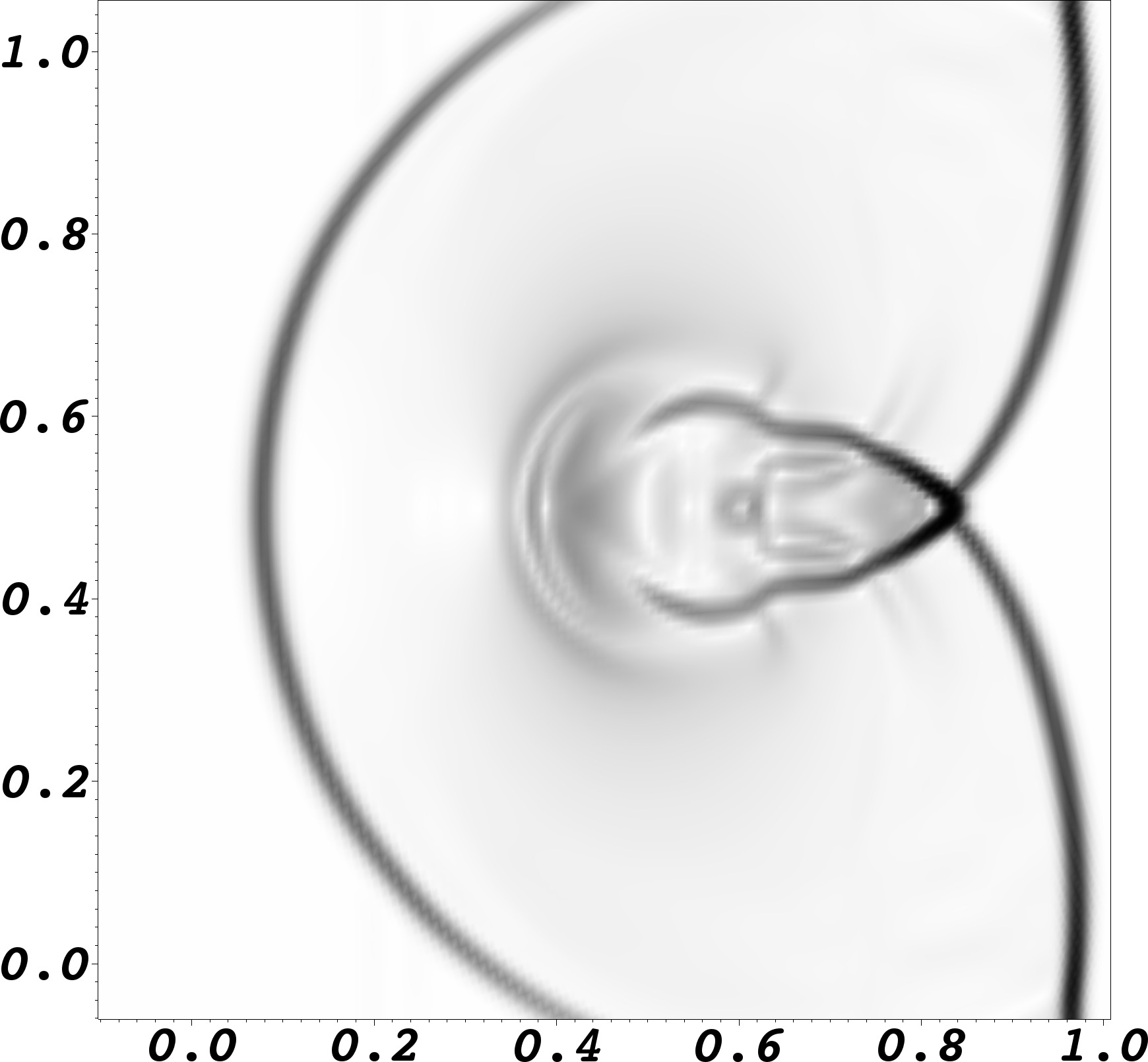}{\figWidth}};
  \draw(0,4.1) node[draw,fill=white,anchor=west,xshift=2pt,yshift=-1pt] {\small $\ln \rho$};
  \draw(5.5,4.1) node[draw,fill=white,anchor=west,xshift=2pt,yshift=-1pt] {\small $\|\Bv\|$};
  \draw(11,4.1) node[draw,fill=white,anchor=west,xshift=2pt,yshift=-1pt] {\small $p$};
  \end{scope}
%
\end{tikzpicture}
\end{center}
    \caption{
    2D cloud-shock interaction. 
Schlieren plots of the logarithm of the density, norm of the magnetic field and pressure at $t=0.06$.
Top row: WENO with the Lax-Friedrichs flux.
Bottom row: WENO with the HLLD flux.
A curvilinear grid of size $256 \times 256$ is used.       
The constrained transport step and positivity-preserving limiter are turned on.
  }
\label{fig:2DCloudShockAltWENONonUniform}
\end{figure}
}

Figure~\ref{fig:2DCloudShockAltWENONonUniform} presents the results on the curvilinear grids.
The presented results are Schlieren plots of the logarithm of the density, norm of the magnetic field and pressure.
The presented results only focus on the region near the location of the initial cloud, which contains most interesting structures.
The schemes with two numerical fluxes are used on a uniform grid of size $256 \times 256$ in the domain $(\xi, \eta)$.
The results on the curvilinear grid are comparable to the results on the uniform grid,
which verifies the solvers on general curvilinear grids.
The HLLD flux also produces less smeared solutions.
Note that the solutions on the curvilinear grid are more smeared because the effective grid  
spacing of the curvilinear grid is much larger than that of the Cartesian grid.
Same as the Cartesian grid results, the HLLD flux also performs better than the Lax-Friedrichs flux in both the bubble region and shock region.

%% file: texSISC/rotor.tex
{ 
\subsection{2D rotor problem}
\label{sec:2DRotor}
	
The initial condition is given as 
\begin{align*}
(\rho,u,v) = \left\{ \begin{array}{ll}
(10, -(y-0.5)/r_0, (x-0.5)/r_0 ), & \text{if} \ r\leq r_0,\\
(1+9 {f}(r), -f(r)(y-0.5)/r , f(r)(x-0.5)/r ), & \text{if} \ r_0<r\leq r_1, \\
(1, 0, 0), & \text{if} \ r > r_1, \\
\end{array}\right.
\end{align*}
and
\begin{align*}
w = 0, \quad B_{1} =2.5/\sqrt{4\pi}, \quad B_{2}=0, \quad B_{3}=0, \quad p=0.5, \quad A=2.5/\sqrt{4\pi}y,
\end{align*}
where $r=\sqrt{(x-0.5)^2+(y-0.5)^2}$, $r_0=0.1$, $r_1=0.115$ and $f(r) = (r_1-r)/(r_1-r_0)$.
Here we use the same initial condition of the second rotor problem test in~\cite{Toth2000}.
The problem is solved on a stationary curve grid determined by the mapping
\begin{subequations}
\label{eq:rotor_mesh}
\begin{align}
x = & \xi  - 0.5 + \epsilon_{x} \cos( \pi(\eta-0.5) )  \sin( \pi(\xi -0.5) ), \\
y = & \eta - 0.5 + \epsilon_{y} \cos( \pi(\xi -0.5) )  \sin( \pi(\eta-0.5) ), 
\end{align}
\end{subequations} 
with $\epsilon_{x}=\epsilon_{y}=0.1$ and $(\xi,\eta) \in [0,1]\times [0,1]$. 
The mach number at $t = 0.295$ on a grid of $256\times256$ are presented in Figure~\ref{fig:rotor}. 
From the zoomed view, we note that there is a significant improvement of the solutions computed by the HLLD flux.


{
\newcommand{\figWidth}{5.cm}
\newcommand{\figWidthb}{4.cm}
\newcommand{\trimfig}[2]{\trimFig{#1}{#2}{.86}{.7}{.45}{.2}}
\newcommand{\trimfigb}[2]{\trimFig{#1}{#2}{.95}{.7}{.48}{.2}}
\begin{figure}[htb]
\begin{center}
\begin{tikzpicture}[scale=1]
\useasboundingbox (0,0) rectangle (16,4.5);  
  \draw(-0.7,-0.5) node[anchor=south west,xshift=0pt,yshift=+0pt]   {\trimfig{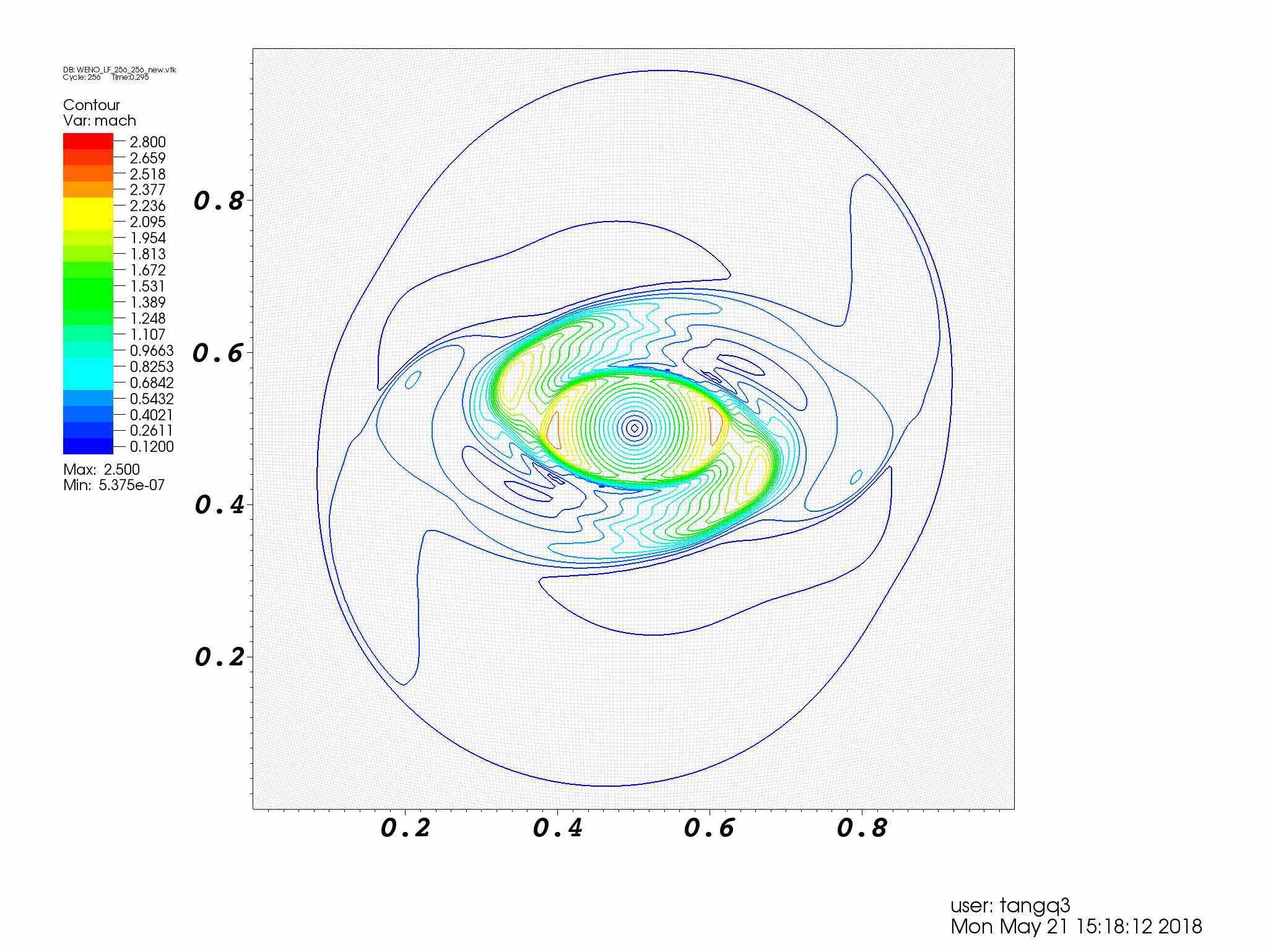}{\figWidth}};
  \draw( 3.9,-0.) node[anchor=south west,xshift=0pt,yshift=+0pt]   {\trimfigb{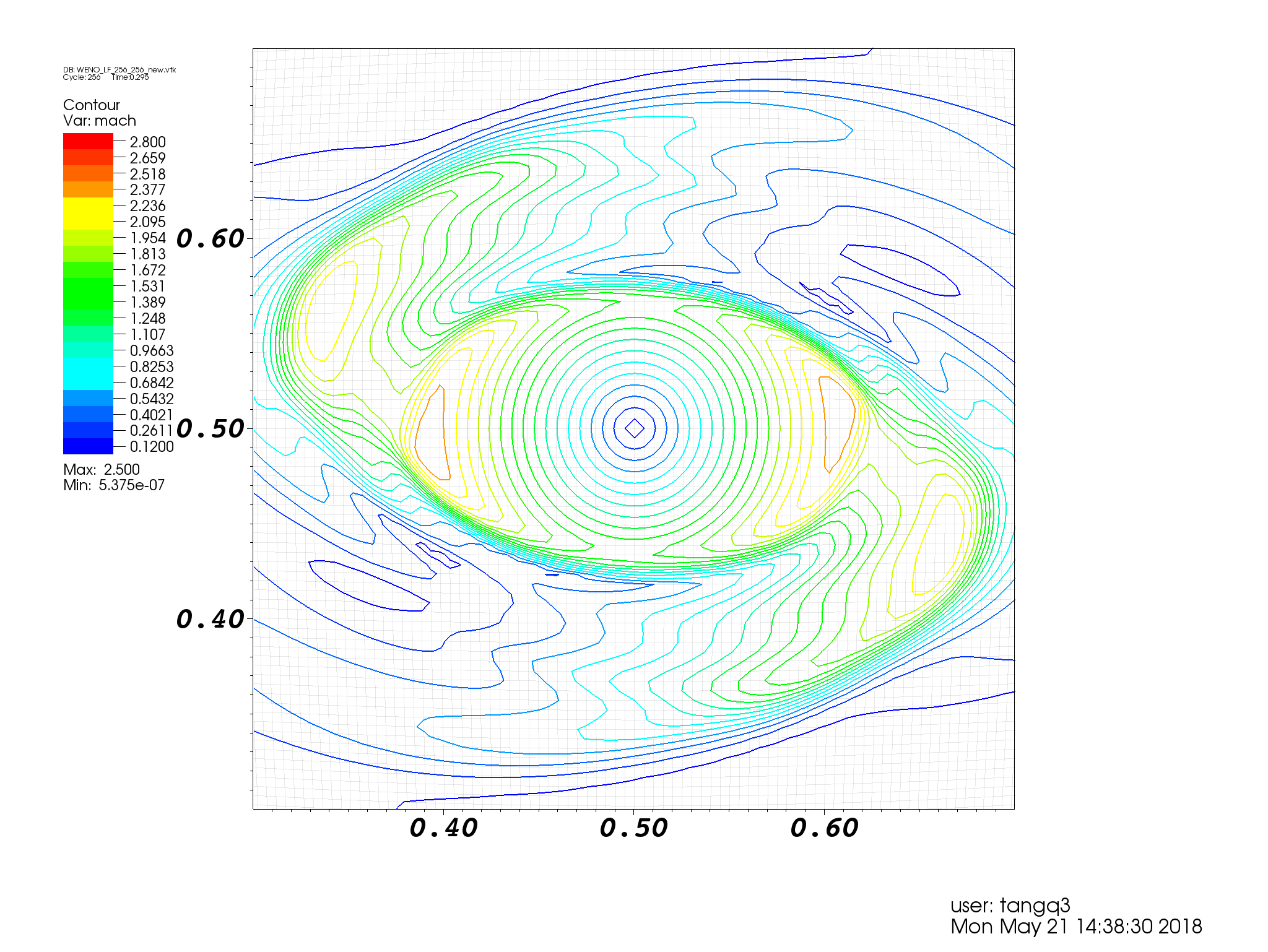}{\figWidthb}};
  \draw(7.9,-0.5) node[anchor=south west,xshift=0pt,yshift=+0pt] {\trimfig{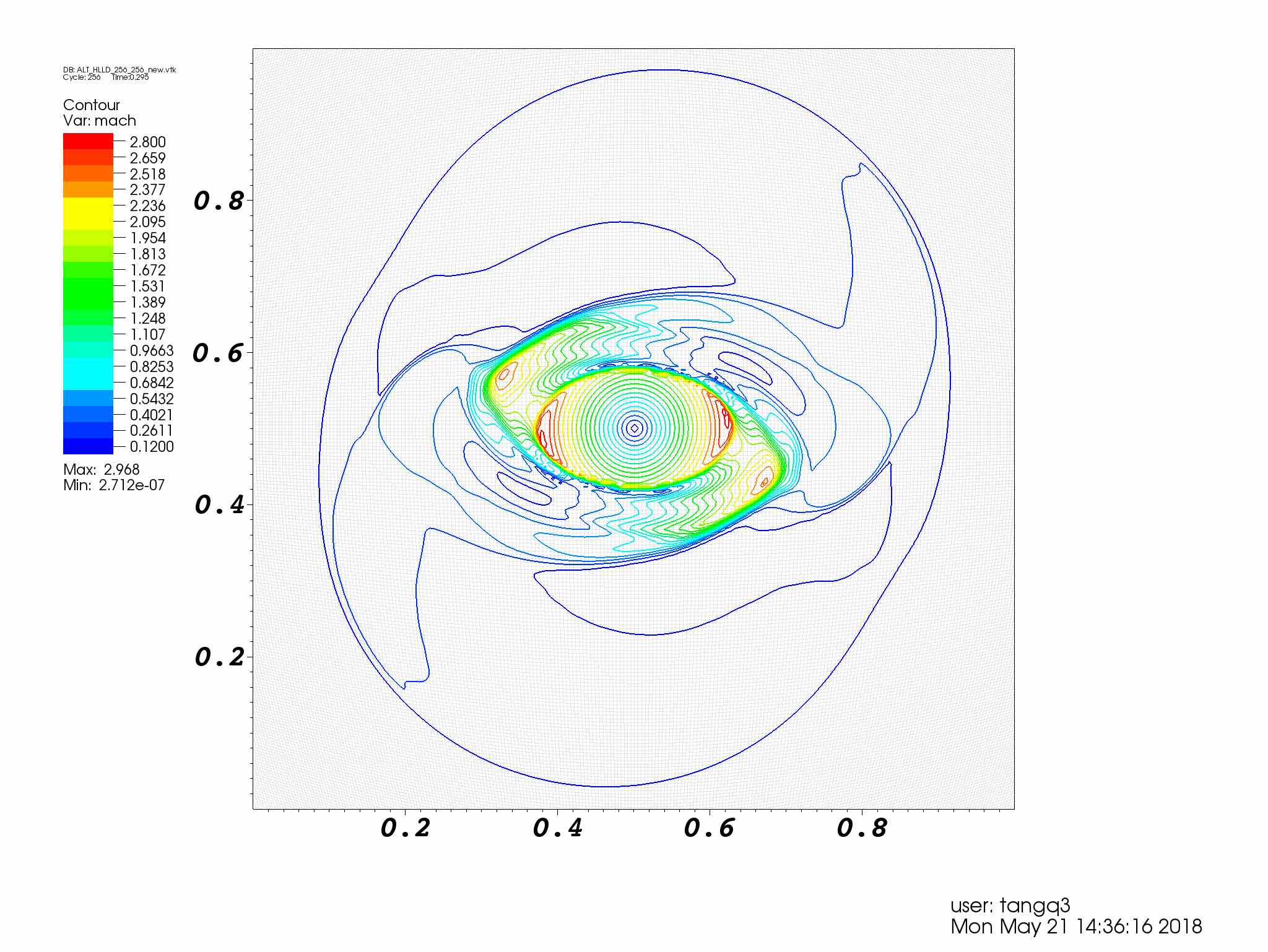}{\figWidth}};
  \draw(12.5,-0.) node[anchor=south west,xshift=0pt,yshift=+0pt] {\trimfigb{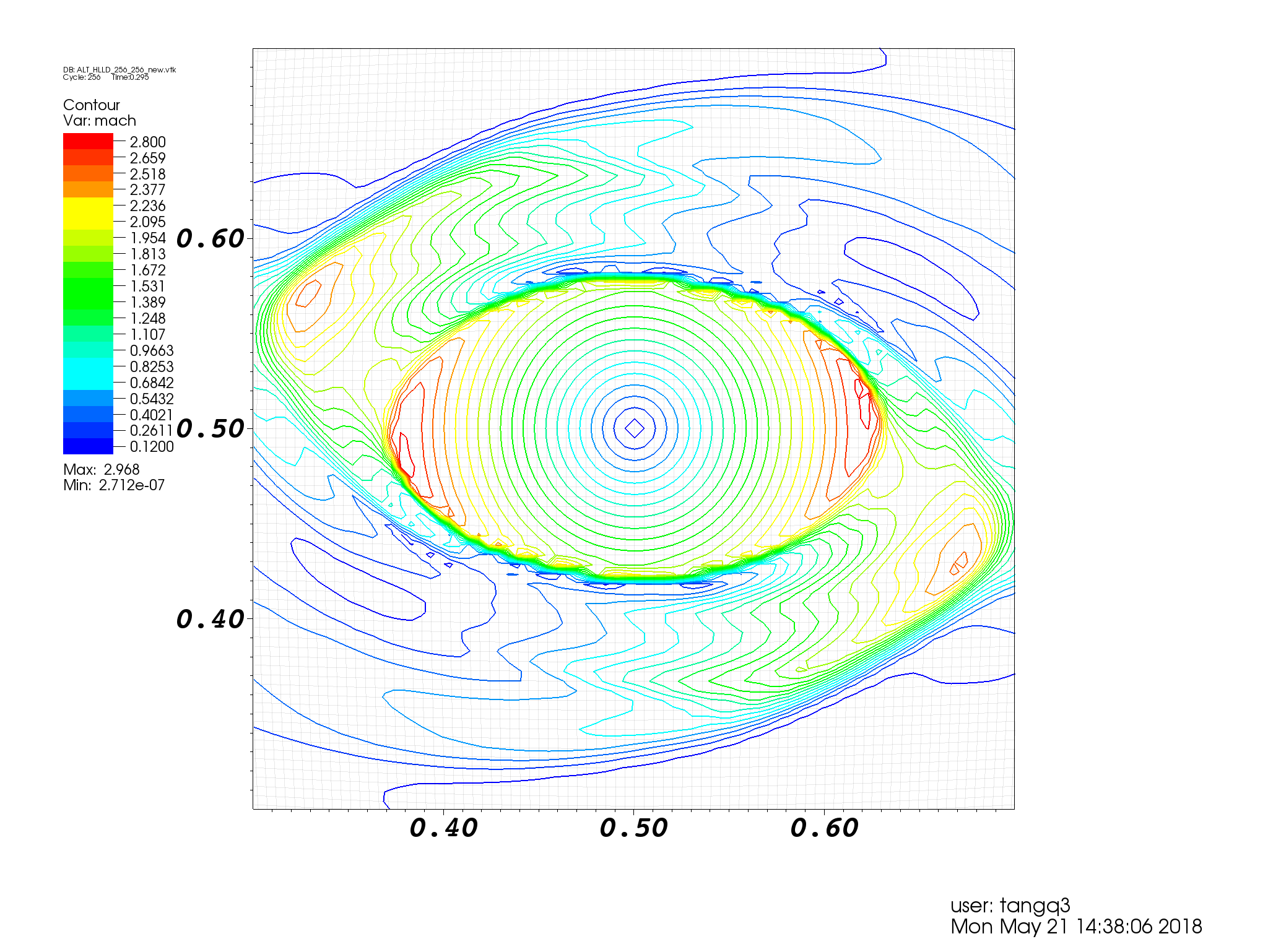}{\figWidthb}};
  \draw(-.3,4.) node[draw,fill=white,anchor=west,xshift=2pt,yshift=-1pt] {\small $\|\uv\|/c$};
  \draw(8.3,4.) node[draw,fill=white,anchor=west,xshift=2pt,yshift=-1pt] {\small $\|\uv\|/c$};
%
\end{tikzpicture}
\end{center}
    \caption{2D rotor problem.
Contour plots of the mach number and their zoomed views are presented with 20 equally spaced contour lines in the range of $[0.12, 2.8]$.  
Left two: WENO with the Lax-Friedrichs flux. 
Right two: WENO with the HLLD flux. 
The computational grid of size $256 \times 256$ is plotted (the light solid lines).
}
\label{fig:rotor}
\end{figure}
}

}

%
%
%
%
%
%
%

%% file: texSISC/blastwave.tex
\subsection{2D blast wave}
\label{sec:2DBlastWave}
	
{ 
	
Next we consider the blast wave problem in 2D.
In this test strong shocks interact with a low-$\beta$ background, which could potentially cause negative density or pressure in numerical simulations. 
The problem has been commonly used to test the positivity-preserving capabilities of numerical methods for MHD equation, see~\cite{Balsara2009, Balsara1999a, Christlieb2015, Li2011} for instance. The initial conditions contain a constant density, velocity and magnetic field
\begin{align*}
(\rho, u,v,w,B_{1}, B_{2}, B_{3} ) (0,x,y) = (1, 0, 0, 0, 50/\sqrt{2\pi}, 50/\sqrt{2\pi},0)
\end{align*}
and a piecewise defined pressure
\begin{align*}
p(0,x,y) = \left\{ \begin{array}{ll}
1000,  & \text{if} \ r\leq 0.1,\\
0.1, &  \text{otherwise},\\
\end{array}
\right.
\end{align*}
where $r$ is the distance to the origin.
The initial magnetic potential is 
$$A(0,x,y) = 50/\sqrt{2\pi} (y-x).$$
Here we use the grid given in~\eqref{eq:rotor_mesh} and the boundary condition identical to the previous test in Section~\ref{sec:2DRotor}.
Figure~\ref{fig:compareBlastWave}
shows the numerical solutions at $t = 0.01$ solved on a grid of $256\times256$, which match well with the previous work. 
Note that the solutions of the HLLD flux are slightly less diffusive than the solutions of the Lax-Friedrichs flux.


{
\newcommand{\figWidth}{4.5cm}
\newcommand{\figWidthb}{4.35cm}
\newcommand{\figWidthc}{4.5cm}
\newcommand{\trimfig}[2]{\trimFig{#1}{#2}{.87}{.7}{.45}{.2}}
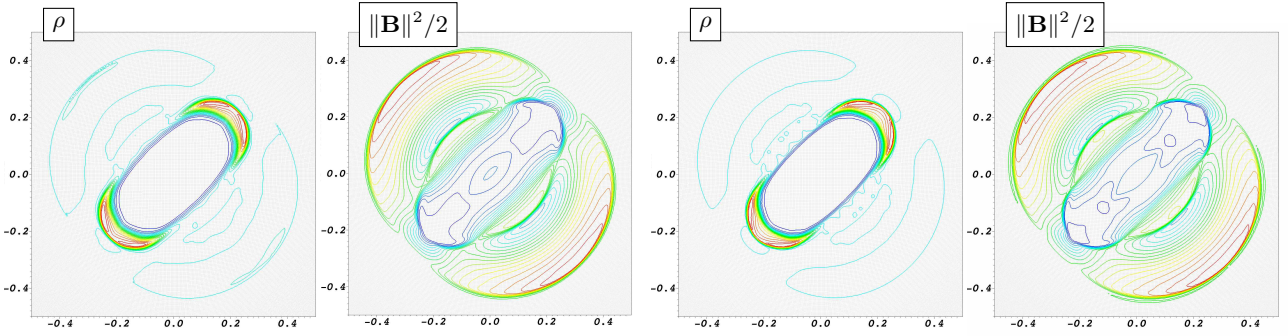
\begin{figure}[htb]
\begin{center}
\begin{tikzpicture}[scale=1]
\useasboundingbox (0,0) rectangle (16,4.2);  
  \draw(-0.7,-0.5) node[anchor=south west,xshift=0pt,yshift=+0pt] {\trimfig{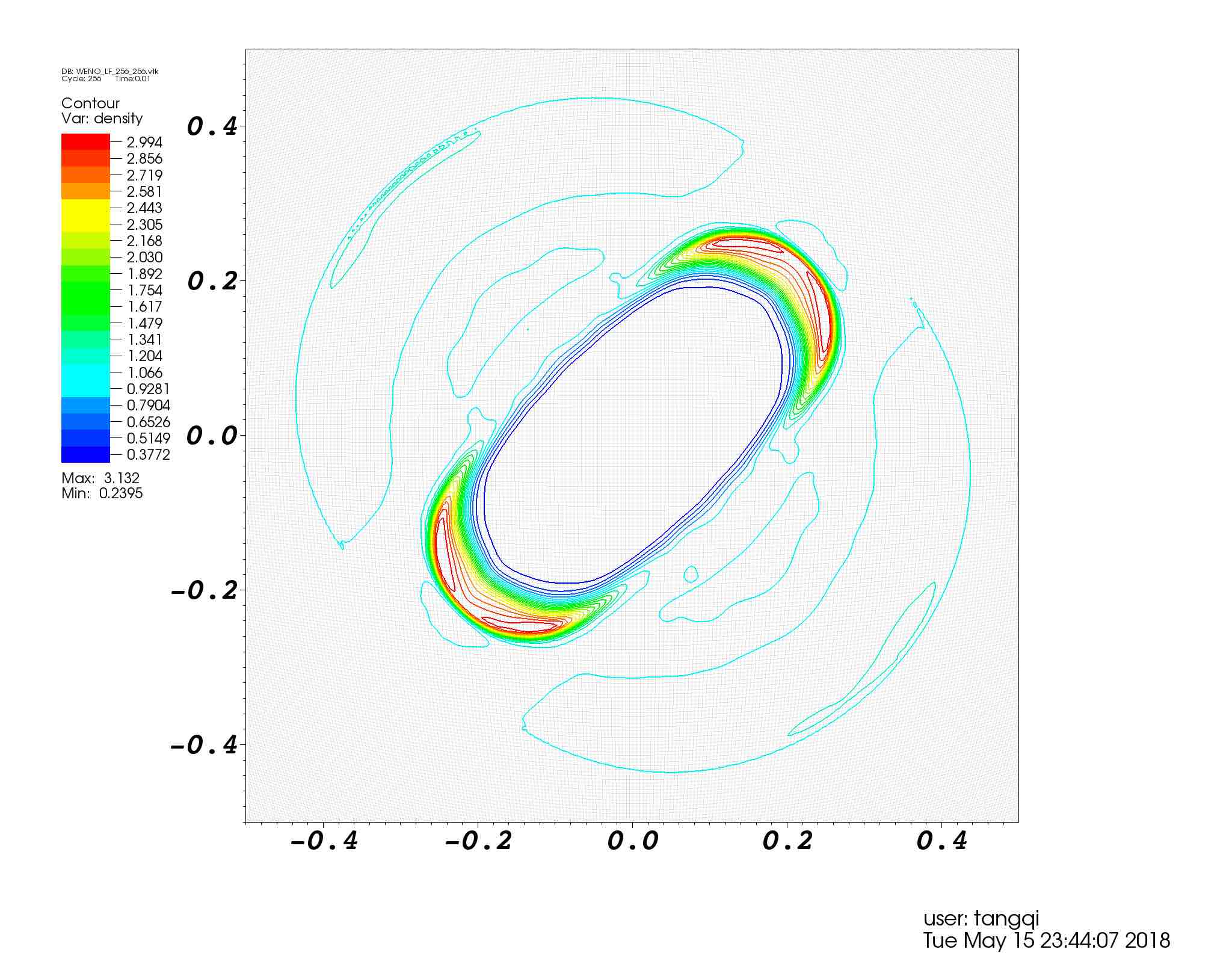}{\figWidth}};
  \draw( 3.5,-0.5) node[anchor=south west,xshift=0pt,yshift=+0pt] {\trimfig{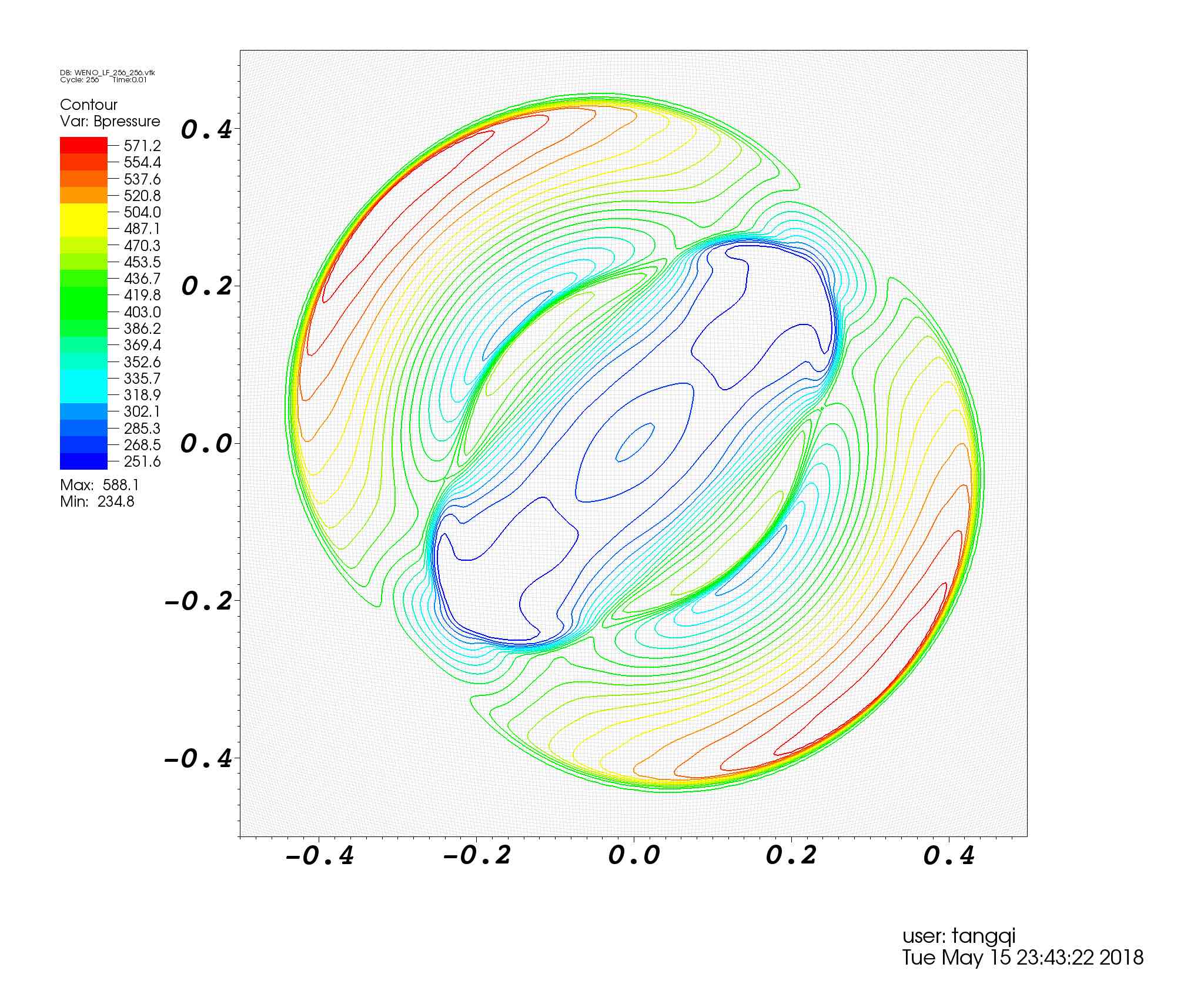}{\figWidthb}};
  \draw(   7.9,-0.5) node[anchor=south west,xshift=0pt,yshift=+0pt] {\trimfig{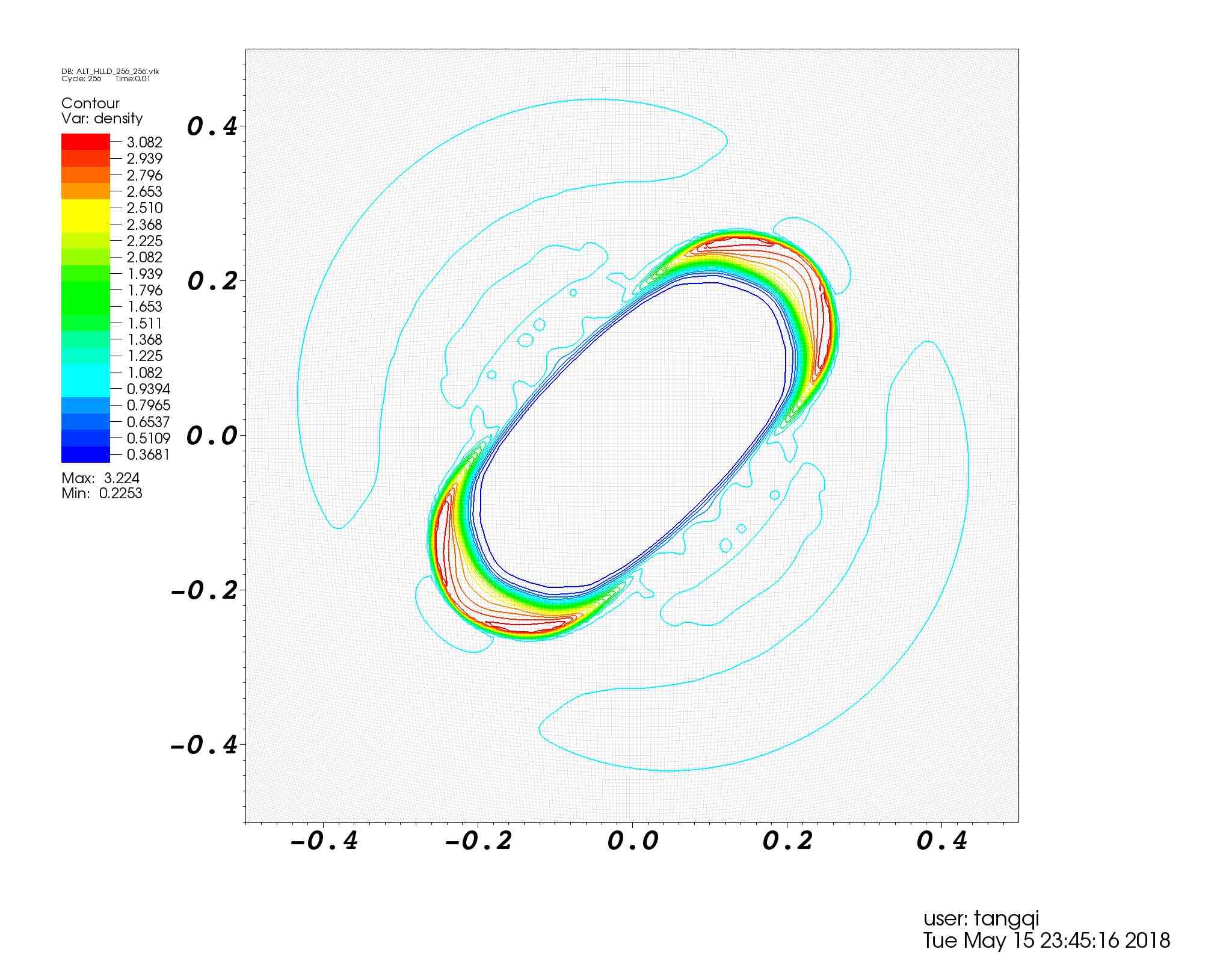}{\figWidth}};
  \draw(  12.1,-0.5) node[anchor=south west,xshift=0pt,yshift=+0pt] {\trimfig{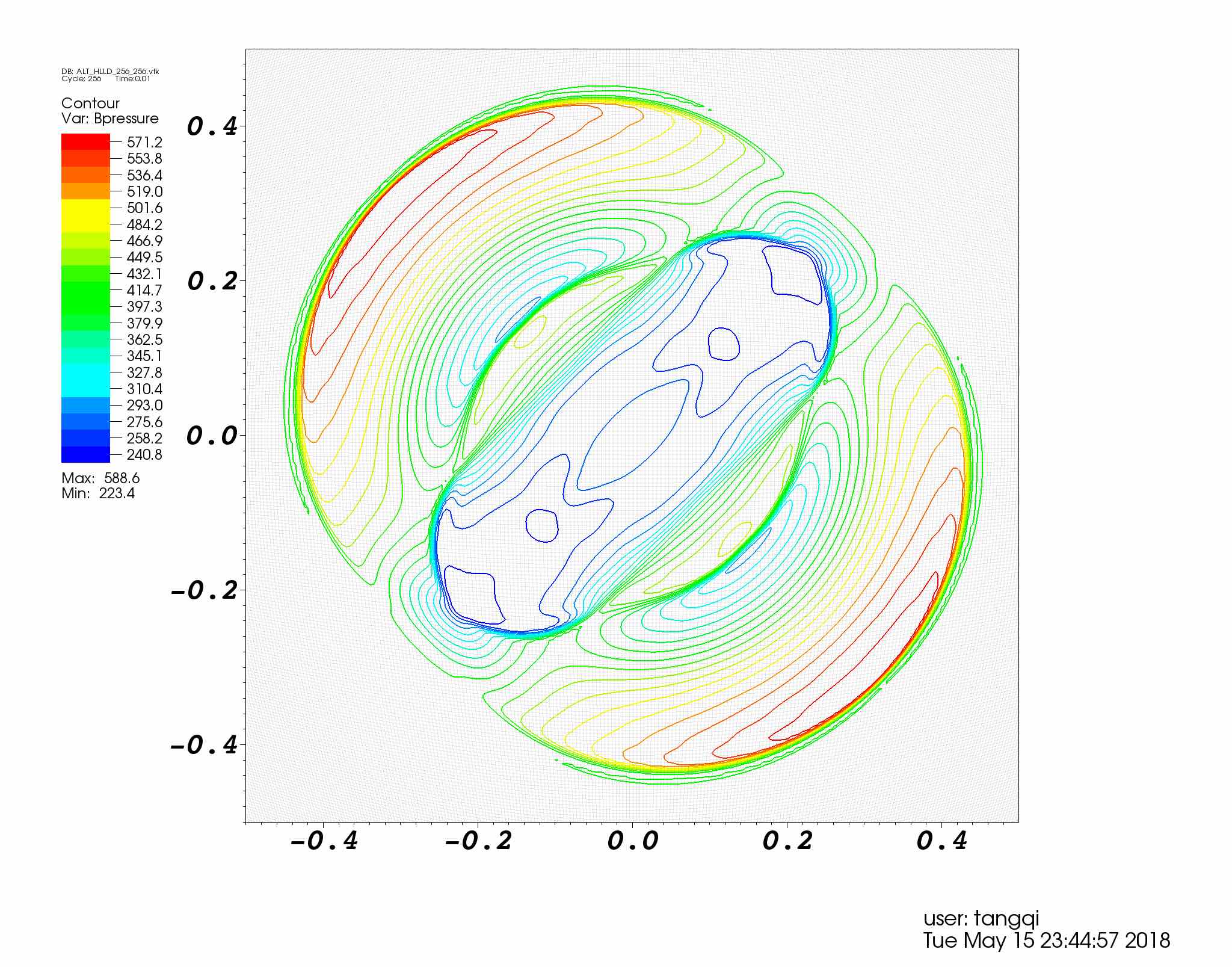}{\figWidthc}};
  \draw(-.1,4.) node[draw,fill=white,anchor=west,xshift=2pt,yshift=-1pt] {\small $\rho$};
  \draw(4.1,4.) node[draw,fill=white,anchor=west,xshift=2pt,yshift=-1pt] {\small $\|\Bv\|^2/2$};
  \draw(8.5,4.) node[draw,fill=white,anchor=west,xshift=2pt,yshift=-1pt] {\small $\rho$};
  \draw(12.7,4.) node[draw,fill=white,anchor=west,xshift=2pt,yshift=-1pt] {\small $\|\Bv\|^2/2$};
%
\end{tikzpicture}
\end{center}
    \caption{2D blast wave.
Contour plots of the density and magnetic potential at $t=0.01$ are presented with 20 equally spaced contour lines.  
Left two: WENO with the Lax-Friedrichs flux. 
Right two: WENO with the HLLD flux. 
The computational grid of size $256 \times 256$ is plotted (the light solid lines).
}
\label{fig:compareBlastWave}
\end{figure}
}
}

%% file: texSISC/bowShock.tex
\subsection{Bow shock flow}
\label{sec:bowShock}

We end our numerical investigation with a stationary bow shock flow. A bow shock flow benchmark has been 
previously considered in~\cite{de1999numerical, de2001,nishida2009} but low-order finite volume methods were used therein.
Here we use a similar problem to examine the performance of
the WENO schemes when applied to problems involving curved physical boundaries.
The compatibility boundary condition derived in Section~\ref{sec:boundaryConditions} is also verified using this example.

The computational domain is determined by the mapping
\begin{align*}
    x & = (r_1 - (r_1-r_0)\xi) \cos(\pi+(1-2\eta)\theta), \\
    y & = (r_2 - (r_2-r_0)\xi) \sin(\pi+(1-2\eta)\theta),
\end{align*}
with $r_1 = 0.3$, $r_2 = 0.65$, $r_0=0.125$, $\theta = 5\pi/12$ and $(\xi, \eta) \in [0, 1]\times [0,1]$.
A constant initial condition was used in~\cite{de1999numerical, de2001,nishida2009} but it is not compatible with the PEC boundary initially.
Typically, such a constant magnetic field will be first projected to satisfy the boundary condition in a solver
but note that its projected field is not divergence-free.
Although such a treatment is very common in the incompressible flow, 
it may cause some troubles for high-order methods in the MHD equations,
since controlling divergence errors is more critical in the MHD equations.
To avoid the potential issue, we modify the initial condition by ramping the constant magnetic field 
in a small annular region of distance $\delta r = 0.125$,
so that both the PEC boundary and divergence-free condition are satisfied exactly in the initial condition.
Assuming $r = \sqrt{x^2 + y^2}$, we impose an initial condition of $\rho = 1$, $p=0.2$, $\uv = (2,0,0)^T$ 
and the magnetic field
\begin{equation*}
    \Bv = \begin{cases}
        ( B_1 , B_2 , 0)^T \quad
      & \text{if $r \leq r_0+\delta r$}, \\
        (0.1, 0 , 0)^T \quad
      & \text{otherwise},
    \end{cases}
\end{equation*}
with
\begin{align*}
    B_1 & = 0.1 \frac{\pi y^2}{ 2\, \delta r \, r} \cos\left(\dfrac{\pi(r-r_0)}{2\, \delta r}\right) 
    +  0.1 \sin\left(\dfrac{\pi(r-r_0)}{2\, \delta r}\right), \\
    B_2 & = - 0.1 \frac{\pi xy}{2 \, \delta r \, r} \cos\left(\dfrac{\pi(r-r_0)}{2\, \delta r}\right).
\end{align*}
The corresponding initial magnetic potential is
\begin{equation*}
    A  = \begin{cases}
       0.1\,y \sin\left(\dfrac{\pi(r-r_0)}{2\, \delta r}\right) \quad
      & \text{if $r \leq r_0+\delta r$}, \\
        0.1\,y \quad
      & \text{otherwise},
    \end{cases}
\end{equation*}
The PEC boundary is applied at $\xi = 1$, i.e., $r = r_0$.
An inflow boundary condition is applied at $\xi = 0$, and the outflow boundary condition is applied at the
other two boundaries.
A uniform grid of size $120\times160$ in the domain $(\xi, \eta)$ is used for all the results presented in this section.
Note that the results of this test are axisymmetric,
and therefore we only present the results in the top half-plane.

It is found that if $\alpha$ in the Lax-Friedrichs flux~\eqref{eq:LFFlux} is estimated from the whole domain, 
the resultant value becomes too large for a reasonable CFL number in this case. 
This issue appears in the both WENO scheme and first-order scheme using the Lax-Friedrichs flux. 
Therefore, the local Lax-Friedrichs flux (the Rusanov flux) is used as one of two Riemann solvers tested for this problem.
Note that the positivity-preserving limiter used in the previous examples relies on the fact that the low-order flux is positivity-preserving. 
But for all the low-order fluxes (the local Lax-Friedrichs and HLL-type solvers) we test, 
a negative solution always appear for a CFL number larger than $0.2$.
Therefore, a slightly smaller CFL number of $0.2$ is used in this case, 
which appears to be enough to eliminate numerical solutions becoming negative.

{
\newcommand{\figWidth}{6.2cm}
\newcommand{\trimfig}[2]{\trimFig{#1}{#2}{.72}{.5}{.33}{.2}}
\begin{figure}[htb]
\begin{center}
\begin{tikzpicture}[scale=1]
  \useasboundingbox (0,0) rectangle (12,5.3);  
  \draw(-.6,-0.7) node[anchor=south west,xshift=0pt,yshift=+0pt] {\trimfig{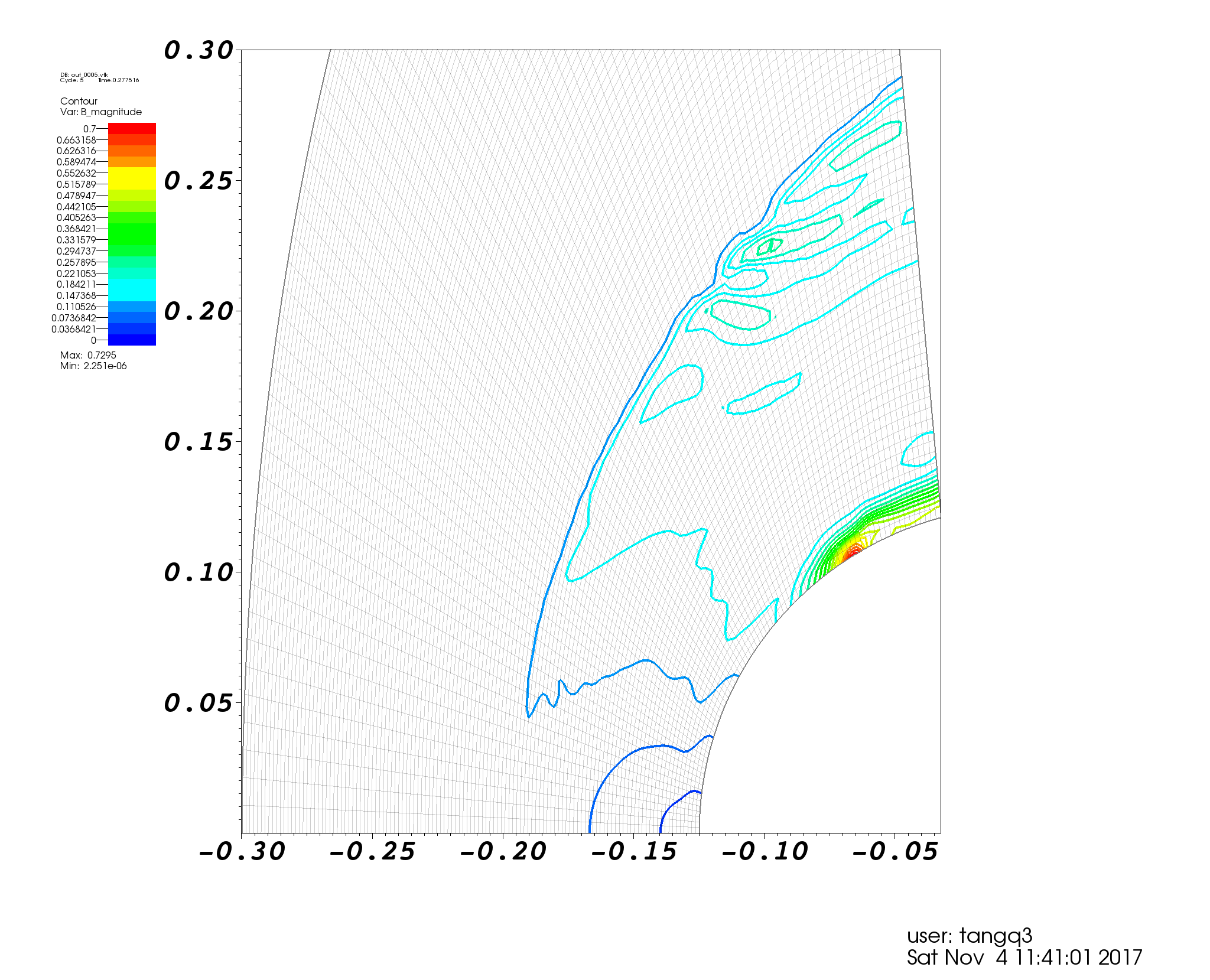}{\figWidth}};
  \draw(7.,-0.7) node[anchor=south west,xshift=0pt,yshift=+0pt] {\trimfig{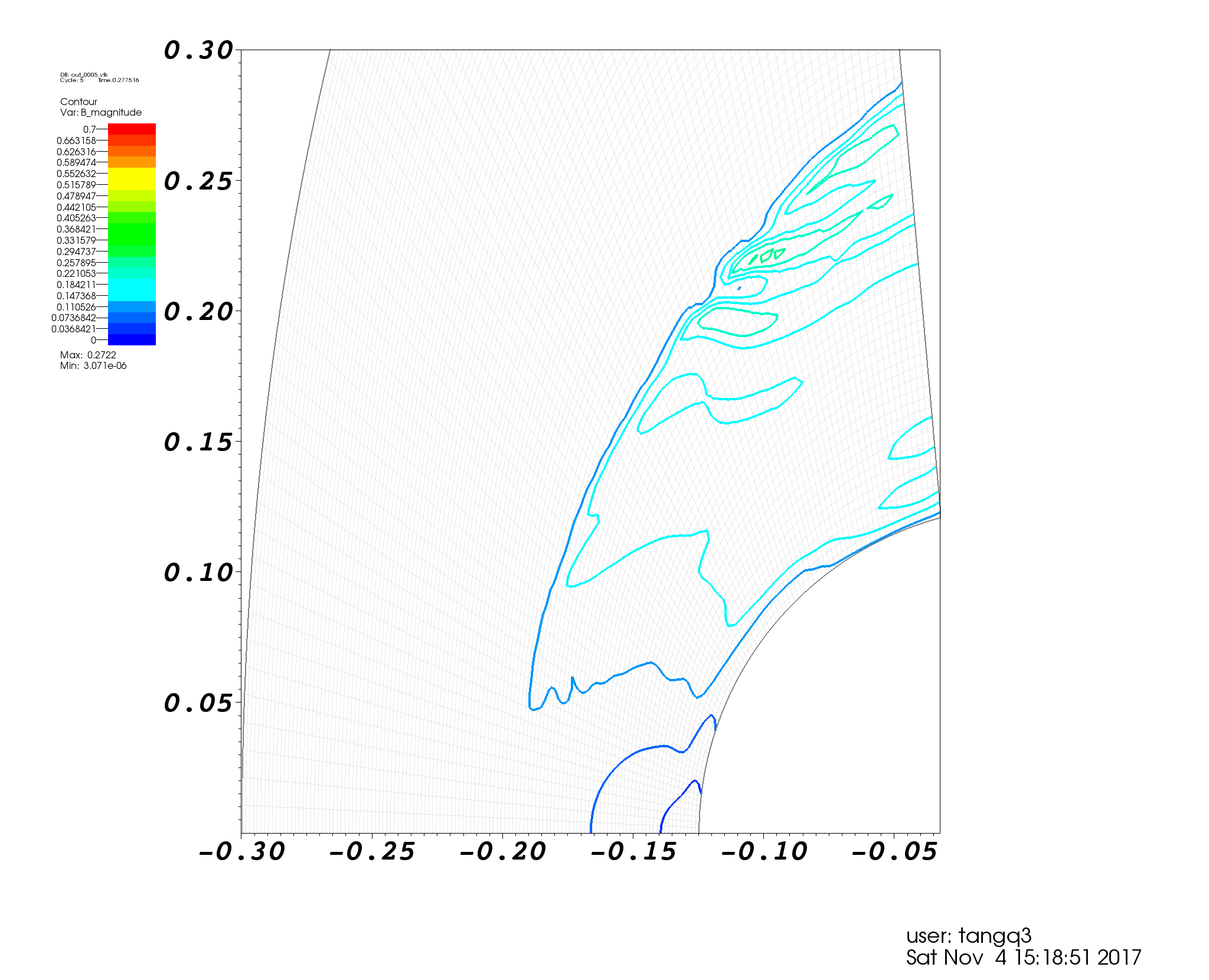}{\figWidth}};
  \draw(.2, 5.) node[draw,fill=white,anchor=west,xshift=2pt,yshift=-1pt] {\small $\|\Bv\|$};
  \draw(7.8,5.) node[draw,fill=white,anchor=west,xshift=2pt,yshift=-1pt] {\small $\|\Bv\|$};
%
\end{tikzpicture}
\end{center}
    \caption{Bow shock flow.
Contour plots of magnetic fields at $t=0.28$ are presented with 20 equally spaced contour lines from $\|\Bv\| =0$ to $\|\Bv\| = 0.7$.  
Left: the solution using the reflective boundary condition.
Right: the solution using the compatibility boundary condition. 
The constrained transport step is not turned on.}
\label{fig:compareBoundary}
\end{figure}
}

We first use this problem to study the numerical boundary condition for the PEC boundary.
In previous work~\cite{de1999numerical, de2001}, a reflective boundary condition was used and
it produced a satisfactory results in those low-order finite volume methods.
This fact has been confirmed in our own implementation, when the problem is simulated 
using the first-order numerical fluxes coupled with the forward Euler method.
However, for the high-order methods derived in the current work, the reflective boundary condition 
generates a spurious magnetic field along the PEC surface.
Figure~\ref{fig:compareBoundary} presents the magnetic fields at $t = 0.28$ generated by the reflective boundary condition and
the compatibility boundary condition we derive in Section~\ref{sec:boundaryConditions}.
Note that the magnetic field along the PEC surface is smooth for the solutions using the compatibility boundary condition
while there is a large magnetic field generated along the surface in the results using the reflective boundary condition. 
This unphysical field becomes even larger as time evolves and eventually leads to a negative pressure,
causing the failure of the solver around $t = 0.5$,
while the solver using the compatibility boundary condition remains stable to produce a stationary bow shock profile, despite
some oscillations due to divergence errors found near the shock front.
Here the WENO scheme with a local Lax-Friedrichs flux is used with the constrained transport method turned off.
The constrained transport method is not turned on for two reasons. 
First, the constrained transport method or other approaches to control the divergence errors, such as the non-conservative source terms used in the previous work~\cite{de1999numerical, de2001,nishida2009}, may diminish this issue to some extent,
by introducing some dissipations to damp the divergence errors.
In order to isolate different issues, those approaches are not used in the study of boundary conditions.
Note that the base scheme using the compatibility condition is stable for the whole simulation,
which clearly indicates the issue is not a direct consequence of the divergence error.
Second, it is not clear how to implement a boundary condition for the potential which is consistent with the reflective boundary condition of the magnetic field.


{
\newcommand{\figWidth}{5.5cm}
\newcommand{\trimfig}[2]{\trimFig{#1}{#2}{.745}{.7}{.4}{.2}}
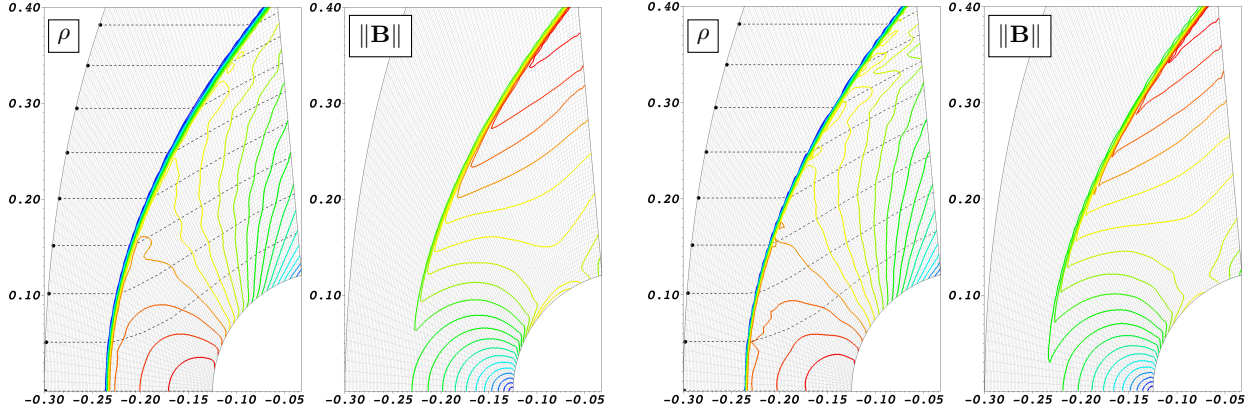
\begin{figure}[htb]
\begin{center}
\begin{tikzpicture}[scale=1]
\useasboundingbox (0,0) rectangle (16,5.);  
  \draw(-0.5,-0.8) node[anchor=south west,xshift=0pt,yshift=+0pt] {\trimfig{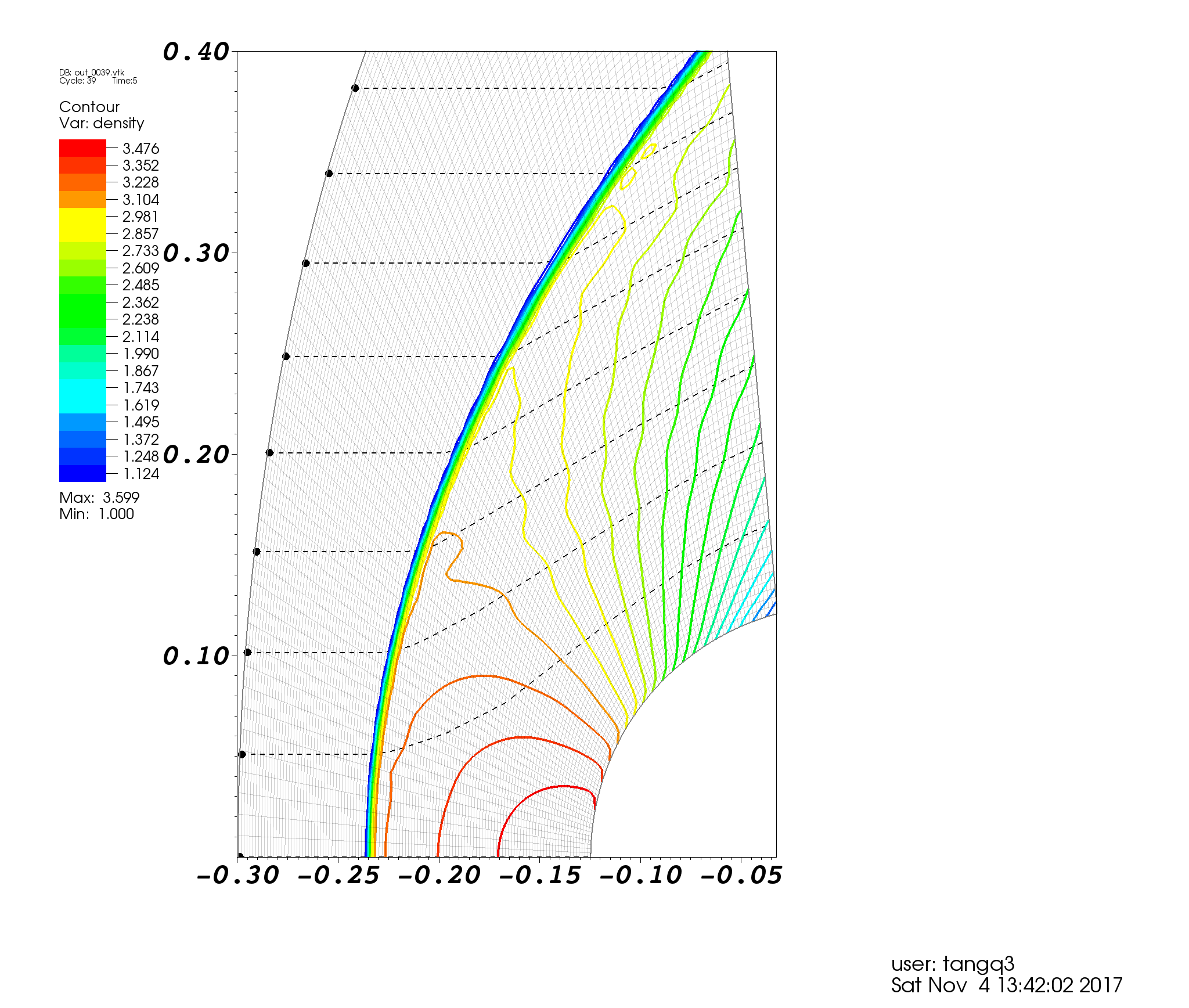}{\figWidth}};
  \draw( 3.5,-0.8) node[anchor=south west,xshift=0pt,yshift=+0pt] {\trimfig{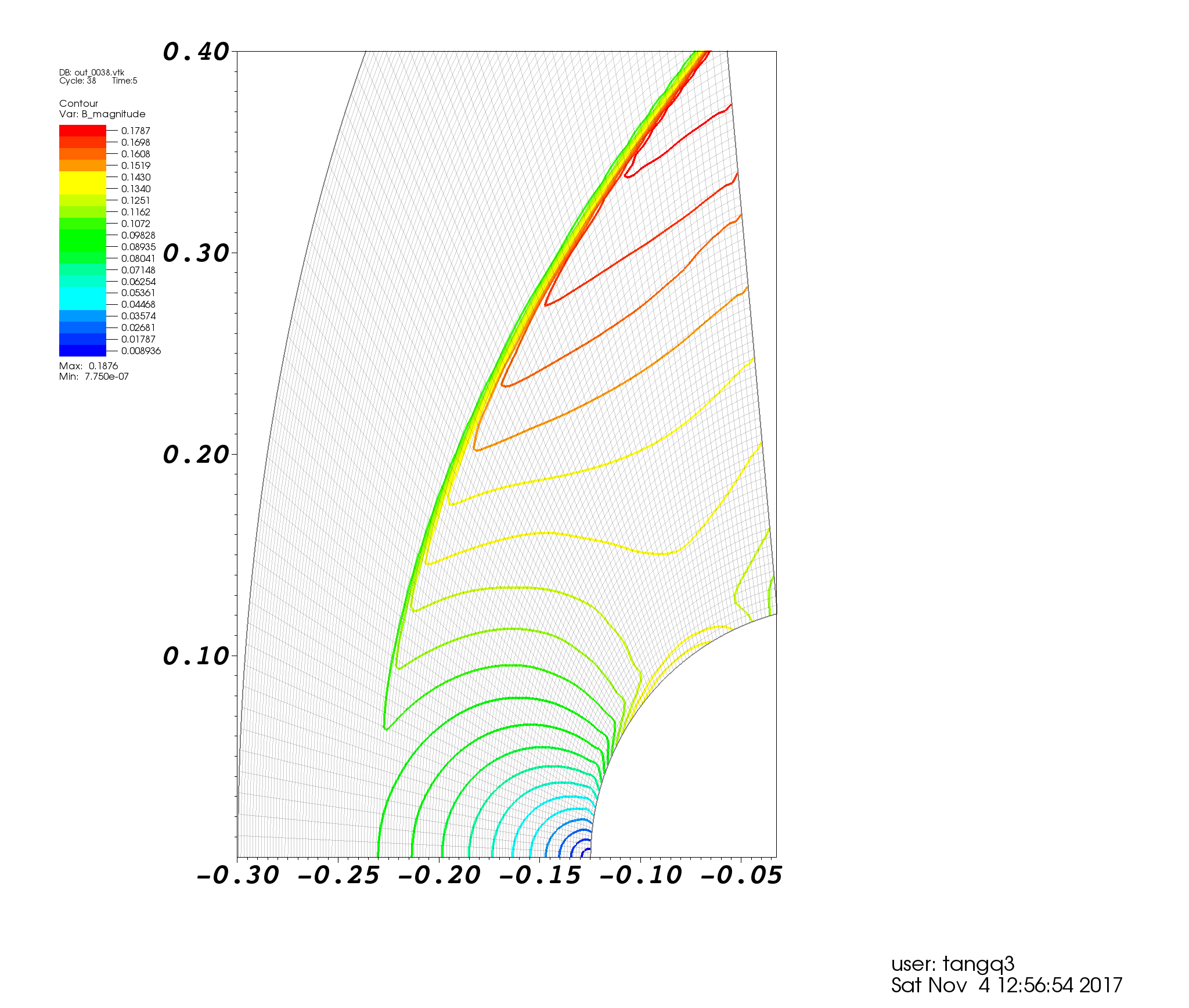}{\figWidth}};
  \draw(   8,-0.8) node[anchor=south west,xshift=0pt,yshift=+0pt] {\trimfig{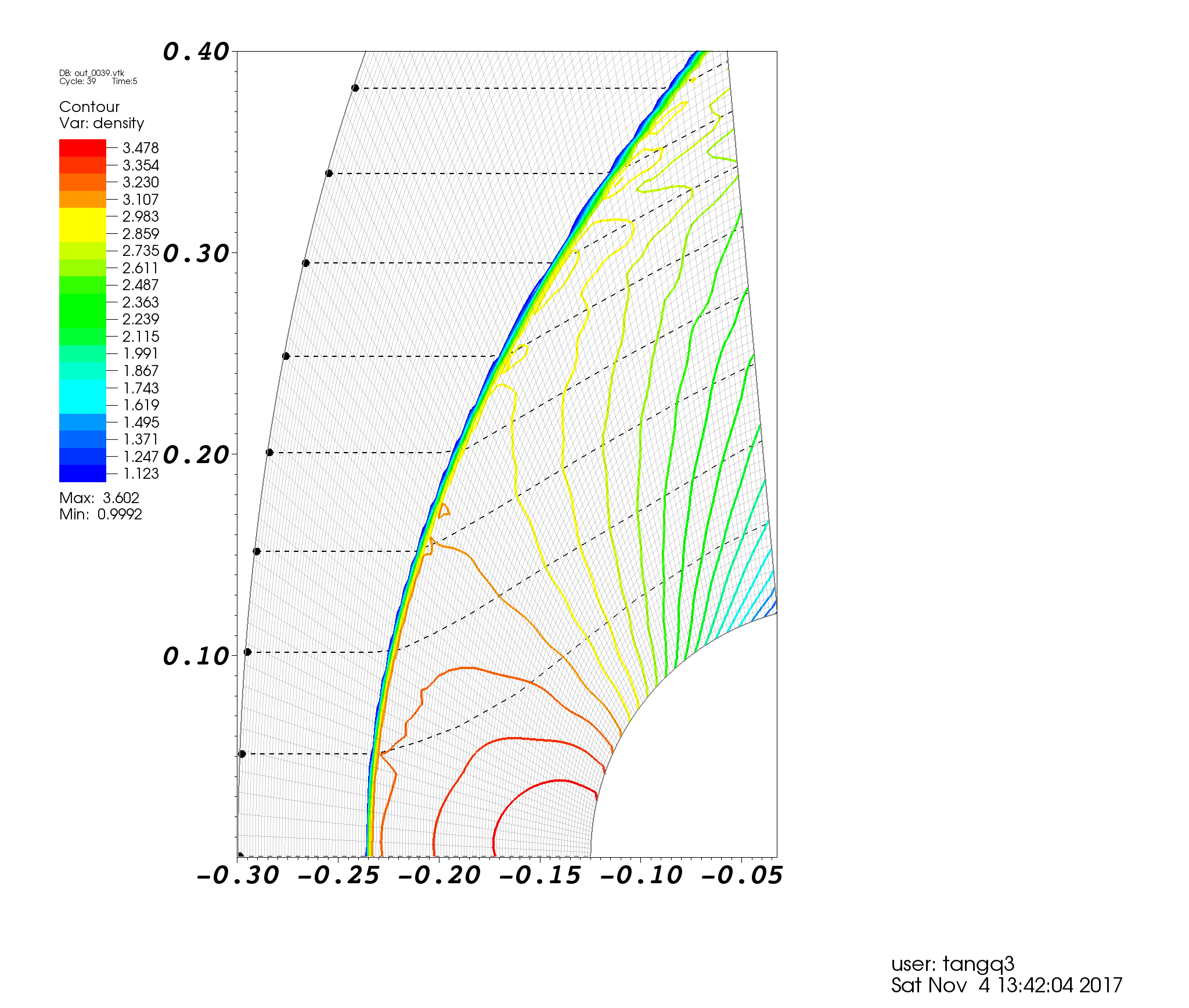}{\figWidth}};
  \draw(  12,-0.8) node[anchor=south west,xshift=0pt,yshift=+0pt] {\trimfig{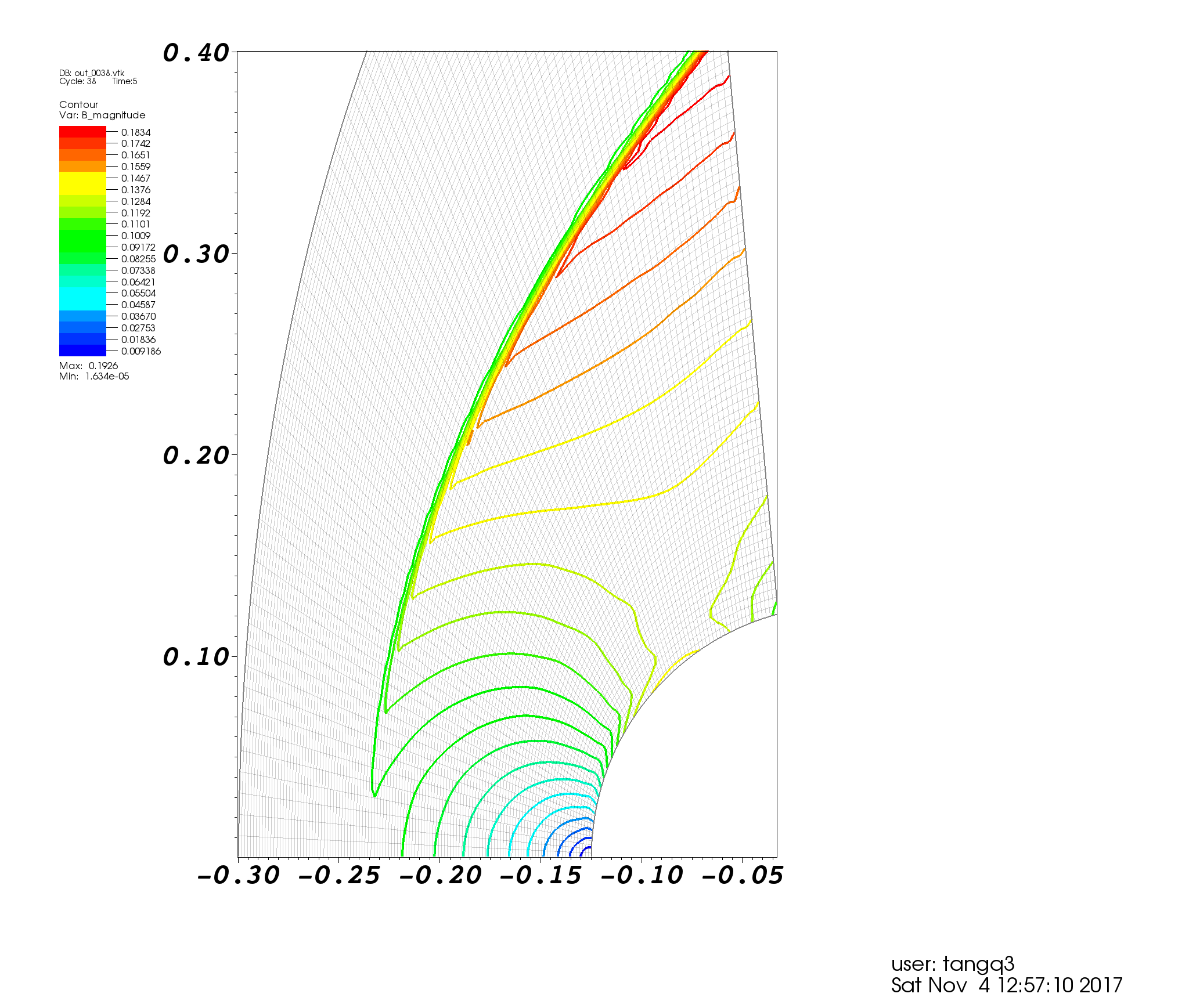}{\figWidth}};
  \draw(.1,4.7) node[draw,fill=white,anchor=west,xshift=2pt,yshift=-1pt] {\small $\rho$};
  \draw(8.6,4.7) node[draw,fill=white,anchor=west,xshift=2pt,yshift=-1pt] {\small $\rho$};
  \draw(4.1,4.7) node[draw,fill=white,anchor=west,xshift=2pt,yshift=-1pt] {\small $\|\Bv\|$};
  \draw(12.6,4.7) node[draw,fill=white,anchor=west,xshift=2pt,yshift=-1pt] {\small $\|\Bv\|$};
%
\end{tikzpicture}
\end{center}
    \caption{Bow shock flow.
Contour plots of the density and magnetic field at $t=5$ are presented with 20 equally spaced contour lines.  
The selected streamlines (the dashed lines) are added to the density plots.
Left two: WENO with the Lax-Friedrichs flux. 
Right two: WENO with the HLLD flux. 
The computational grid of size $120 \times 160$ are also plotted (the light solid lines).
The constrained transport step is turned on.}
\label{fig:compareBowShockFull}
\end{figure}
}

The full benchmark problem is then simulated to obtain a stationary bow shock.
As illustrate in Figure~\ref{fig:compareBoundary}, the divergence errors can 
still lead to spurious oscillations around the bow shock front (see the contours around $y = 0.25$ for instance).
For the results presented in Figure~\ref{fig:compareBowShockFull}, the constrained transport step is turned on and
the WENO schemes with the compatibility boundary condition are used to simulate the problem up to $t = 5$.
The compatibility boundary condition for the potential described in Section~\ref{sec:boundaryConditions} is also used for 
the PEC boundary.
A steady bow shock profile are observed in the both density and magnetic field in Figure~\ref{fig:compareBowShockFull}.
The contour plots are complemented by the plots of several instantaneous streamlines (the dashed lines in the density plot),
similarly to the previous work~\cite{de1999numerical, de2001,nishida2009}.
Our results match well with the results therein and the bow shock is much less smeared compared to their results.
Note that the constrained transport method produces a much smoother profile of the magnetic field.
During the simulations, it is found that the compatibility boundary condition may 
lead to a negative pressure at the ghost points, which does not require any extra care 
(note that the positivity-preserving limiter is not used in this case).
The fundamental reason is that the pressure computed from 
the step~\eqref{eqn:compatibilityBoundary} is only used to convert the normal 
derivative of the pressure to the tangential derivatives along the surface.
Since the pressure at ghost points are only used to constructed the numerical fluxes, the negative pressure at those points
will not cause the solutions at the computational domain become unphysical.
Finally, the HLLD flux is also slightly advantageous over the local Lax-Friedrichs flux in terms of resolving 
the bow shock, but some oscillations are observed in its density plot,
which is not surprising since the HLLD flux is much less dissipative.





%% file: texSISC/conclusions.tex
\section{Conclusions}
\label{sec:conclusions}

In this work, we have extended an alternative flux formulation of the WENO scheme to the ideal MHD equations.
Several Riemann solvers including a HLLD Riemann solver are used to approximate 
the leading (low-order) term in the numerical flux.
The higher-order terms in the numerical flux are approximated by limited central differences of the physical flux and 
the limiter is based on the smoothness indicators in WENO interpolations.  
An unstaggered constrained transport method is
used to control the divergence error of the magnetic field and a positivity-preserving limiter is 
implemented to increase the robustness of the scheme.  
The resulting scheme is applicable to general curvilinear meshes.  
To solve some benchmark problems involving a curved PEC boundary, 
we also derive a numerical compatibility boundary condition for both conserved quantities and magnetic potential.

Several numerical benchmark problems are used to validate the resultant scheme and to confirm the fourth-order accuracy of the scheme.
The results show that when a low dissipative solver such as the HLLD solver is used in the base scheme,
shocks and other complex features are much better captured than
those in the solutions obtained using the base scheme with the Lax-Friedrichs solver. 
Through solving the benchmark problems on different meshes, we further demonstrate  
the robustness of the scheme on general curvilinear meshes.
Using a bow shock flow benchmark, we confirm that the compatibility boundary condition produces
a solution consistent with the solutions obtained by other low-order finite volume methods.
It also demonstrates that a common implementation of the PEC boundary through a reflective boundary condition
may cause failure of high-order methods.

There are several future directions for the current work.
We are interested in extending the current scheme to three dimensions, 
which requires some work particularly in the constrained transport step.
A compatibility boundary condition for a vector potential requires some careful derivations,
in order to simulate a three-dimensional bow shock problem.
In addition, a single-stage single-step scheme based on the alternative flux
formulation of the WENO scheme can be further derived.
Finally, the applications of the current scheme will be further explored in areas involving complex geometry.




%% file: texSISC/WENO.tex
\section{WENO interpolation for a system}
\label{sec:WENO}

For the sake of completeness, a WENO interpolation used in the current work is described for a hyperbolic system.
Note that the WENO interpolation used here is performed on the local characteristic variables instead of on the components of $\qv$.
The resulting approximations to $\qv_{i+1/2}^{\pm}$  are fifth-order accurate.


 \begin{enumerate}
\item Compute an average state
$\qv_{i+1/2}$.  
In the current work, we use the arithmetic mean of primitive variables by setting
$\psi_{i+1/2} = (\psi_{i} + \psi_{i+1})/2$
with $\psi$ in the range of  $\{\rho, \uv, p, \Bv \}$, 
and the conserved variables $\qv = \{ \rho, \rho \uv, \Ec, \Bv \}$ are then recovered from $\psi$.
 	
\item Compute the right and left eigenvectors of the Jacobian $\partial \fv / \partial \qv$,
and denote their matrices by
$$R_{i+1/2}=R(\qv_{i+1/2}) ,\quad R^{-1}_{i+1/2}=R^{-1}(\qv_{i+1/2}).$$

\item Project the conserved quantities $\qv$, which is in the stencil of computing the numerical flux $\fv_{i+1/2}$, to the local characteristic variables $\vv$,
\begin{equation}
\vv_{j} = R^{-1}_{i+1/2} \qv_{j}, \quad \text{for } j = i-2, \ldots, i+3.
\end{equation}
 	
\item Perform a scalar WENO interpolation on each component of the characteristic variable $\vv_{j}$ to obtain the corresponding component of $\vv_{i+1/2}^{\pm}$. 
Here, the procedure of a fifth-order WENO interpolation to obtain the $k$-th component $v^{-}_{k,i+1/2}$ is described:
\begin{enumerate}
	\item Choose one big stencil as $S = \{x_{i-2},\ldots,x_{i+2}\}$, and three small stencils as $S^{(0)} = \{x_{i}, x_{i+1}, x_{i+2}\}$, $S^{(1)} = \{x_{i-1}, x_{i}, x_{i+1}\}$, and $S^{(2)} = \{x_{i-2}, x_{i-1}, x_{i}\}$. 
On those four stencils,  the standard interpolation gives
	\begin{subequations}
        \begin{alignat}{3}
	& v^{(0)}_{k,i+1/2} = \frac{3}{8}v_{k,i} +\frac{3}{4}v_{k,i+1} -\frac{1}{8}v_{k,i+2},\\
	& v^{(1)}_{k,i+1/2} = -\frac{1}{8}v_{k,i-1} +\frac{3}{4}v_{k,i} +\frac{3}{8}v_{k,i+1},\\
	& v^{(2)}_{k,i+1/2} = \frac{3}{8}v_{k,i-2} -\frac{5}{4}v_{k,i-1} +\frac{15}{8}v_{k,i},\\
	& v^\text{big}_{k,i+1/2} = d_{0} v^{(0)}_{k,i+1/2} +d_{1} v^{(1)}_{k,i+1/2} +d_{2} v^{(2)}_{k,i+1/2},
	\end{alignat}
	\end{subequations}
	with the linear weights being $d_{0}={5}/{16}$, $d_{1}={5}/{8}$ and $d_{2}={1}/{16}$.
	\item Compute nonlinear weights $\omega_r$ from the linear weights $d_r$, 
	\begin{equation}
	\omega_{r}=\frac{\alpha_{r}}{\alpha_{0}+\alpha_{1}+\alpha_{2}}, \quad \alpha_{r}=\frac{d_{r}}{(\beta_{r}+\epsilon)^2}, \quad r=0,1,2,
	\end{equation}
	where $\epsilon = 10^{-6}$ is used to avoid division by zero, and the \emph{smoothness indicators} are given by
	\begin{subequations}
	\label{eq:beta}
        \begin{alignat}{3}
	&\beta_{0} = \frac{13}{12}\left( v_{k,i}-2v_{k,i+1}+v_{k,i+2}\right)^2 +\frac{1}{4}\left( 3v_{k,i}-4v_{k,i+1}+v_{k,i+2}\right)^2,\\
	&\beta_{1} = \frac{13}{12}\left( v_{k,i-1}-2v_{k,i}+v_{k,i+1}\right)^2 +\frac{1}{4}\left( v_{k,i}-v_{k,i+1}\right)^2,\\
	&\beta_{2} = \frac{13}{12}\left( v_{k,i-2}-2v_{k,i-1}+v_{k,i}\right)^2 +\frac{1}{4}\left( v_{k,i-2}-4v_{k,i-1}+3v_{k,i}\right)^2.
	\end{alignat}
	\end{subequations}
	 \item The WENO interpolation is defined by
	  $$v^{-}_{k,i+1/2}=\sum_{r=0}^{2}\omega_{r}v^{(r)}_{k,i+1/2}.$$
\end{enumerate}
Note that the process to obtain $\vv^{+}_{i+1/2}$ is mirror-symmetric to the procedure described above.

\item Finally, project $\vv_{i+1/2}^{\pm}$ back to the conserved quantities,
\begin{equation}
\qv_{i+1/2}^{\pm} = R_{i+1/2} \vv_{i+1/2}^{\pm}.
\end{equation}
 \end{enumerate}

%% file: texSISC/riemannSolver.tex
\section{HLL-type Riemann solvers for MHD equations}
  \label{sec:riemannSolver}


Several HLL-type Riemann solvers for ideal MHD equations are described for the Riemann problem given by
the initial  conditions~(\ref{eq:RiemannProblemInitialCondition}).
Since the HLLD solver is an extension of the HLL and HLLC
solvers, these two solvers are first reviewed in~\ref{sec:MHDHLL} and~\ref{sec:MHDHLLC}, 
and our version of the HLLD solver is then introduced in~\ref{sec:MHDHLLD}.
To save the space, the discussion here only focuses on the approximation solutions, 
and the corresponding numerical flux, which can be easily worked out, are therefore not described.
 
  \subsection{The HLL approximate Riemann solver}
  \label{sec:MHDHLL}
  The approximate solution $\tilde{\qv}$ in the HLL solver,
  consisting of three states, is given by
  \begin{equation}
    \tilde{\qv}(t, \xv\cdot\nv) =
    \begin{cases}
      \qv_{\LL},  &\text{if } (\xv\cdot\nv) / t \leq S_{\LL}, \\
      \qv_{\hll}, &\text{if } S_{\LL} \leq (\xv\cdot\nv) / t \leq S_{\RR}, \\
      \qv_{\RR}, &\text{if } S_{\RR} \leq (\xv\cdot\nv) / t,
    \end{cases}
  \end{equation}
  where $S_{\LL}$ and $S_{\RR}$ are the smallest and largest of all the signal speeds. 
  The intermediate state satisfies
  \begin{equation}
    \label{eq:HLL}
    \qv_{\hll} = \frac{S_{\RR}\qv_{\RR} - S_{\LL} \qv_{\LL} + \Fv_{\LL} - \Fv_{\RR}}{S_{\RR} - S_{\LL}},
  \end{equation}
  which is given by the \emph{consistency condition} of conservation laws 
  \begin{equation}
      \label{eq:consistencyCondition}
    \int_{\xi_{\LL}}^{\xi_{\RR}} \tilde{\qv}(t, \xi) \, d\xi = \xi_{\RR} \qv_{\RR} - \xi_{\LL} \qv_{\LL} + T( \Fv_{\LL} - \Fv_{\RR} ),
  \end{equation}
  with $\xi_{\LL} = T S_{\LL}$ and $\xi_{\RR} = T S_{\RR}$.
  $\qv_\text{HLL}$ is used  in the numerical flux of the \emph{subsonic} case where $S_{\LL} \le 0 \le S_{\RR}$.




  \subsection{The HLLC approximate Riemann solver}
  \label{sec:MHDHLLC}

  The approximate solution, consisting of two intermediate states connected by a contact discontinuity, is given by
  \begin{equation}
    \tilde{\qv}(t, \xv\cdot\nv) = 
    \begin{cases}
      \qv_{\LL},  & \text{if } (\xv\cdot\nv)/t \leq S_{\LL}, \\
      \qv_{\LL}^{*},  & \text{if } S_{\LL} \leq (\xv\cdot\nv)/t \leq S_{\MM}, \\
      \qv_{\RR}^{*},  & \text{if } S_{\MM} \leq (\xv\cdot\nv)/t \leq S_{\RR}, \\
      \qv_{\RR}, & \text{if } S_{\RR} \leq (\xv\cdot\nv)/t,
    \end{cases}
  \end{equation}
  where $S_{\MM}$ is the estimated speed of an entropy wave, and $\qv_{\LL}^{*}$ and
  $\qv_{\RR}^{*}$ are the intermediate states. Consider the subsonic case where $S_{\LL} \leq 0 \leq S_{\RR}$. In~\cite{Li2005} it was assumed
  \begin{gather}
    \label{eq:HLLCSM}
 S_{\MM} =    \uv_{\LL}^{*}\cdot\nv = \uv_{\RR}^{*}\cdot\nv = \uv_{\hll}\cdot\nv, \\
    \label{eq:HLLCptotEqual}
    \hLstar{\ptot} = \hRstar{\ptot},
  \end{gather}
since the normal velocities and total pressures are same across a contact discontinuity.
 It was also assumed
  \begin{equation}
    \label{eq:HLLCNormalB}
    \Bv_{\LL}^{*} \cdot \nv = \Bv_{\RR}^{*} \cdot \nv = \Bv_{\hll} \cdot \nv.
  \end{equation}
  Note that the RH condition~\eqref{eq:RHCondition} across
  $\hA{S}$ implies
  \begin{equation}
    \label{eq:RHCondition1}
    \hA{S}\hAstar{\qv} - \hAstar{\Fv} = \hA{S} \hA{\qv} - \hA{\Fv},
  \end{equation}
  where  $\alpha = L$ and  $R$.
  Applying~\eqref{eq:RHCondition1} to $\rho$ gives
    $\hA{S} \hAstar{\rho} - \hAstar{\rho} \hAstar{\uv} \cdot \nv =\hA{S} \hA{\rho} - \hA{\rho} \hA{\uv} \cdot \nv$,
  which, with~\eqref{eq:HLLCSM}, implies
  \begin{equation}
    \label{eq:HLLCrho}
    \hAstar{\rho} = \hA{\rho} \dfrac{\hA{S} - \hA{\uv}\cdot\nv}{\hA{S} - S_{\MM}}.
  \end{equation}
  Applying~\eqref{eq:RHCondition1} to $\rho\uv$ gives
  \begin{equation}
    \label{eq:HLLCRHrhou}
    \hA{S} \hAstar{\rho}\hAstar{\uv} - {\big[ \hAstar{\rho} {\left( \hAstar{\uv} \cdot \nv \right)} \hAstar{\uv} + \hAstar{\ptot}\nv - {\left( \hAstar{\Bv} \cdot \nv \right) \hAstar{\Bv}}\big]}
       = \hA{S} \hA{\rho}\hA{\uv} - {\big[ \hA{\rho} {\left(
              \hA{\uv} \cdot \nv \right)} \hA{\uv} +
          \hA{\ptot}\nv - {\left( \hA{\Bv} \cdot \nv \right) \hA{\Bv}}\big]}.
  \end{equation}
  Taking the normal component of~\eqref{eq:HLLCRHrhou} and using~\eqref{eq:HLLCrho} give
  \begin{equation}
    \label{eq:HLLCptot}
    \hAstar{\ptot} = \hA{\ptot} + \hA{\rho}{\left(\hA{S} - \hA{\uv} \cdot \nv\right)} {\left(S_{\MM} - \hA{\uv} \cdot \nv\right)} + {\left( \hAstar{\Bv} \cdot \nv \right)}^2 - {\left( \hA{\Bv} \cdot \nv \right)}^2.
  \end{equation}
  Using the assumption~\eqref{eq:HLLCSM} to eliminate $\hAstar{\uv} \cdot \nv $ in~\eqref{eq:HLLCRHrhou} gives
  \begin{equation}
    \label{eq:HLLCrhou}
    \hAstar{\rho}\hAstar{\uv} = \dfrac{\hA{\rho}\hA{\uv}(\hA{S} - \hA{\uv}\cdot\nv) + {\left(\hAstar{\ptot} - \hA{\ptot}\right)}\nv + {\left(\hA{\Bv}\cdot\nv\right)}\hA{\Bv} - {\left(\hAstar{\Bv}\cdot\nv\right)}\hAstar{\Bv}  }{\hA{S} - S_{\MM}}
  \end{equation}
  Similarly, applying~\eqref{eq:RHCondition1} to $\Ec$ and using the assumption~\eqref{eq:HLLCSM}  give
  \begin{equation}
    \label{eq:HLLCE}
    \hAstar{\Ec} = \dfrac{ \hA{\Ec} {\left( \hA{S} - \hA{\uv}\cdot\nv \right)} + \hAstar{\ptot} S_{\MM} - \hA{\ptot} {\left( \hA{\uv}\cdot\nv \right)} + {\left( \hA{\Bv} \cdot \nv  \right)}  {\left( \hA{\Bv} \cdot \hA{\uv} \right)} - {\left( \hAstar{\Bv} \cdot \nv  \right)}  {\left( \hAstar{\Bv} \cdot \hAstar{\uv} \right)} }{\hA{S} - S_{\MM}}.
  \end{equation}
  Here $\hAstar{\Bv}$ still needs to be determined.
  Note that the consistency condition in the HLLC solver gives
  \begin{equation}
    \label{eq:HLLCConsistency}
    \dfrac{S_{\MM} - S_{\LL}}{S_{\RR} - S_{\LL}}  \hLstar{\qv} + \dfrac{S_{\RR} - S_{\MM}}{S_{\RR} - S_{\LL}}  \hRstar{\qv} = \qv_{\hll}.
  \end{equation}
  Applying it to $\rho\uv$ and using~\eqref{eq:HLL} and~\eqref{eq:HLLCptotEqual} give
  ${\left( \hLstar{\Bv} \cdot \nv \right)} \hLstar{\Bv} = {\left( \hRstar{\Bv} \cdot \nv \right)} \hRstar{\Bv}$.
  Therefore,~\cite{Li2005} suggested to set
  \begin{equation}
    \label{eq:HLLCBstarEqual}
    \hLstar{\Bv} = \hRstar{\Bv} = \Bv_{\hll}.
  \end{equation}
  Finally, the intermediate states are
  completely determined by~\eqref{eq:HLLCSM},~\eqref{eq:HLLCrho},~\eqref{eq:HLLCptot},
  ~\eqref{eq:HLLCrhou},~\eqref{eq:HLLCE}~and~\eqref{eq:HLLCBstarEqual}.  



  \subsection{The HLLD approximate Riemann solver}
  \label{sec:MHDHLLD}

In the previous version of the HLLD Riemann solver~\cite{Miyoshi2005},
it was assumed that the normal vector to the discontinuity interface
is parallel to one of the coordinate axes and
the magnetic field components normal to the discontinuity interface are identical across the discontinuity.  
Neither of the assumptions holds on a curvilinear mesh.  
Therefore, we consider Riemann problems with  the initial
conditions~(\ref{eq:RiemannProblemInitialCondition}) where the
normal vector $\nv$ to the discontinuity interface may not be parallel to the coordinate axes and there are jumps in the normal magnetic field. 
Note that the divergence-free condition in multiple dimensions can be saved by the jumps in the tangential directions.  
We present our version of the HLLD solver that can handle these issues.
The approximate solution,  consisting of four
intermediate states which are connected by two rotational discontinuities and
one contact or tangential discontinuity, is given by
\begin{equation}
  \label{eq:HLLDApproximateSolution}
  \tilde{\qv}(t, \xv\cdot\nv) = 
  \begin{cases}
    \hL{\qv},  & \text{if } (\xv\cdot\nv)/t \leq \hL{S}, \\
    \hLstar{\qv},  & \text{if } \hL{S} \leq (\xv\cdot\nv)/t \leq \hLstar{S}, \\
    \hLstartwo{\qv},  & \text{if } \hLstar{S} \leq (\xv\cdot\nv)/t \leq S_{\MM}, \\
    \hRstartwo{\qv},  & \text{if } S_{\MM} \leq (\xv\cdot\nv)/t \leq \hRstar{S}, \\
    \hRstar{\qv}, & \text{if } \hRstar{S} \leq (\xv\cdot\nv)/t \leq \hR{S}, \\
    \hR{\qv}, & \text{if } \hR{S} \leq (\xv\cdot\nv)/t,
  \end{cases}
\end{equation}
where  $\hLstar{S}$ and $\hRstar{S}$ are the estimated speeds of rotational
discontinuities, and $S_{\MM}$ is the estimated speed of a contact
or tangential discontinuity.  


  We first consider the subsonic case.  Following the
  assumptions~\eqref{eq:HLLCSM}--\eqref{eq:HLLCNormalB}, we similarly assume
  \begin{gather}
    \label{eq:HLLDSM}
 S_{\MM} = \hLstar{\uv}\cdot\nv = \hLstartwo{\uv}\cdot\nv = \hRstartwo{\uv}\cdot\nv = \hRstar{\uv}\cdot\nv = \uv_{\hll}\cdot\nv, \\
    \label{eq:HLLDptotEqual}
    \hLstar{\ptot} = \hLstartwo{\ptot} = \hRstartwo{\ptot} = \hRstar{\ptot}, \\
    \label{eq:HLLDNormalB}
    \hLstar{\Bv}\cdot\nv = \hLstartwo{\Bv}\cdot\nv =
    \hRstartwo{\Bv}\cdot\nv = \hRstar{\Bv}\cdot\nv =
    \Bv_{\hll}\cdot\nv.
  \end{gather}
  Recall that rotational discontinuities are linearly degenerate and
  correspond to the \Alfven~waves of speed
  $\uv\cdot\nv \mp \sqrt{{\left( \Bv \cdot \nv \right)}^2
    / \rho}$.  Since densities are identical across rotational discontinuities, their speeds are 
  \begin{gather}
    \label{eq:HLLDSstar}
    \hAstar{S} = S_{\MM} \mp \sqrt{{{\left( \Bv_{\hll} \cdot \nv \right)}^2}/{\hAstar{\rho}}},
  \end{gather}
  where $-$ and $+$ correspond to $\alpha=\LL$ and $\RR$ respectively.
  Note that the four-intermediate-state solution~\eqref{eq:HLLDApproximateSolution} degenerates
  to a two-intermediate-state solution when either $S^*_\alpha$ is close to $S_\alpha$ for $\alpha = L$ and $R$, 
  or $\hLstar{S}$ and $\hRstar{S}$ are both close to $S_{\MM}$,  
  in which case the HLLD solver degenerates to the HLLC solver and no extra cares are needed.
  Note $\rho_\alpha^*$ can be determined similarly to~\eqref{eq:HLLCrho} and remains the same as
  \begin{equation}
    \label{eq:HLLDrho}
    \hAstar{\rho} = \hA{\rho} \dfrac{\hA{S} - \hA{\uv}\cdot\nv}{\hA{S} - S_{\MM}}.
  \end{equation}


  Now consider the non-degenerate case when all the speeds are sufficiently spread apart. 
  We first obtain $\hAstar{\qv}$ from the RH condition
  across $\hA{S}$.  
  \eqref{eq:HLLCptot} and~\eqref{eq:HLLCrhou} also remain the same as 
  \begin{align}
    \label{eq:HLLDptot}
    \hAstar{\ptot} & = \hA{\ptot} + \hA{\rho}{\left(\hA{S} - \hA{\uv} \cdot \nv\right)} {\left(S_{\MM} - \hA{\uv} \cdot \nv\right)} + {\left( \hAstar{\Bv} \cdot \nv \right)}^2 - {\left( \hA{\Bv} \cdot \nv \right)}^2, \\
    \label{eq:HLLDrhoustar}
    \hAstar{\rho}\hAstar{\uv} & = \dfrac{\hA{\rho}(\hA{S} -
      \hA{\uv}\cdot\nv) + {\left(\hAstar{\ptot} -
          \hA{\ptot}\right)}\nv +
      {\left(\hA{\Bv}\cdot\nv\right)}\hA{\Bv} -
      {\left(\hAstar{\Bv}\cdot\nv\right)}\hAstar{\Bv}
    }{\hA{S} - S_{\MM}}.
  \end{align}
  Applying the RH condition across $\hA{S}$ to $\Bv$ gives
  \begin{equation}
    \label{eq:HLLDRHB}
\hA{S}\hAstar{\Bv} - \hAstar{\Big[ {\left( {\uv}\cdot\nv \right)} \Bv - {\left( \Bv\cdot\nv \right)} \uv \Big]}
    =
    \hA{S}\hA{\Bv} - \hA{\Big[ {\left( \uv\cdot\nv \right)} \Bv - {\left( \Bv\cdot\nv \right)} \uv \Big]}.
  \end{equation}
  Here the notation $[\cdot]_\alpha^*$ stands for the flux of the corresponding solution state.
  The equations~\eqref{eq:HLLDrhoustar}
  and~\eqref{eq:HLLDRHB} form a linear system 
  for $\hAstar{\uv}$ and $\hAstar{\Bv}$, and the solutions are
  \begin{align}
    \label{eq:HLLDustar}
  &  \begin{aligned}
      \hAstar{\uv}
      & =\Big[ \Big({ {\left( \hA{\Bv}\cdot\nv \right)} {\left(
            \hA{S} - S_{\MM} \right)} - {\left(\Bv_{\hll}\cdot\nv \right)} {\left( \hA{S} -
            \hA{\uv}\cdot\nv \right)} }\Big)\hA{\Bv} \\
      & + \Big({\hA{\rho} {\left( \hA{S} - \hA{\uv}\cdot\nv  \right)} {\left( \hA{S} - S_{\MM}  \right)} - {\left( \Bv_{\hll} \cdot \nv  \right)} {\left( \hA{\Bv} \cdot \nv  \right)}}\Big) \hA{\uv}
  + { {\left( \hAstar{\ptot} - \hA{\ptot}  \right)}  {\left( \hA{S} - S_{\MM}  \right)}  } \nv \Big]/D,
    \end{aligned} \\
    \label{eq:HLLDBstar2}
  &  \begin{aligned}
      \hAstar{\Bv}
      & = \Big[ \Big({\hA{\rho} {\left( \hA{S} - \hA{\uv}\cdot\nv  \right)}^2 - {\left( \Bv_{\hll} \cdot \nv  \right)} {\left( \hA{\Bv}\cdot\nv  \right)}  }\Big) \hA{\Bv} \\
      & + \Big({ \hA{\rho} {\left( \hA{S} - \hA{\uv}\cdot\nv \right)} {\left( \hA{\Bv}\cdot\nv - \Bv_{\hll}\cdot\nv  \right)} }\Big) \hA{\uv}  - { {\left( \hAstar{\ptot} - \hA{\ptot} \right)} {\left(
      \Bv_{\hll}\cdot\nv \right)} } \nv\Big]/D,
    \end{aligned}
  \end{align}
  with $D = {\hA{\rho} {\left( \hA{S} - \hA{\uv}\cdot\nv  \right)} {\left( \hA{S} - S_{\MM}  \right)} - {\left( \Bv_{\hll} \cdot \nv  \right)}^2}$. Applying the RH condition across $\hA{S}$ to $\Ec$ gives
  \begin{equation}
    \label{eq:HLLDE}
    \hAstar{\Ec} =
    \dfrac{ \hA{\Ec} {\left( \hA{S} - \hA{\uv}\cdot\nv \right)} + \hAstar{\ptot} S_{\MM} - \hA{\ptot} {\left( \hA{\uv}\cdot\nv \right)} + {\left( \hA{\Bv} \cdot \nv  \right)}  {\left( \hA{\Bv} \cdot \hA{\uv} \right)} - {\left( \hAstar{\Bv} \cdot \nv  \right)}  {\left( \hAstar{\Bv} \cdot \hAstar{\uv} \right)} }{\hA{S} - S_{\MM}},
  \end{equation}
  which has the same form as~\eqref{eq:HLLCE}.  Note that whereas
  the consistency condition~\eqref{eq:HLLCConsistency} on $\rho\uv$ is used to obtain $\hAstar{\Bv}$
  in the HLLC solver, here the RH condition across $\hA{S}$ is applied
  to $\rho\uv$ and $\Bv$ for determining $\hAstar{\rho\uv}$ and
  $\hAstar{\Bv}$, and a consistency condition similar
  to~\eqref{eq:HLLCConsistency}  is saved for determining $\hAstartwo{\left(\rho\uv\right)}$ and $\hAstartwo{\Bv}$ later.

  Next consider the inner intermediate states $\hAstartwo{\qv}$. Applying the RH condition across $\hAstar{S}$ to $\rho$ gives
    $\hAstar{S} \hAstartwo{\rho} - \hAstartwo{\left( \rho \uv \right)} \cdot \nv
    =
    \hAstar{S} \hAstar{\rho} - \hAstar{\left( \rho \uv \right)} \cdot \nv$,
  which, with~\eqref{eq:HLLDSM}, gives
  \begin{equation}
    \label{eq:HLLDrhostartwo}
    \hAstartwo{\rho} = \hAstar{\rho}.
  \end{equation}
  Note that this is consistent with the fact that densities are identical across 
  rotational discontinuities, which is assumed in deriving the estimates of 
  $\hLstar{S}$ and $\hRstar{S}$ in~\eqref{eq:HLLDSstar}.
  Since $\hLstar{S}$ and $\hRstar{S}$ are sufficiently spread apart, we assume $\Bv_{\hll}\cdot\nv \neq 0$.
Therefore, applying the RH condition across $S_{\MM}$ to $\rho\uv$ implies
    $\hLstartwo{\Bv} = \hRstartwo{\Bv}$
  and applying it to $\Bv$ implies
    $\hLstartwo{\uv} = \hRstartwo{\uv}$.
  Note that the consistency condition here becomes
  \begin{equation}
    \label{eq:HLLDConsistency1}
    {\left( \hR{S} - \hRstar{S} \right)} \hRstar{\qv} +
    {\left( \hRstar{S} - S_{\MM} \right)} \hRstartwo{\qv} +
    {\left( S_{\MM} - \hLstar{S} \right)} \hLstartwo{\qv} +
    {\left( \hLstar{S} - \hL{S} \right)} \hLstar{\qv}
    - \hR{S}\hR{\qv} + \hL{S}\hL{\qv}
    + \hR{\Fv} - \hL{\Fv} = 0,
  \end{equation}
  and the RH condition across $\hAstar{S}$ becomes
  \begin{equation}
    \label{eq:HLLDRH}
    \hA{S} {\left( \hAstar{\qv} - \hA{\qv}  \right)} = \hAstar{\Fv} - \hA{\Fv}
  \end{equation}
  Substituting~\eqref{eq:HLLDSstar} and~\eqref{eq:HLLDRH}
  into~\eqref{eq:HLLDConsistency1} gives
  \begin{equation}
    \label{eq:HLLDConsistency3}
    \abs{\Bv_{\hll}\cdot\nv} {\left( \dfrac{\hRstartwo{\qv}}{\sqrt{\hRstar{\rho}}} + \dfrac{\hLstartwo{\qv}}{\sqrt{\hLstar{\rho}}} \right) }
    + \hRstar{\Fv} - \hLstar{\Fv}
    - \hRstar{S}\hRstar{\qv} + \hLstar{S}\hLstar{\qv} = 0.
  \end{equation}
  Applying~\eqref{eq:HLLDConsistency3} to $\rho\uv$
  and using~\eqref{eq:HLLDrhostartwo} and~$\hLstartwo{\uv} = \hRstartwo{\uv}$ give 
  \begin{equation*}
    \begin{aligned}
        \phantom{+}\abs{\Bv_{\hll}\cdot\nv} & {\left( \sqrt{\hRstar{\rho}} + \sqrt{\hLstar{\rho}} \right) } \hAstartwo{\uv}  \\ 
        + & \hRstar{\Big[  {\left( \uv\cdot\nv \right)} \rho\uv + \ptot \nv - {\left( \Bv \cdot \nv  \right)} \Bv \Big]}  - \hLstar{\Big[ {\left( \uv\cdot\nv \right)} \rho\uv
          + \ptot \nv - {\left( \Bv \cdot \nv \right)}\Bv
      \Big]} - \hRstar{S} \hRstar{\left( \rho\uv \right)} +
      \hLstar{S} \hLstar{\left( \rho\uv \right)} = 0,
    \end{aligned}
  \end{equation*}
  which, with~\eqref{eq:HLLDSM}--\eqref{eq:HLLDSstar}, implies
  \begin{equation}
    \label{eq:HLLDustartwo}
    \hAstartwo{\uv} = 
    \dfrac{ \sqrt{\hLstar{\rho}} \hLstar{\uv}  + \sqrt{\hRstar{\rho}} \hRstar{\uv} + \sign{\left( \Bv_{\hll}\cdot\nv \right)} {\left( \hRstar{\Bv} - \hLstar{\Bv} \right)} }{ \sqrt{\hLstar{\rho}} + \sqrt{\hRstar{\rho}} }.
  \end{equation}
  Applying~\eqref{eq:HLLDConsistency3} to $\Bv$ and using $\hLstartwo{\Bv} = \hRstartwo{\Bv}$ give
  \begin{equation*}
      \phantom{+} \abs{\Bv_{\hll}\cdot\nv} {\left( \dfrac{1}{\sqrt{\hRstar{\rho}}} + \dfrac{1}{\sqrt{\hLstar{\rho}}} \right) } \hAstartwo{\Bv}
       + \hRstar{\Big[ {\left( \uv\cdot\nv \right)} \Bv - {\left( \Bv\cdot\nv \right)} \uv \Big]} 
      - \hLstar{\Big[ {\left( \uv\cdot\nv \right)} \Bv - {\left( \Bv\cdot\nv \right)} \uv \Big]} -
      \hRstar{S}\hRstar{\Bv} + \hLstar{S}\hLstar{\Bv} = 0,
  \end{equation*}
  which, with~\eqref{eq:HLLDSM},
  {}\eqref{eq:HLLDNormalB} and~\eqref{eq:HLLDSstar} gives
  \begin{equation}
    \label{eq:HLLDBstartwo}
    \hAstartwo{\Bv} = 
    \dfrac{\sqrt{\hRstar{\rho}}\hLstar{\Bv} + \sqrt{\hLstar{\rho}}\hRstar{\Bv} + \sign{\left( \Bv_{\hll}\cdot\nv  \right)} \sqrt{\hLstar{\rho}\hRstar{\rho}}{\left( \hRstar{\uv} - \hLstar{\uv}  \right)} }{\sqrt{\hLstar{\rho}} + \sqrt{\hRstar{\rho}}}.
  \end{equation}
  The RH condition across $\hAstar{S}$ on $\Ec$ gives
  \begin{equation*}
      \hAstar{S}\hAstartwo{\Ec} - \hAstartwo{\Big[ (\Ec+\ptot)(\uv\cdot\nv) - (\uv\cdot\Bv) (\Bv\cdot\nv) \Big]} =
      \hAstar{S}\hAstar{\Ec} - \hAstar{\Big[ (\Ec+\ptot)(\uv\cdot\nv) - (\uv\cdot\Bv) (\Bv\cdot\nv) \Big]},
  \end{equation*}
  which, with~\eqref{eq:HLLDSM}--\eqref{eq:HLLDSstar},
  gives
  \begin{equation}
    \label{eq:HLLDEstartwo}
    \hAstartwo{\Ec} = \hAstar{\Ec} \mp \sqrt{\hAstar{\rho}} \left( \hAstar{\uv}\cdot\hAstar{\Bv} - \hAstartwo{\uv}\cdot\hAstartwo{\Bv} \right) \sign(\Bv_{\hll} \cdot \nv),
  \end{equation}
  where $-$ and $+$ correspond to $\alpha=\LL$ and $\RR$ respectively.

  Finally the intermediate states in the approximate solution are
  completely determined in~\eqref{eq:HLLDSM}--\eqref{eq:HLLDptot}, {}\eqref{eq:HLLDustar}--\eqref{eq:HLLDrhostartwo}, {}\eqref{eq:HLLDustartwo}--\eqref{eq:HLLDEstartwo}.
